\documentclass[aps,physrev,preprint,groupedaddress]{revtex4-2}

\usepackage{amsmath, amssymb}
\usepackage{tikz} 
\newlength{\figscale} 
\usepackage{graphicx}
\usepackage{color}
\usepackage{hyperref} 
\usepackage[capitalize,noabbrev]{cleveref} 
\usepackage[version = 4]{mhchem} 

\DeclareRobustCommand{\healthyfiringicon}{%
  \tikz[baseline=-0.5ex]{\pic[scale=0.8] at (0,0) {healthyfiring};}%
}

\DeclareRobustCommand{\healthytiredicon}{%
  \tikz[baseline=-0.5ex]{\pic[scale=0.8] at (0,0) {healthytired};}%
}

\newtheorem{Lemma}{Lemma}[section]

\newtheorem{Remark}[Lemma]{Remark}

 {\begin{trivlist} \item[]{\textbf{#1.} }}%
{\hspace*{\fill}$\rule{.4\baselineskip}{.4\baselineskip}$\end{trivlist}}

\newcommand{\vb}{\mathbf{v}}
\newcommand{\wb}{\mathbf{w}}
\newcommand{\ab}{\mathbf{a}}
\newcommand{\hide}[1]{}
\newcommand{\red}{\color{red}}         
     
\definecolor{maroon}{rgb}{0.5,0,0.1}
\definecolor{dark green}{rgb}{0, .5, 0}
\definecolor{pink}{rgb}{1, 0.08, 0.58}
\newcommand{\qw}[1]{\color{maroon}#1}

\newcommand{\gw}[1]{\color{dark green}#1}
\newcommand{\ck}[1]{\color{pink}#1}

\tikzset{
    healthy/.style={
        draw=red,
        fill=red!10,
    },
    fast/.style={
        draw=green,
        fill=green!10,
    },
    slow/.style={
        draw=blue,
        fill=blue!10,
    },
}

\tikzset{
    healthy/.pic={\draw[healthy,fill=none] (0,0) circle[x radius=0.25,y radius=0.15];},
    healthyfiring/.pic={\draw[red,fill=red!50] (0,0) circle [x radius=0.25,y radius=0.15];},
    healthytired/.pic={{\draw[healthy,fill=none] (0,0) circle[x radius=0.25,y radius=0.15];\path[fill=red!40] (0,0) circle (0.1);}},
    fast/.pic={\draw[fast,fill=none] (0,0) circle [x radius=0.25,y radius=0.15];},
    fastfiring/.pic={\draw[fast,fill=green!50] (0,0) circle[x radius=0.25,y radius=0.15];},
    fasttired/.pic={{\draw[fast,fill=none] (0,0) circle[x radius=0.25,y radius=0.15];\path[fill=green!40] (0,0) circle (0.1);}},
    slow/.pic={\draw[slow,fill=none] (0,0) circle[x radius=0.25,y radius=0.15];},
    slowfiring/.pic={\draw[slow,fill=blue!50] (0,0) circle[x radius=0.25,y radius=0.15];},
    slowtired/.pic={{\draw[slow,fill=none] (0,0) circle[x radius=0.25,y radius=0.15];\path[fill=blue!40] (0,0) circle (0.1);}},
}
\DeclareRobustCommand{\healthyicon}{%
  \tikz[baseline=-0.5ex]{\pic[scale=0.8] at (0,0) {healthy};}%
}

\begin{document}


\title{A Heterogeneous FitzHugh-Nagumo Model Exhibiting Local and Regional Dynamics that  Mimic Spontaneous Atrial Fibrillation}



\author{Alexander Hattoum M.D.}
\email[]{ahattoummd@gmail.com}
\affiliation{West Virginia University Hospitals}
\author{Camden Kilton}
\email[]{ck003122@ohio.edu}
\affiliation{Ohio University}
\author{Martin J. Mohlenkamp}
\email[]{mohlenka@ohio.edu}
\affiliation{Ohio University}

\author{Graham Walther}
\email[]{gw514022@ohio.edu}
\affiliation{Ohio University}
\author{Qiliang Wu}
\email[]{wuq@ohio.edu}
\affiliation{Ohio University}
\author{Todd R. Young}
\email[]{youngt@ohio.edu}
\affiliation{Ohio University}


\date{\today}

\begin{abstract}
We study two-dimensional networks of FitzHugh-Nagumo oscillators in regimes that mimic cardiac tissue. The networks include small layers of cells that have fast-recovery characteristics and adjacent layers of cells that have slow-recovery. Numerical simulations show that the boundary between the fast- and slow-recovery cells can initiate fibrillatory-like behavior when activated by a single, healthy wave in the network. We discuss possible connections between these structures and the onset of atrial fibrillation, with the goal of improving our understanding and treatment of this costly condition. 

Our main conclusion is that a single, normal excitation of small areas of heterogeneity embedded in healthy tissue can  initiate arrhythmias and the two dimensional network can then sustain irregular firing patterns in the form of double swirling-back curves and spiral waves. This occurs with a very simple model of excitable cells and we  give a clear mechanistic  description of the arrhythmia. 
The model illustrates how increasing heterogeneity in heart tissue can progressively cause first intermittent and then persistent arrhythmic behavior. 
Another conclusion is that the center of problematic firings can drift away from the heterogeneity into healthy tissue, thus presenting phantom targets for treatment.

\end{abstract}


\maketitle

\section{Introduction}
\label{sec:introduction}

When a heart cell (myocyte) is triggered, the voltage across its cell membrane follows a pattern called a (single-cell) cardiac action potential.
First (phase 0 and 1) it fires with a rapid spike in its voltage (called a depolarization), sufficient to trigger the next cell, and then its voltage decays through a plateau (phase 2) and recovery (phase 3) during which it is inhibited from firing again, with the inhibition weakening at the end of recovery as the voltage returns to its starting value (phase 4).
In the atrium of a healthy human heart, a beat starts when a voltage spike produced by the sinus (sinoatrial/SA) node triggers the neighboring cells.
A wave of firing then proceeds through the atrium, followed closely by a wave of recovery that ensures that the stimulation from the sinus node causes only a single beat. 
Atrial fibrillation is a heart condition in which the firing wave of a beat becomes disordered, with cells in different regions firing multiple times and out of the proper sequence.
Clinically, atrial fibrillation has a massive health care burden including increased stroke risk, 
hospitalization, and decline in overall cardiac function. 
In this paper, we consider a FitzHugh-Nagumo model of excitable media, with a two-dimensional arrays of cells, each with parameters chosen to model heart cells in the atrium.
Our goal is to produce a model that matches the phenomena observed during atrial fibrillation, with the hope of improving understanding of these phenomena, and eventually improving detection and treatment.

The electrical behavior of a cell can be changed (remodeled) when biological effects modify its ion channels. 
We call cells with a longer plateau (elevated voltage) phase and a longer total recovery phase {\em slow-recovery cells}, and cells with faster voltage decay and a shorter recovery phase {\em fast-recovery cells}.
It is well-established that slow-recovery cells can develop due to organic free radicals (also called reactive oxygen species) \cite{P-C-F-B-W-K:2018,O-H-K-L-L-W-C-K:2007}, which can be created by diabetes \cite{YA-AT-SA:2019}, smoking \cite{CA-TO-BE:2021}, and inflammation in general \cite{CA-TO-BE:2021}; due to the pressure of hypertension \cite{C-M-K-K-C-K-G-M-B:1986}; and due to the aldosterone that causes some forms of hypertension \cite{L-B-Z-Z-M-L-X:2015,L-W-Z-Z:2018}.
It is also well-established that fast-recovery cells can develop due to the parasympathetic action of the vagus nerve \cite{D-G-W-H-H-Detal:2001,A-U-V-N-H-Tetal:2008,G-P-W-G-etal:2019} as a counter-response to inflammation \cite{G-D-D-T-V:2018,JOH-WEB:2009} or hypertension \cite{LAT-SAL:2018}; (temporarily) due to direct physical strain on the cell membrane from acute spikes in blood pressure \cite{NIN-SAI:2008,FRA-BOD:2011} as often caused by episodes of sleep apnea; (temporarily) due to higher-frequency stimulation of the cell (activating the cell's restitutional property) \cite{NOL-DAH:1968}; and (chronically) due to prolonged higher-frequency stimulation of the cell caused by atrial fibrillation itself \cite{D-F-V-J-W-C-K-R:2005}.
Both types of changes occur heterogeneously \cite{C-M-K-K-C-K-G-M-B:1986,RUB-ZIP:2005,D-G-W-H-H-Detal:2001}, especially due to the heterogeneous distribution of vagus nerve fibers \cite{A-N-A-M:1958,C-T-M-C-R-Letal:2005,A-U-V-N-H-Tetal:2008,G-P-W-G-etal:2019,L-C-W-Z-Z:2012}, so one might have regions or layers of slow-recovery cells and regions or layers of fast-recovery cells intermixed among healthy cells. 
Sharp boundaries between cell types can exist due to inherent anatomical differences between regions or by dynamic effects \cite{B-S-M-M-W-H-B-B:1980,C-M-K-K-C-K-G-M-B:1986,ANT-BUR:2011}.
For example, free radicals could create slow-recovery cells and inflammation in a region, then the vagus nerve react to the inflammation and create fast-recovery cells but only along the nerves, leaving a strip of fast-recovery cells.
It is believed that heterogeneity is a key to atrial fibrillation \cite{MO-RH-AB:1964,B-S-M-M-W-H-B-B:1980,ANT-BUR:2011,A-A-K-Z-M-Betal:2019,H-S-V-H:2021,R-R-B-B-T:2023,LI-QU-HU:2023}, since uniformly changed cells would generally still produce a wave that is coherent, even if it is abnormal in other ways.

While it has markedly varied phenomenological displays from patient to patient, tissue cluster to tissue cluster, and even event to event within the same heart, atrial fibrillation is characterized by disorganized firings and recoveries that persist for a significant period of time in a self-sustaining manner, that are irregular, but that often include certain geometric patterns such as double swirling-back curves and wandering spirals \cite{S-C-B-Z-B-Wetal:2020}.
Other breakdowns in normal rhythm (arrhythmias), such as focal or reentrant-circuit based rhythms, have organized waves and clearly defined arrhythmia mechanisms, which enable effective treatment \cite{SASYNI:1984,SCH-YAN:2005,N-K-S-C-R-M:2012}.
For atrial fibrillation, antiarrhythmic medications and procedural interventions are only partially effective \cite{AMU-CUR:2021}. 
The most successful technique is pulmonary vein isolation via ablation (killing selected cells), which blocks firing waves in that region and can terminate atrial fibrillation in around 70\% of people with intermittent (paroxymal) fibrillation \cite{K-B-F-Metal:2016,C-H-C-Ketal:2018,AMU-CUR:2021}, outperforming systemic antiarrhythmic drug therapies \cite{AMU-CUR:2021}.
Identifying ablation targets beyond pulmonary vein isolation is an area of active research \cite{AMU-CUR:2021,R-B-B-A-T-R:2021,R-R-B-B-T:2023}.
Mapping spiral wave centers and targeting them with ablation showed great initial promise in terminating instances of atrial fibrillation \cite{N-K-S-C-R-M:2012}, but has not shown statistically significant long-term benefits \cite{B-V-H-S-R-S-S-S:2025}. 
The target for ablation or other treatments is still uncertain, since the fundamental mechanisms causing atrial fibrillation are still not well-understood.

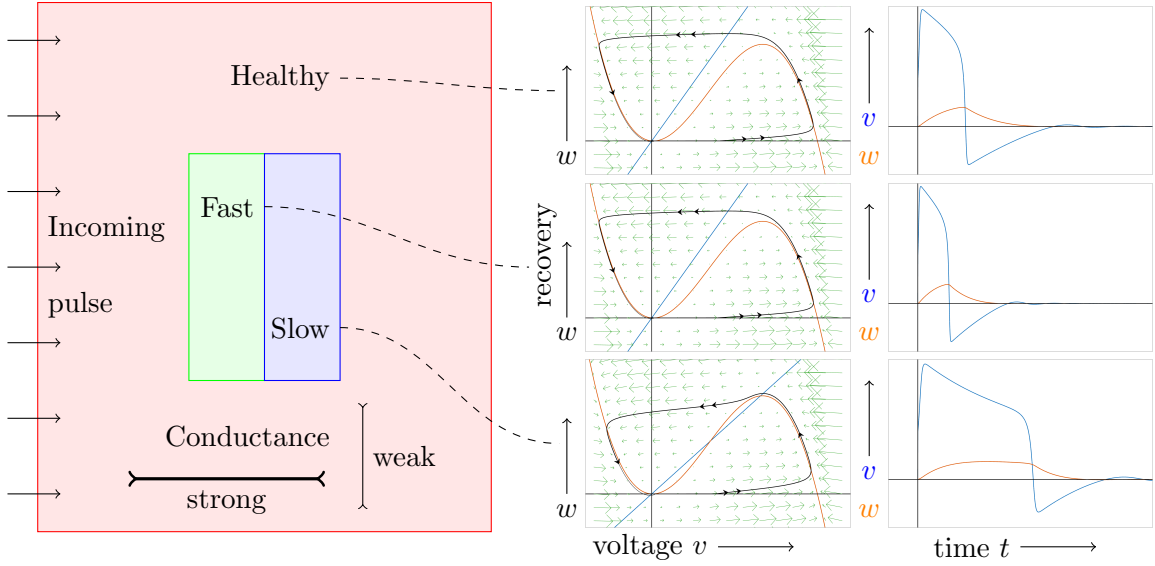
\begin{figure}
    \centering
    \begin{tikzpicture}
        \path[draw,healthy] (0,0) rectangle (6,7); 
        \path[draw,fast] (2,2) rectangle (3,5); 
        \path[draw,slow] (3,2) rectangle (4,5); 
        \node[anchor=south east] at (4,1) {Conductance};
        \path[draw,very thick,>-<] (1.2,0.7) -- (2.5,0.7) node[anchor=north]{strong} -- (3.8,0.7);
        \path[draw,>-<,thin] (4.3,0.3) -- (4.3,1) node[anchor=west]{weak} -- (4.3,1.7);
        \foreach \num in {0.5,1.5,2.5,3.5,4.5,5.5,6.5} 
            \path[draw,->] (-0.4,\num) -- (0.3,\num);
        \node[anchor=west] at (0,4) {Incoming};
        \node[anchor=west] at (0,3) {pulse};
        \path[draw,dashed] (4,6) node[anchor=east]{Healthy} to [out=0,in=180] (6.9,35/6);
        \path[draw,dashed] (3,4.3) node[anchor=east]{Fast} to [out=0,in=180] (6.5,7/2);
        \path[draw,dashed] (4,2.7) node[anchor=east]{Slow} to [out=0,in=180] (6.9,7/6);
        \setlength{\figscale}{3.5cm}
        \def\x {9};
        \def\xx {13};
        \node at (\x,35/6) {\includegraphics[width=\figscale]{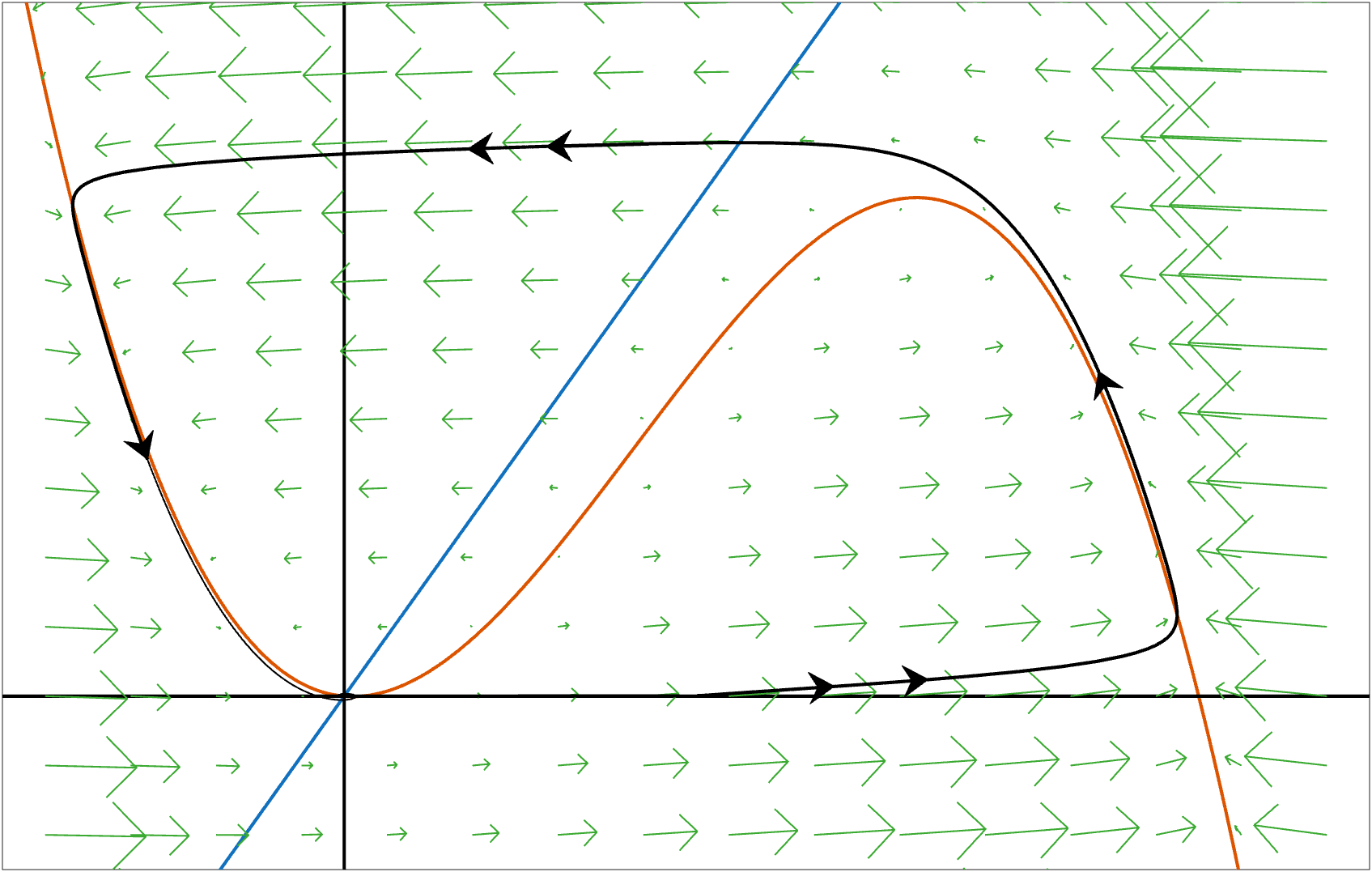}};
        \node at (\xx,35/6) {\includegraphics[width=\figscale]{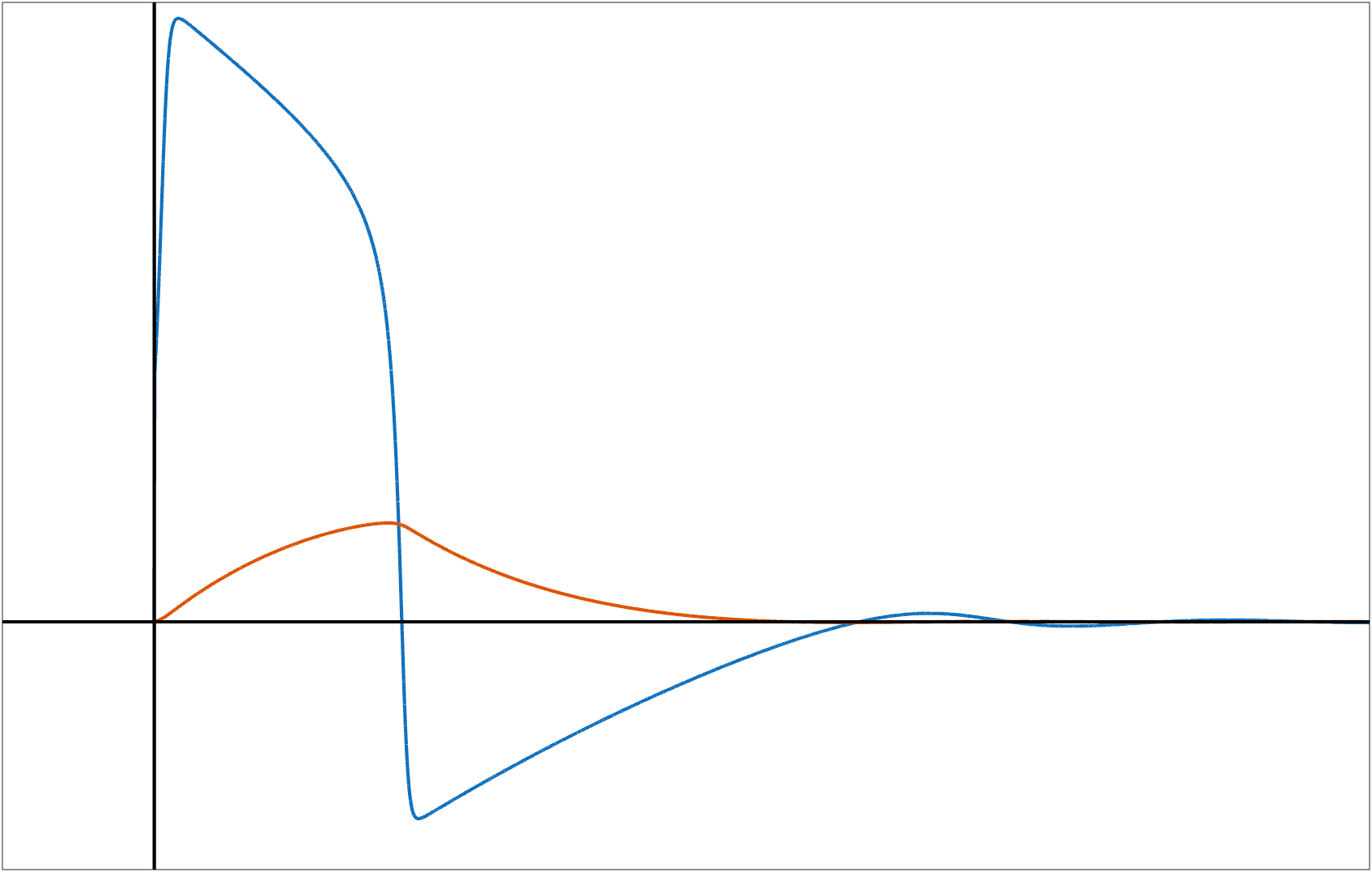}};
        \node at (\x,7/2) {\includegraphics[width=\figscale]{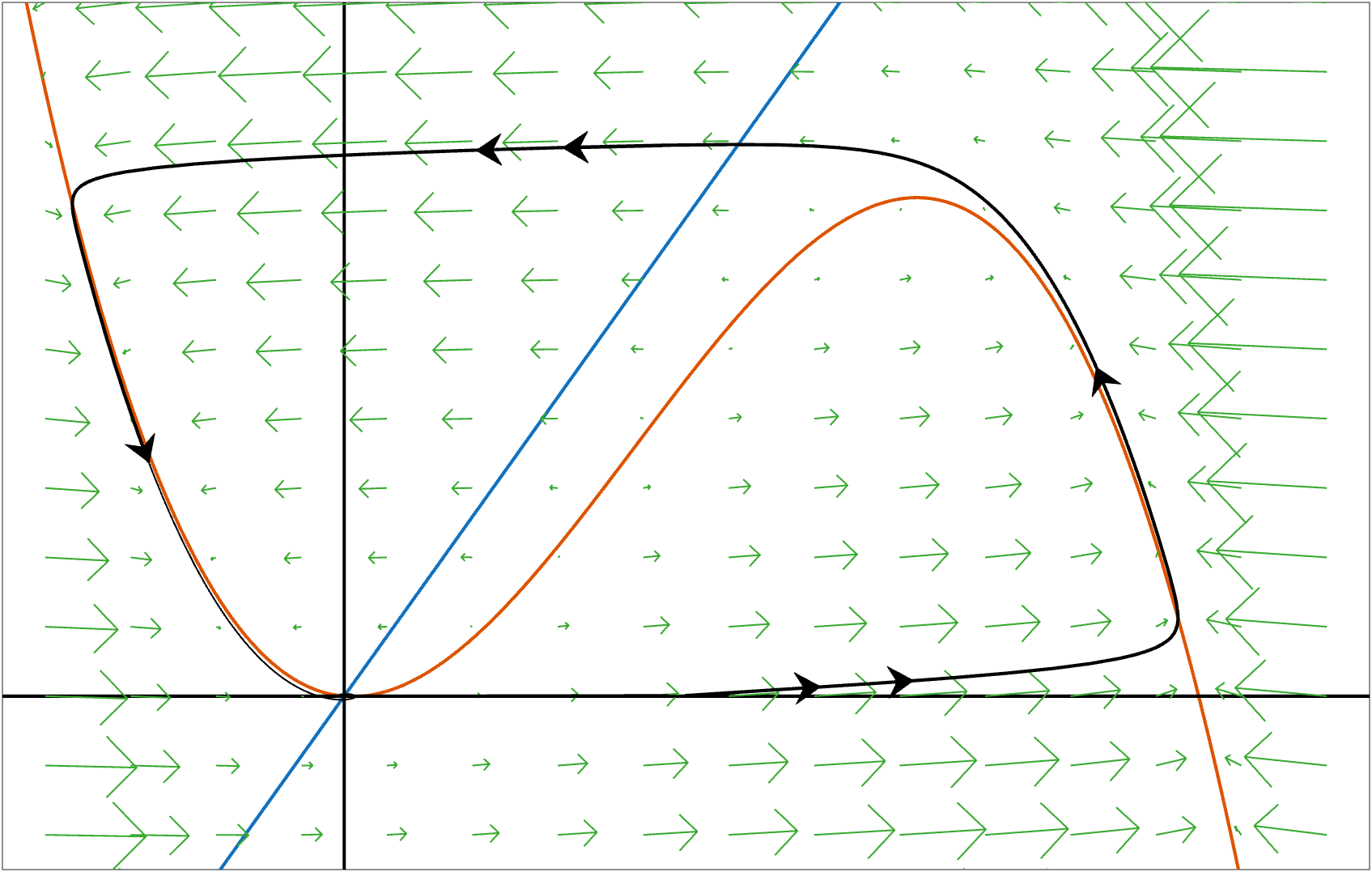}};
        \node at (\xx,7/2) {\includegraphics[width=\figscale]{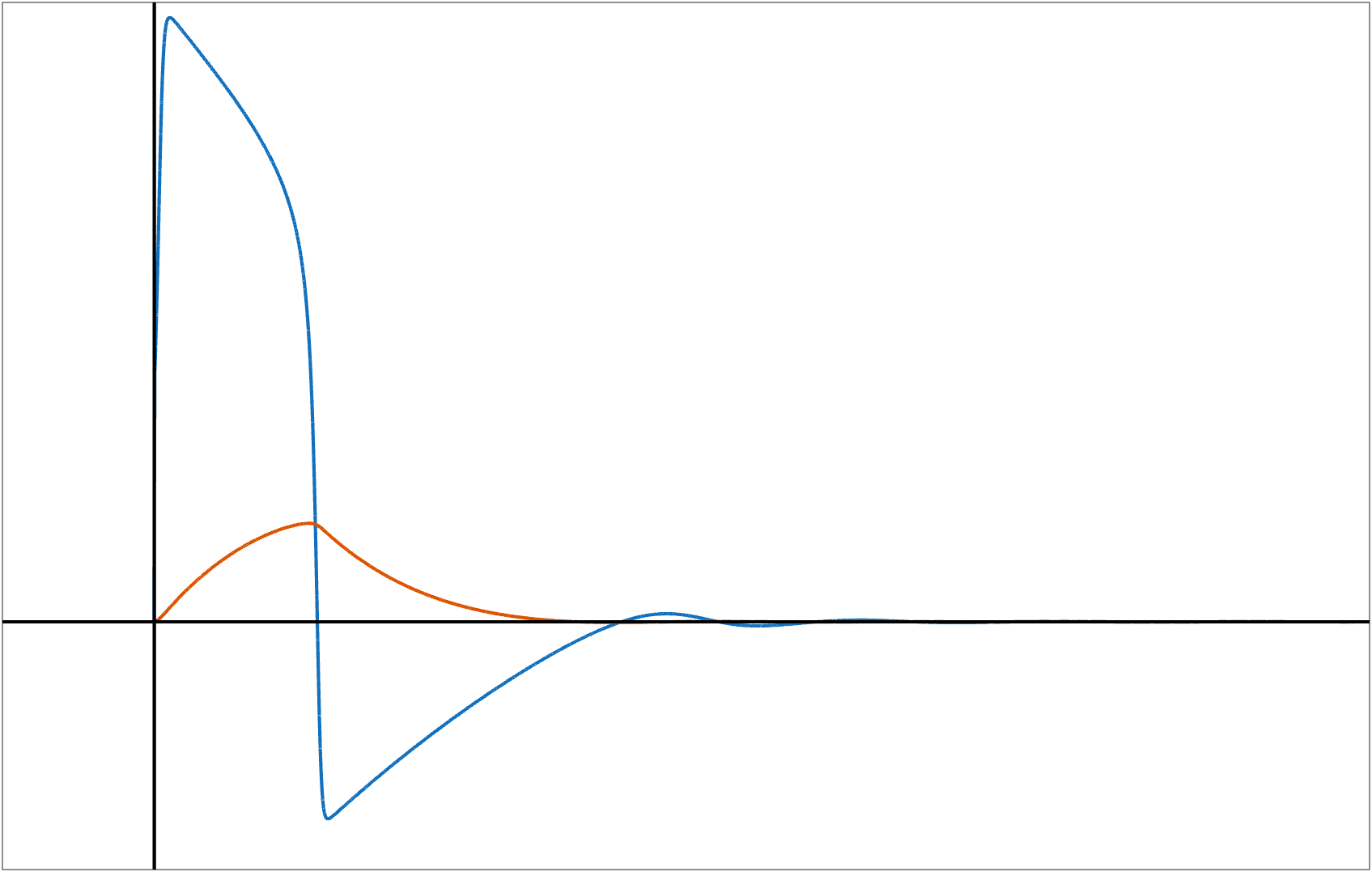}};
        \node at (\x,7/6) {\includegraphics[width=\figscale]{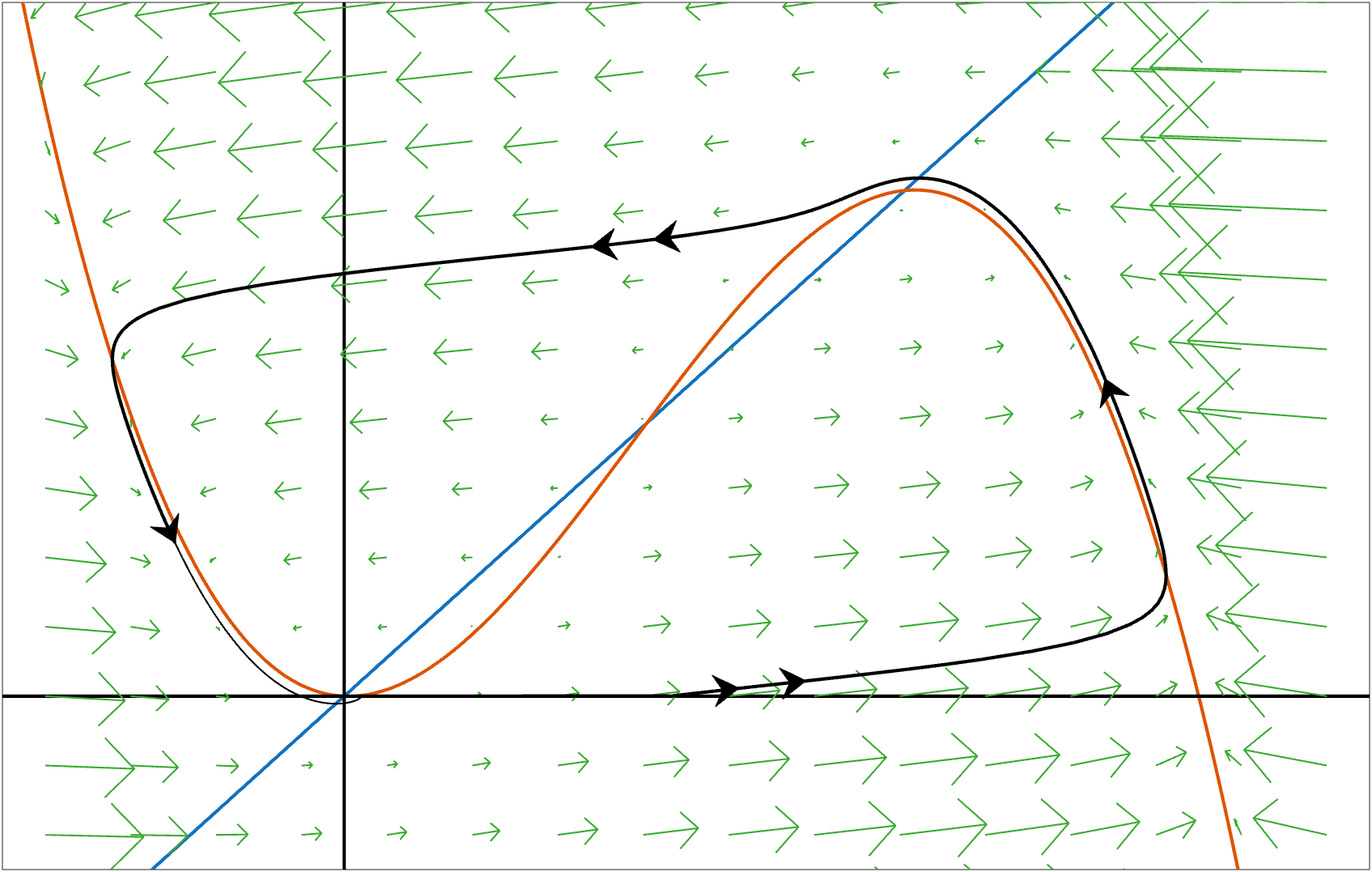}};
        \node at (\xx,7/6) {\includegraphics[width=\figscale]{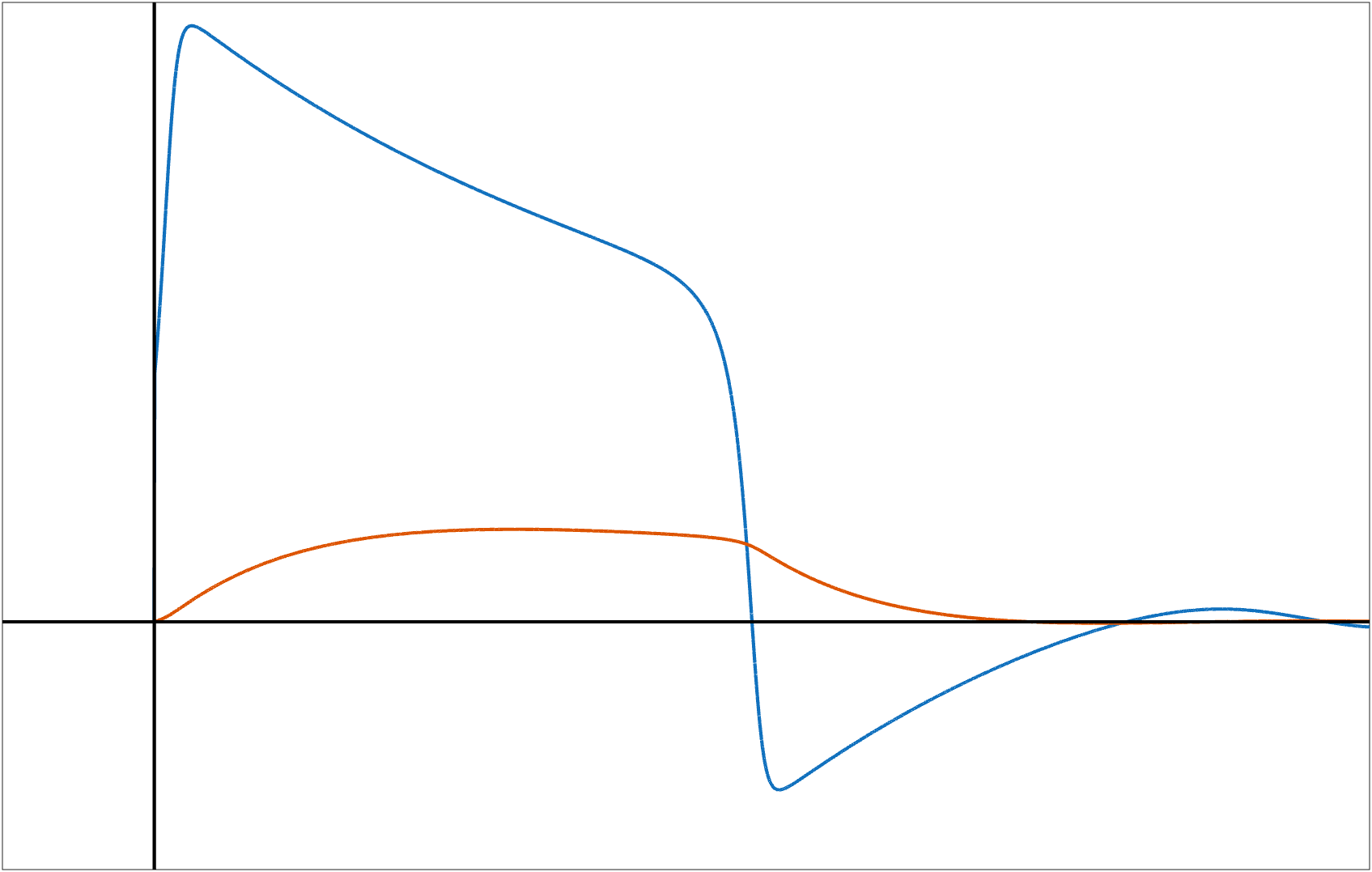}};
        \draw[->] (7,0.5+14/3) node[anchor=north] {$w$} -- ++(0,1);
        \draw[->] (7,0.5+7/3) node[anchor=north] {$w$} node[anchor=south west,rotate=90] {recovery} -- ++(0,1) ;
        \draw[->] (7,0.5) node[anchor=north] {$w$} -- ++(0,1);
        \foreach \y in {0,7/3,14/3}
        \draw[->] ({(\x+\xx)/2},\y+0.5) node[anchor=north,color=orange] {$w$} ++(0,0.5) node[anchor=north,color=blue] {$v$}-- ++(0,1);
        \draw[->] (9,-0.2) node[anchor=east] {voltage $v$} -- ++(1,0);
        \draw[->] (13,-0.2) node[anchor=east] {time $t$} -- ++(1,0);
    \end{tikzpicture}
    \caption{ {\bf Left:} The basic geometry of our model. 
    Within a region of healthy cells, there is a layer of fast-recovery cells and an adjacent layer of slow-recovery cells.
    The incoming pulse, which is traveling along the direction of stronger conductance, hits the fast-recovery cell layer first.
    {\bf Right:} The prototypical behavior of each type of cell. 
    A healthy cell fires (voltage $v$ spikes), then the recovery variable $w$ increases, then the voltage decreases, then the recovery variable decrease; once both $v$ and $w$ are near zero, the cell can fire again.
    A fast-recovery cell follows the same basic sequence as a healthy cell, but recovers faster and is ready to fire again sooner.
    A slow-recovery cell follows the same basic sequence as a healthy cell, but maintains high voltage longer and is ready to fire again later. 
    }
    \label{fig:basicgeometry}
\end{figure}

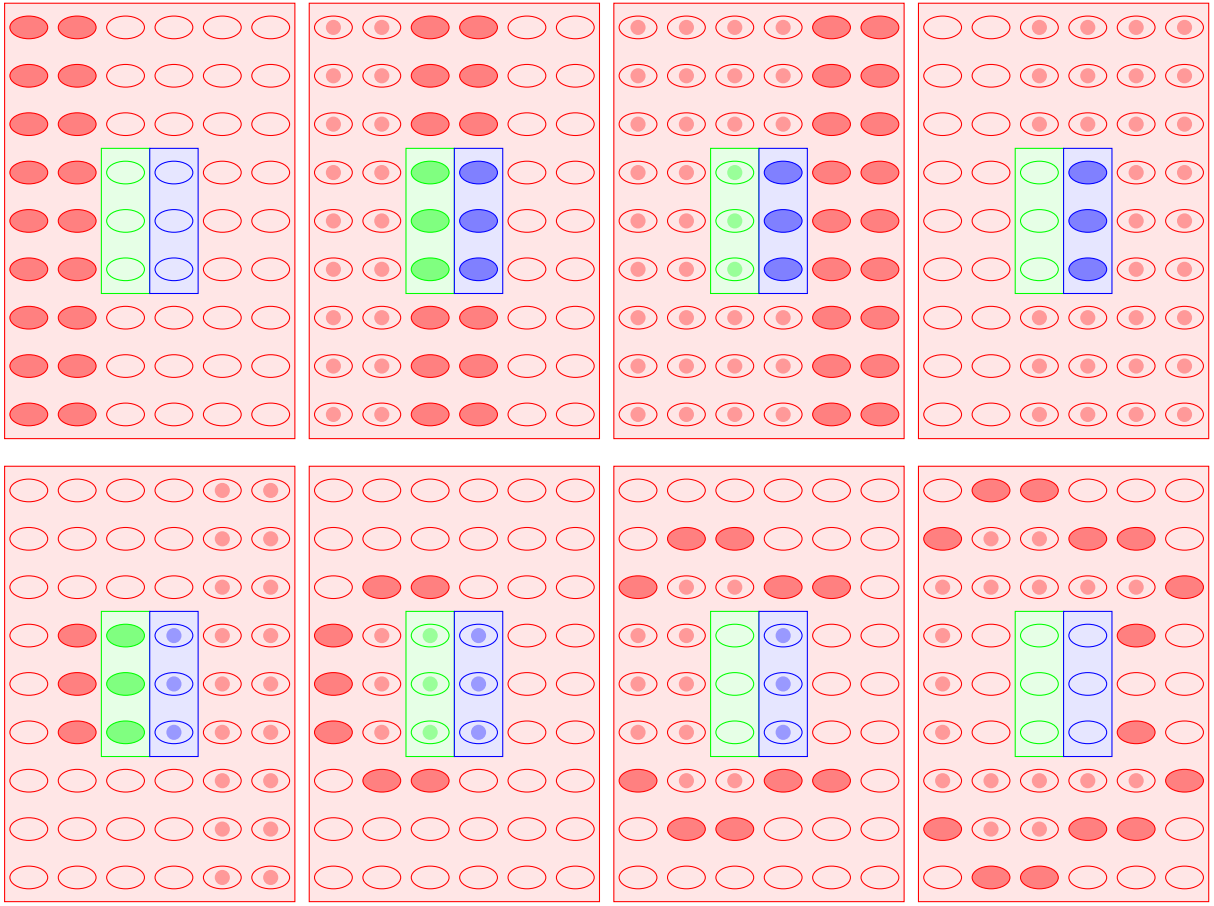
\begin{figure}[pt]
    \begin{tikzpicture} 
    \matrix[row sep=10, column sep=5, cells={scale=0.64},
    execute at begin cell={
        \path[draw,healthy] (-0.5,-1.5) rectangle +(6,9); 
        \path[draw,fast] (1.5,1.5) rectangle +(1,3); 
        \path[draw,slow] (2.5,1.5) rectangle +(1,3); }
    ] at (0,0)
        {
        \foreach \x in {0,1} \foreach \y in {-1,0,1,...,6,7} \pic at (\x,\y) {healthyfiring};
        \foreach \x in {2,3} \foreach \y in {-1,0,1,5,6,7} \pic at (\x,\y) {healthy};
        \foreach \y in {2,3,4} 
            {\pic at (2,\y) {fast};
            \pic at (3,\y) {slow};}
        \foreach \x in {4,5} \foreach \y in {-1,0,1,...,6,7} \pic at (\x,\y) {healthy};
        &
        \foreach \x in {0,1} \foreach \y in {-1,0,1,...,6,7} \pic at (\x,\y) {healthytired};
        \foreach \x in {2,3} \foreach \y in {-1,0,1,5,6,7} \pic at (\x,\y) {healthyfiring};
        \foreach \y in {2,3,4} 
            {\pic at (2,\y) {fastfiring};
            \pic at (3,\y) {slowfiring};}
        \foreach \x in {4,5} \foreach \y in {-1,0,1,...,6,7} \pic at (\x,\y) {healthy};
        &
        \foreach \x in {0,1} \foreach \y in {-1,0,1,...,6,7} \pic at (\x,\y) {healthytired};
        \foreach \x in {2,3} \foreach \y in {-1,0,1,5,6,7} \pic at (\x,\y) {healthytired};
        \foreach \y in {2,3,4} 
            {\pic at (2,\y) {fasttired};
            \pic at (3,\y) {slowfiring};}
        \foreach \x in {4,5} \foreach \y in {-1,0,1,...,6,7} \pic at (\x,\y) {healthyfiring};
        &
        \foreach \x in {0,1} \foreach \y in {-1,0,1,...,6,7} \pic at (\x,\y) {healthy};
        \foreach \x in {4,5} \foreach \y in {-1,0,1,...,6,7} \pic at (\x,\y) {healthytired};
        \foreach \x in {2,3} \foreach \y in {-1,0,1,5,6,7} \pic at (\x,\y) {healthytired};
        \foreach \y in {2,3,4} \pic at (2,\y) {fast};
        \foreach \y in {2,3,4} \pic at (3,\y) {slowfiring};
        \\
        \foreach \y in {-1,0,1,...,6,7} \pic at (0,\y) {healthy};
        \foreach \y in {-1,0,1,5,6,7} \pic at (1,\y) {healthy};
        \foreach \y in {2,3,4} \pic at (1,\y) {healthyfiring};
        \foreach \x in {2,3} \foreach \y in {-1,0,1,5,6,7} \pic at (\x,\y) {healthy};
        \foreach \y in {2,3,4} \pic at (2,\y) {fastfiring};
        \foreach \y in {2,3,4} \pic at (3,\y) {slowtired};
        \foreach \x in {4,5} \foreach \y in {-1,0,1,...,6,7} \pic at (\x,\y) {healthytired}; 
        &
        \foreach \y in {-1,0,1,5,6,7} \pic at (0,\y) {healthy};
        \foreach \y in {2,3,4} \pic at (0,\y) {healthyfiring};
        \foreach \y in {-1,0,6,7} \pic at (1,\y) {healthy};
        \foreach \y in {1,5} \pic at (1,\y) {healthyfiring};
        \foreach \y in {2,3,4} \pic at (1,\y) {healthytired};
        \foreach \y in {1,5} \pic at (2,\y) {healthyfiring};
        \foreach \y in {-1,0,6,7} \pic at (2,\y) {healthy};
        \foreach \y in {2,3,4} \pic at (2,\y) {fasttired};
        \foreach \y in {-1,0,1,5,6,7} \pic at (3,\y) {healthy};
        \foreach \y in {2,3,4} \pic at (3,\y) {slowtired};        
        \foreach \x in {4,5} \foreach \y in {-1,0,1,...,6,7} \pic at (\x,\y) {healthy};
        &
        \foreach \y in {-1,0,6,7} \pic at (0,\y) {healthy};
        \foreach \y in {1,5} \pic at (0,\y) {healthyfiring};
        \foreach \y in {2,3,4} \pic at (0,\y) {healthytired};
        \foreach \y in {-1,7} \pic at (1,\y) {healthy};
        \foreach \y in {1,2,3,4,5} \pic at (1,\y) {healthytired};
        \foreach \y in {0,6} \pic at (1,\y) {healthyfiring};
        \foreach \y in {-1,7} \pic at (2,\y) {healthy};
        \foreach \y in {0,6} \pic at (2,\y) {healthyfiring};
        \foreach \y in {1,5} \pic at (2,\y) {healthytired};
        \foreach \y in {2,3,4} \pic at (2,\y) {fast};
        \foreach \y in {-1,0,6,7} \pic at (3,\y) {healthy};
        \foreach \y in {1,5} \pic at (3,\y) {healthyfiring};
        \foreach \y in {2,3,4} \pic at (3,\y) {slowtired};      
        \foreach \x in {4} \foreach \y in {-1,0,2,3,4,6,7} \pic at (\x,\y) {healthy};
        \foreach \x in {4} \foreach \y in {1,5} \pic at (\x,\y) {healthyfiring};
        \foreach \x in {5} \foreach \y in {-1,0,1,...,6,7} \pic at (\x,\y) {healthy};
        &
        \foreach \y in {-1,7} \pic at (0,\y) {healthy};
        \foreach \y in {0,6} \pic at (0,\y) {healthyfiring};
        \foreach \y in {1,2,3,4,5} \pic at (0,\y) {healthytired};
        \foreach \y in {-1,7} \pic at (1,\y) {healthyfiring};
        \foreach \y in {0,1,5,6} \pic at (1,\y) {healthytired};
        \foreach \y in {2,3,4} \pic at (1,\y) {healthy};
        \foreach \y in {-1,7} \pic at (2,\y) {healthyfiring};
        \foreach \y in {0,1,5,6} \pic at (2,\y) {healthytired};
        \foreach \y in {2,3,4} \pic at (2,\y) {fast};
        \foreach \y in {-1,7} \pic at (3,\y) {healthy};
        \foreach \y in {0,6} \pic at (3,\y) {healthyfiring};
        \foreach \y in {1,5} \pic at (3,\y) {healthytired};
        \foreach \y in {2,3,4} \pic at (3,\y) {slow};
        \foreach \y in {-1,3,7} \pic at (4,\y) {healthy};
        \foreach \y in {0,2,4,6} \pic at (4,\y) {healthyfiring};
        \foreach \y in {1,5} \pic at (4,\y) {healthytired};
        \foreach \y in {-1,0,2,3,4,6,7} \pic at (5,\y) {healthy};
        \foreach \y in {1,5} \pic at (5,\y) {healthyfiring};
        \\
        };
    \end{tikzpicture}
    \caption{
    Cartoon of the sequence of events in a Fibrillatory Initiating Boundaries (FIBs) model with the geometry of \cref{fig:basicgeometry}.
    An open ellipse 
    \healthyicon \
    indicates a region of recovered cells ready to fire; the longer horizontal axis is a reminder that conductance is stronger in the horizontal direction.
    A filled ellipse \healthyfiringicon \
    indicates a region of cells that have fired and still have high-enough voltage to trigger a recovered neighboring cell.
    A partly-filled ellipse \healthytiredicon \
    indicates a region of cells that are in the recovery process and do not have high-enough voltage to trigger neighboring cells and are not recovered enough to be triggered themselves.
    In the first row, a wave passes through from left to right. 
    The slow-recovery cells still have high voltage and the fast-recovery cells are ready to fire again. 
    In the second row, the fast-recovery cells fire and initiate a second wave, which travels around the slow-recovery cells and forms swirling-back curves.
    By the end of the second row, the wave has reached the back of the slow-recovery layer, which is now ready to fire again and initiate a third wave.}
    \label{fig:basicsequence}
\end{figure}

We construct a simple model that includes healthy and unhealthy cells with realistic recovery times and realistic geometric arrangements of the unhealthy cells.
We study in detail the basic geometry shown in \cref{fig:basicgeometry}, which is simply a layer of fast-recovery cells adjacent to a layer of slow-recovery cells, both contained within a larger region of healthy cells.
\Cref{fig:basicsequence} gives a cartoon version of the dynamics we observe in the basic geometry of \cref{fig:basicgeometry}.
After a single, normal wave passes, the slow-recovery cells still have high voltage and the fast-recovery cells are ready to fire again.
The fast-recovery cells fire a second time and initiate a second wave.
Note, crucially, that the fast-recovery cells fired due to the sustained high voltage in the plateau phase (phase 2) of the slow-recovery cells, rather than due to an initial spike in voltage (phase 0) of a neighboring cell as in a normal beat.
This second wave travels around the slow-recovery cells and forms swirling-back curves.
Note, again crucially, that it is the ``broken ends'' of this second wave that create the swirling-back curve, not the front of the wave. 
The slow-recovery cells are then ready to fire again, and a third wave can start.
Depending on interactions of timing and geometry, persistent but irregular generation of waves may occur. 
These dynamics have the key features matching the phenomena observed during atrial fibrillation.
Therefore, we conjecture that some instances of atrial fibrillation are caused by unhealthy cells in this geometry.

We call our model {\em fibrillatory initiating boundaries}\ (FIBs) because its key assumption is that there is a {\em boundary} between the layers of fast- and slow-recovery cells.
With appropriate parameters, this mismatch in recovery times allows the slow-recovery cells to trigger the fast-recovery cells and start the second wave, which is the beginning of the arrhythmia.
A second assumption is that the 2-dimensional geometry requires the second wave to take enough time to reach back to the slow-recovery cells that they can recover. 
This time-delay allows the fast-recovery cells to indirectly trigger the slow-recovery cells and provides a mechanism for the arryhthmia to sustain itself and become a fibrillation.
Once the slow-recovery cells have fired a second time, the fact that they are slow-recovery is not essential to maintaining the fibrillation, but it can provide a backup mechanism to restart the fibrillation if it falters.
We emphasize that the FIB mechanism is spontaneous in the sense that only one normal stimulation is needed.

Heterogeneity has been thought to be a key to atrial fibrillation at least since 1964, when it was used in a computer model \cite{MO-RH-AB:1964} that started the `wavelet' theory of fibrillation.
Two mechanisms have been described and studied.
The first mechanism is slow-recovery cells maintaining elevated voltage in their plateau phase and re-trigger neighboring cells \cite{ANT-BUR:2011,H-S-V-H:2021,R-R-B-B-T:2023}. 
In models with randomly distributed cell types, scattered slow-recovery cells can cause disordered refirings \cite{MO-RH-AB:1964,H-S-V-H:2021}.
In models with sharp boundaries between cell types, the slow-recovery region can trigger repeated firing of the other region \cite{T-D-V-P-P:2018,A-A-K-Z-M-Betal:2019,LI-QU-HU:2023}.
The second mechanism is a difference in cell properties causing a difference in wave speed, which allows a wave to lose coherence and break apart \cite{A-A-K-Z-M-Betal:2019,G-P-W-G-etal:2019}.
If the recovery phase of a faster wave passes the firing phase of an adjacent slower wave, then the broken end of the slower wave causes sideways and backwards firings behind the faster wave.
%
%
Thus, several aspects of our model and its observed dynamics are supported by previous work in the literature.

In \cref{sec:FHNmodel} we describe the FitzHugh-Nagumo model we use and the choices of parameters to model heart cells and their connections into a two-dimensional array.
In \cref{sec:protosimulation}, we describe two-dimensional simulations using this model and highlight one that illustrates the dynamics sketched as a cartoon in \cref{fig:basicsequence}.

In \cref{sec:parameterinteractions}, we study the dependence on the parameters used for the fast-recovery and slow-recovery cells.
By varying the parameters that determine the speed at which the  slow-recovery cells and fast-recovery cells recover, we produce heat maps showing when persistent wave generation occurs.
We first observe, as noted in the cartoon sequence, that the slow-recovery cells must maintain high voltage long enough and the fast-recovery cells recover fast enough that the fast-recovery cells are triggered to fire again via the plateau phase stimulation.
In \cref{sec:localtiming}, we extract these times from single-cell models of fast-recovery and slow-recovery cells and use them to qualitatively validate  features of the heat maps.
Second, we observe that if the slow-recovery cells retain high voltage too long, then they are repeatedly stimulated by the adjacent fast-recovery cells and never fully recover;
in this non-physiologic state, these slow-recovery cells would eventually die (due to \ce{Ca^{2+}} overload \cite{C-Z-K-H-M-Metal:2005}). 
In \cref{sec:localnonphysslow}, we study a two-cell model coupling a fast-recovery cell to a slow-recovery cell and show that increasing the recovery time causes a reverse period-doubling bifurcation that creates a stable periodic orbit where the fast-recovery cell fires and recovers repeatedly, while the slow-recovery cell oscillates at a high voltage. 
In the two-cell model, a plateau phase stimulation  can occur when $b$ is between the initial saddle-node and the reverse period-doubling bifurcations, opening the possibility of interesting dynamics in a grid of cells.
All of these observations help explain the structure seen in the heat maps.

In \cref{sec:variationgeometry}, we report on simulations that test the robustness of the FIB mechanism by varying the conductance, by varying the size and proportions of the rectangular region of unhealthy cells, and by considering an L-shaped FIB.
We find that the qualitative behavior is unchanged for conductance anisotropies up to about 4, but for higher anisotropies the self-sustaining behavior is not robust.
We find that taller FIBs increase the robustness of the FIB mechanism by both helping the slow-recovery cells maintain voltage (by reducing charge  leakage to non-slow-recovery neighbors) and by lengthening the travel time of the wave from the fast-recovery layer to the back of the slow-recovery layer.
For a small FIB, the region in parameter space where sustained oscillations occur is rather thin and fractal-like, whereas for the taller FIB it is solid and large. 
This is consistent with the usual progression of disease from normal function to occasional AFib to  intermittent AFib to persistent AFib as regions of heterogeneity grow. 
We find that wider slow-recovery layers can lead to non-recovery, while wider fast-recovery layers increase robustness by providing an additional refiring mechanism within the fast-recovery layer.
Finally, we find that an L-shaped FIB is also more robust, through a combination of the mechanisms for a taller FIB and for a wider fast-recovery layer.
Overall, these tests show that these variations change the simulation results in sensible ways, so we conclude that the FIB mechanism is robust and progressive.

In \cref{sec:exoticsimulations}, we report on simulations using  more exotic geometries, artificial initial conditions, and an additional cell type (scars).
With these, we can produce spiral-like waves, phantom fibrillatory centers, and multiple wave(let)s interacting.
In several cases, the simulation results look like results from other models or observations from real atria. 
Generally, the FIB structure we propose causes excitation waves with ``broken ends" or points of phase discontinuity, which have been identified as the centers of spiral waves. 
In many situations, the broken ends appear as part of swirling-back curves and in some instance with more asymmetry as the center of spiral waves.
We observe that a homogeneous network of healthy cells in our model can sustain spiral waves, that the apparent center of a disturbance can drift away from its initiating FIB source, and that a FIB can mask the presence of other FIBs. 
These facts serve as a cautionary note in patient care: ablation can easily miss the real cause of an arrhythmia.

In this paper we describe the phenomena that the model exhibits and give mathematical explanations for them.
In \cref{sec:conclusion} we summarize our observations and conclusions.
A future paper will address approaches for detecting this model geometry in actual hearts and the clinical implications of the dynamics we described.


\section{FitzHugh–Nagumo Model of Heart Cells}
\label{sec:FHNmodel}

\subsection{Model for a Single Cell}

We model the electrical behavior of individual heart cells using a FitzHugh-Nagumo (FHN) model:
\begin{equation}\label{eqn:FHN}
\begin{split}
\frac{dv}{dt} &= \sigma \left( v ( v - a)(c -v) - w \right)  \\
\frac{dw}{dt} & = \epsilon (v - b w).
\end{split}
\end{equation}
There are several equivalent versions of the model in the literature \cite{RO-GE-GI:2000,LI-QU-HU:2023,GO-LA-RO:2025}.
Here, $v$ is the cell's voltage and $w$ is a `recovery variable'. 
The parameters $\sigma >0$ and $\epsilon>0$ set the overall speed of $v$ and $w$.
The parameter $a>0$ is a firing threshold for $v$ and $c>a$ is a high-voltage limit for $v$.
The parameter $b>0$ sets the slope of the $w$-nullcline. 
Increasing $b$ makes recovery slower and eventually creates a new stable equilibrium in the model where $v$ can become stuck in a high voltage (depolarized) state.

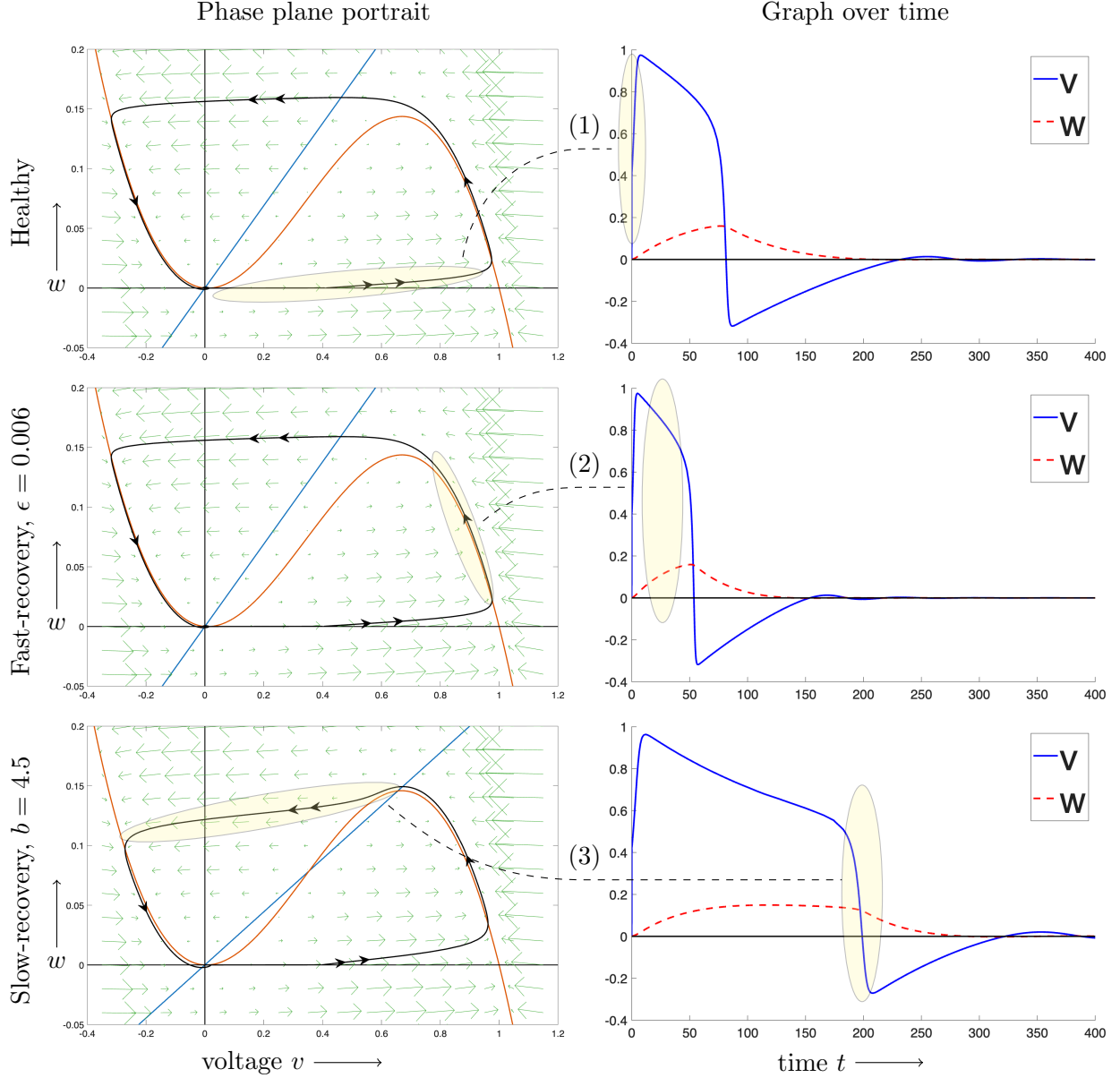
\begin{figure}
\begin{tikzpicture}
        \setlength{\figscale}{0.45\textwidth}
        \def\x {4};
        \def\xx {12};
        \def\y {15}
        \node[anchor=south] at (\x,\y) {Phase plane portrait};
        \node[anchor=south] at (\xx,\y) {Graph over time};
        \node[rotate=90,anchor=south] at (0,{\y*5/6}) {Healthy};
        \node at (\x,{\y*5/6}) {{\includegraphics[width=\figscale]{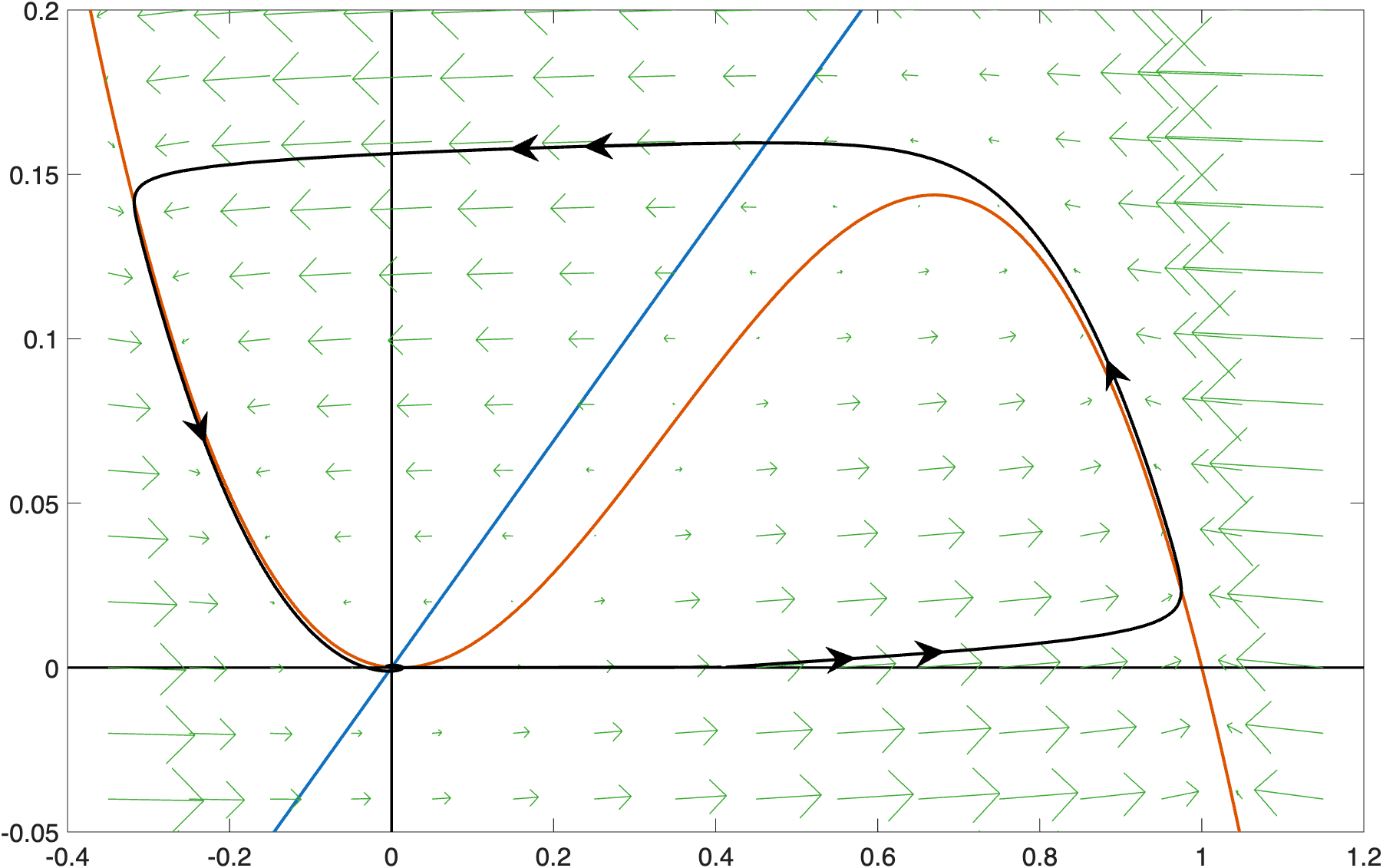}}};
        \node at (\xx,{\y*5/6}) {{\includegraphics[width=\figscale]{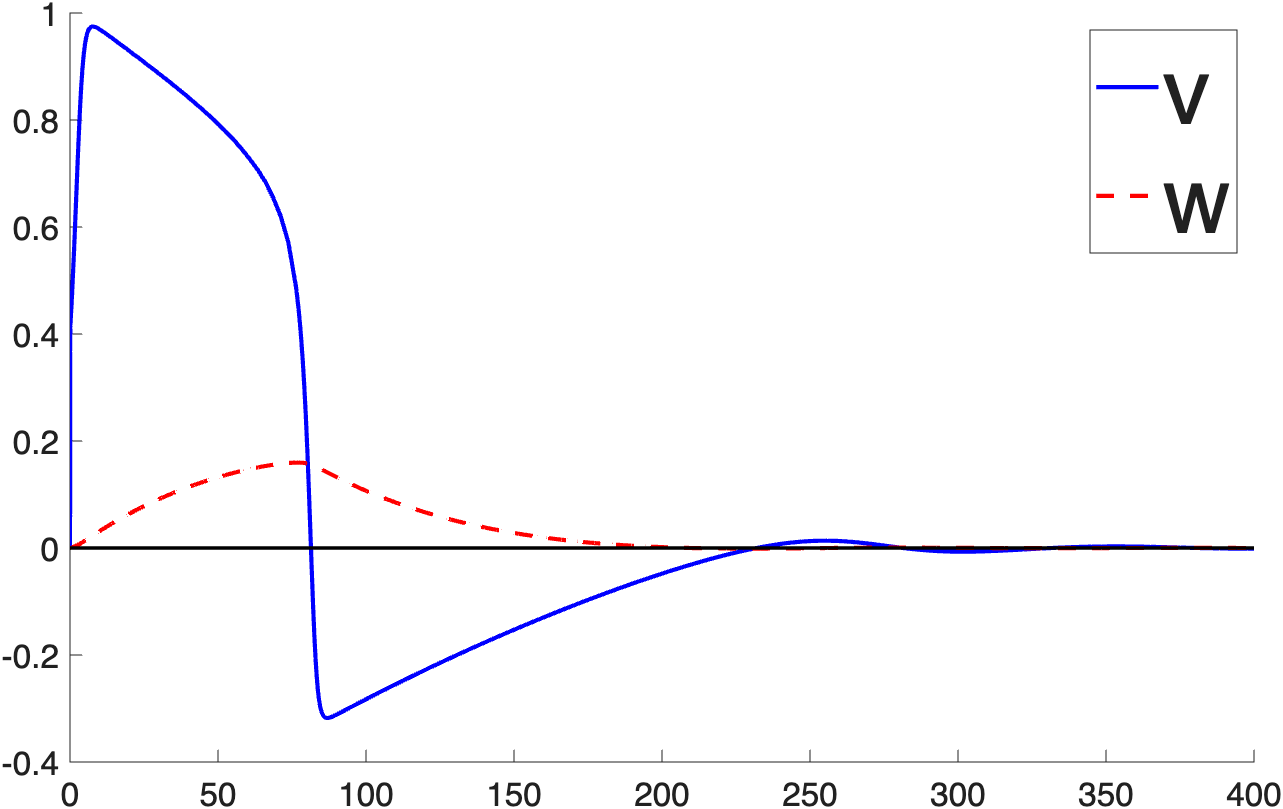}}};
        \draw[fill=yellow!50,nearly transparent] ({\x+0.5},{\y*5/6 -1.2}) circle[dashed,x radius=2,y radius=0.2,rotate=5];
        \path[draw,dashed] ({\x+2.2},{\y*5/6 -0.8}) to [out=80,in=180] ({\xx-4},{\y*5/6+0.8}) node[anchor=south] {(1)} --({\xx-3.6},{\y*5/6+0.8}); 
        \draw[fill=yellow!50,nearly transparent]  ({\xx-3.3},{\y*5/6+0.8}) circle[dashed,y radius=1.4,x radius=0.2];
        \node[rotate=90,anchor=south] at (0,{\y*1/2}) {Fast-recovery, $\epsilon=0.006$};
        \node at (\x,{\y*1/2}) {{\includegraphics[width=\figscale]{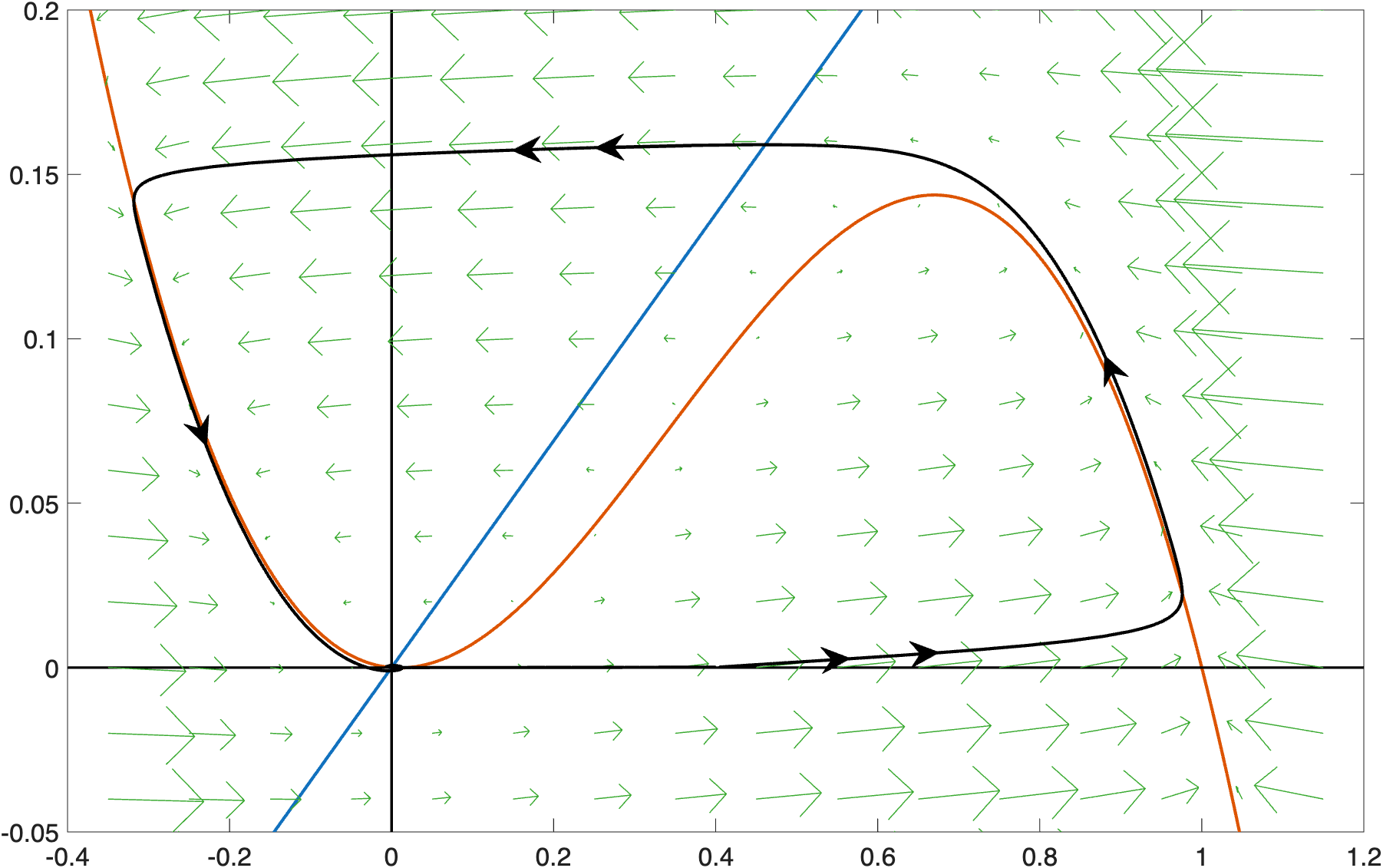}}};
        \node at (\xx,{\y*1/2}) {{\includegraphics[width=\figscale]{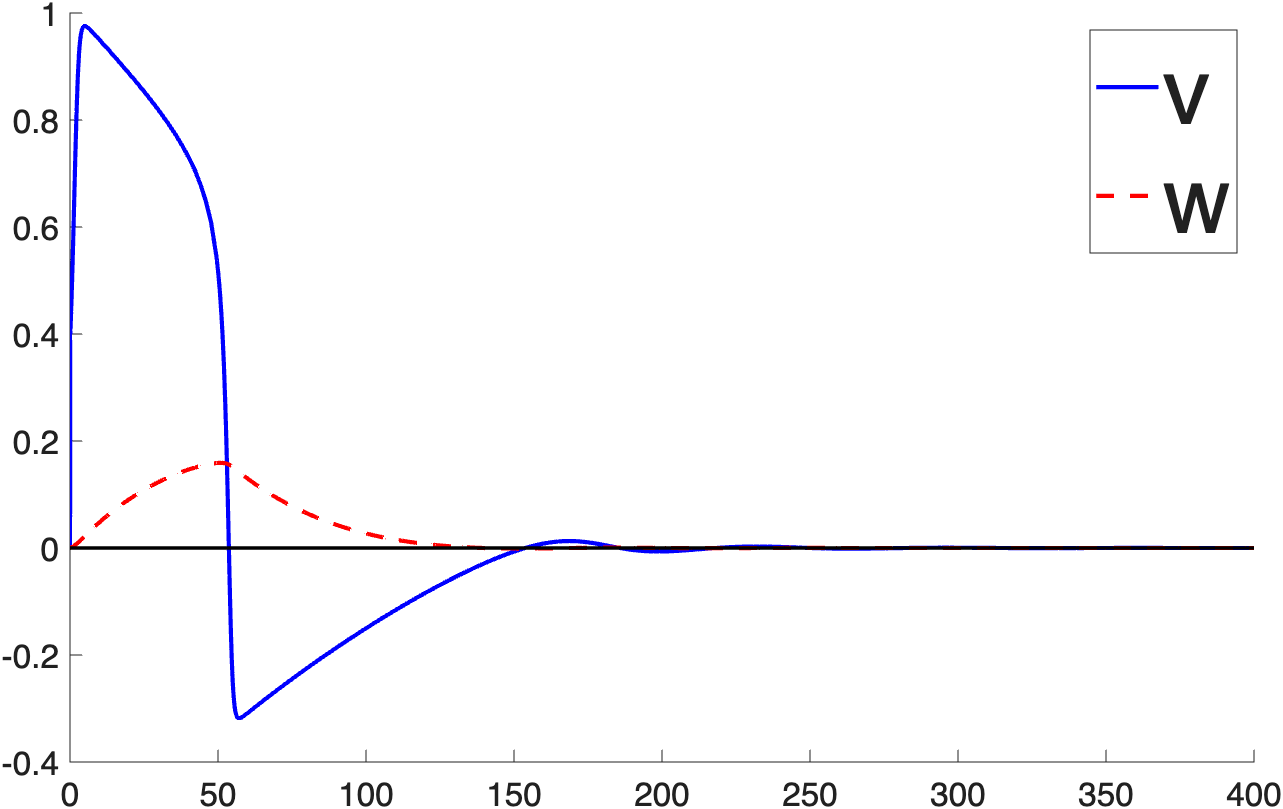}}};
        \draw[fill=yellow!50,nearly transparent] ({\x+2.2},{\y*1/2 +0.2}) circle[dashed,x radius=0.2,y radius=1.2,rotate=20];
        \path[draw,dashed] ({\x+2.5},{\y*1/2 +0.3}) to [out=40,in=180] ({\xx-4},{\y*1/2+0.8}) node[anchor=south] {(2)} --({\xx-3.4},{\y*1/2+0.8}); 
        \draw[fill=yellow!50,nearly transparent]  ({\xx-2.85},{\y*1/2+0.6}) circle[dashed,y radius=1.8,x radius=0.3];
        \node[rotate=90,anchor=south] at (0,{\y*1/6}) {Slow-recovery, $b=4.5$};
        \node at (\x,{\y*1/6}) {{\includegraphics[width=\figscale]{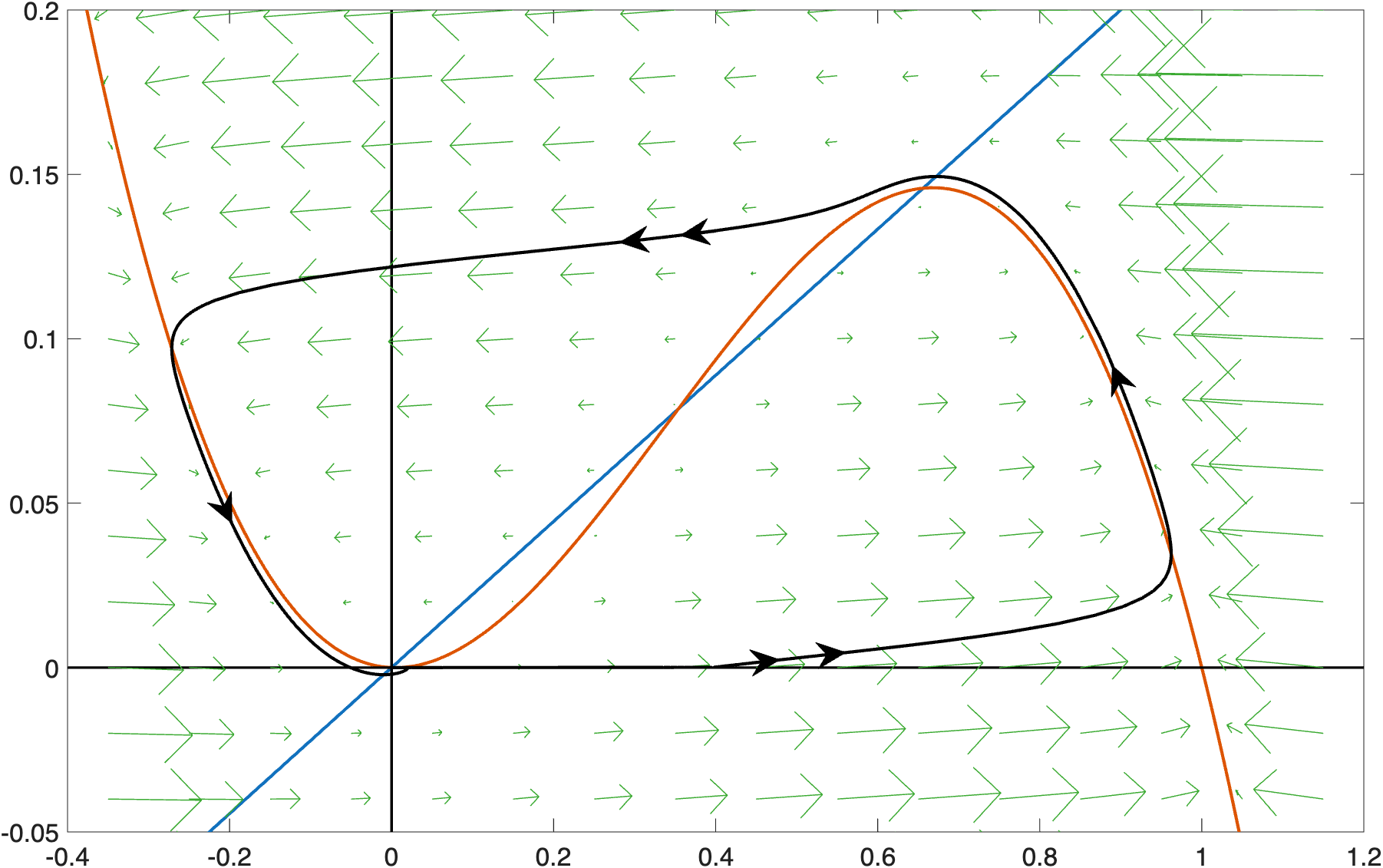}}};
        \node at (\xx,{\y*1/6}) {{\includegraphics[width=\figscale]{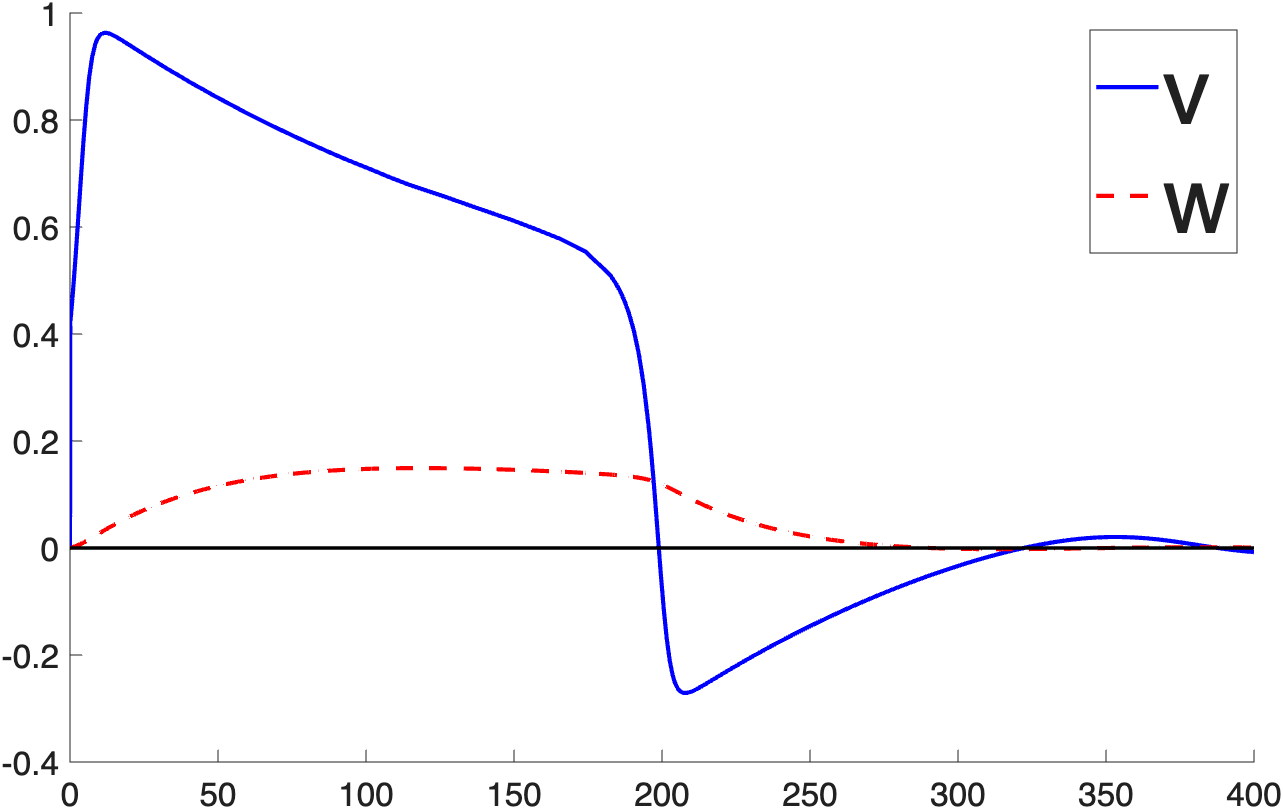}}};
        \draw[fill=yellow!50,nearly transparent] ({\x-0.8},{\y*1/6 +1.0}) circle[dashed,x radius=2.1,y radius=0.25,rotate=10];
        \path[draw,dashed] ({\x+1.1},{\y*1/6 +1.1}) to [out=320,in=180] ({\xx-4},{\y*1/6}) node[anchor=south] {(3)} --({\xx-0.2},{\y*1/6}); 
        \draw[fill=yellow!50,nearly transparent]  ({\xx+0.1},{\y*1/6-0.2}) circle[dashed,y radius=1.6,x radius=0.3];
 \draw[->] (\x,-0.2) node[anchor=east] {voltage $v$} -- ++(1,0);
        \draw[->] (\xx,-0.2) node[anchor=east] {time $t$} -- ++(1,0);
        \foreach \yy in {{\y/6},{\y/2},{5*\y/6}}
           \draw[->] (0.2,{\yy-1}) node[anchor=north] {$w$}-- ++(0,1);
\end{tikzpicture}
\caption{Phase plane portraits and graphs over time for cells in the FitzHugh-Nagumo model.
In all cell types, (1) 
a rapid spike in $v$  occurs, followed by (2) growth in $w$, followed by (3) 
rapid decrease in $v$, and concluded with return to $(v,w)=(0,0)$. 
 {\bf Left}: Phase plane portrait of the FHN model for a cell. 
The blue line is the nullcline for $w$ (where $w'=0$) and the orange cubic is the nullcline for $v$ (where $v'=0$). 
The black curve is the trajectory followed by the model after being subjected to an initial impulse.
{\bf Right}: Voltage $v$ and recovery variable $w$ of the cell in response to the same initial impulse. 
\label{fig:pp_trace_3}}
\end{figure}

The FHN cell model is known as an `excitable' system. 
In \cref{fig:pp_trace_3} we show the phase planes for a cell of each type.
For a healthy cell, $(v,w)=(0,0)$ is an asymptotically stable equilibrium and is the only fixed point or periodic solution of the system, and it is globally stable. 

\Cref{fig:pp_trace_3} also shows the affect of a single input pulse, as both trajectories superimposed on the phase planes and graphs over time. 
When the voltage $v$ is perturbed (e.g.\ by a current $I$) so that it increases above a threshold (represented by $a$), then the FHN model will produce a rapid spike in the voltage, $v$.
When this occurs, the cell is said to have exhibited an `action potential'.
It is this spike that triggers mechanical contraction of the cell (which is not modeled here).
After that, the trajectory will approximately follow the $v$-nullcline in the phase plane, with decreasing $v$ and increasing $w$, to near the local maximum point of the $v$-nullcline. 
From there, the trajectory will move quickly to the left, with $v$ rapidly decreasing and becoming negative, to the vicinity of the other part of $v$-nullcline.  
From there it will follow the leftmost branch of the $v$-nullcline downward, with increasing $v$ and decreasing $w$, back to $(0,0)$. 
The final two phases of this process are referred to as `recovery'.  
During recovery, a voltage equal to the distance between the left two portions of the $v$-nullcline would be needed to make the cell fire again.

\subsection{Grids of Coupled Cells}

In a functioning heart, cells are connected electrically through gap junctions that allow the rapid passing of ions and effectively act as conductive connections between neighboring cells.
The single-cell model \cref{eqn:FHN} thus becomes a FHN coupled model
\begin{equation}\label{eqn:FHNcoupled}
\begin{split}
\frac{dv_i}{dt} &= \sigma_i \left( v_i ( v_i - a_i)(c_i -v_i) - w_i +\sum_j g_{ij}(v_j-v_i)\right)  \\
\frac{dw_i}{dt} & = \epsilon_i (v_i - b_i w_i)\,,
\end{split}
\end{equation}
where $g_{ij}$ is a conductance coefficient between cells $i$ and $j$ (zero if they are not connected). 
Heart cells are elongated and have more gap junctions on the ends of cells than on the sides, leading to stronger conductance between cells that are connected at their ends and thus $g_{ij}$  is larger for pairs that are connected at their ends.


\begin{Remark}
  While a standard nondimensionalization procedure can eliminate two additional parameters in the single-cell model \cref{eqn:FHN}, for instance, by setting $\sigma = c = 1$, such a procedure does not extend naturally to a network of cells of heterogeneous types, where parameter values may differ across cells. 
\end{Remark}

\subsection{Calibration of Parameters to Model Heart Cells}



The FitzHugh-Nagumo model is phenomenological and we chose parameters for the model to capture the main qualitative features of cardiac action potentials and recovery. These choices are consistent with choices in the literature \cite{LI-QU-HU:2023, L-I-D-Z-B:2009}.
Note that the time units are not scaled to match typical heartbeats.
\Cref{tab:parameters} shows the parameter values we use for healthy, fast-recovery, and slow-recovery cells.
\begin{table}
\centering
\begin{tabular}{|c|c|c|c|c|c|} \hline
\  Type \ & $a$ & $b$ & $c$ & $\epsilon$ & $\sigma$  \\ \hline
 \ Healthy\ & \ 0.02 \ & \ 2.9 \  & \ 1.0 \ & \ 0.004 \ & \ 1.0 \  \\
\ Fast-recovery \ & 0.02 & 2.9 & 1.0 & \ 0.006+ \ & 1.6  \\ 
\ Slow-recovery \ & 0.01 & \ 4.5+ \ & 1.0 & 0.004 & 0.6 \\ 
\hline
\end{tabular}
\caption{Parameter choices for 3 types of cells within the FHN model \cref{eqn:FHN}. The $+$ indicates that we try a range of values above the given number.
}
\label{tab:parameters}
\end{table}
Changes to the $\sigma$ parameter may be interpreted as changes in the capacitance of cells, which correspond directly to physical changes in the myocytes \cite{SA-CH-LE:2010}.
We note that the FitzHugh-Nagumo model that we use does not include changes in recovery time due to the frequency of stimulation (restitution effects \cite{MO-RH-AB:1964,NOL-DAH:1968}), so if one wants to model the effects of multiple fast stimulations (pacing), then the $\epsilon$ parameter of the fast-recovery cells should be set to achieve the post-restitution recovery time.

The fast-recovery cells differ from healthy cells in two parameters.
First, they have a larger $\sigma$ parameter, which directly scales the derivative of the voltage variable $v$, making it rise or fall more quickly.
Second, they have a larger $\epsilon$ parameter, which directly scales the derivative of the recovery variable $w$, making it also rise or fall more quickly.
In our experiments, we will vary the parameter $\epsilon$.
Combined, these speed up the entire cycle, while qualitatively preserving the shape and leaving the nullclines unchanged.

The slow-recovery cells differ from healthy cells in three parameters.
First, they have a lower $\sigma$ parameter, which directly scales the derivative of the voltage variable $v$, making it rise or fall more slowly.
Second, they have a larger $b$ parameter, which controls the slope of the $w$ nullcline and allows the cells to stay longer at higher voltages.
Third, they have a lower $a$ parameter, which extends their plateau phase by making the trajectory pass close by a slow, unstable focus and, incidentally, makes them easier to cause to fire.
In our experiments, we will vary this $b$ parameter.
Combined, these slow down the entire cycle, and especially extend the time in which $v$ is high.
The decrease in $a$ slightly changes the $v$-nullcline, while the increase in $b$ significantly reduces the slope of the $w$-nullcline and for large values of $b$ creates two new intersections of the nullclines.
These are new equilibrium points in the single cell model. If we increase $b$ further, one of these fixed points will become stable. 
We analyze the bifurcations involved in this process in \cref{sec:localnonphysslow}.




When working on a two-dimensional grid of cells, we set the horizontal and vertical conductance coefficients to be $g_{\mathrm{horz}} = 0.027$ and $g_{\mathrm{vert}} = 0.015$, giving an anisotropy ratio of 1.8. 
Clinical estimates of anisotropy ratios vary from greater than 2 in 2007 \cite{VALDER:2007}, to similar to 1.8 in 2013 \cite{M-Z-A-E-etal:2013}, to smaller than 1.8 in tissues where atrial fibrillation was occurring in 2024 \cite{S-T-K-T-G:2024}.
While we are trying to be consistent with this literature, the choice was essentially phenomenological.


\section{Key observations in a prototypical simulation}
\label{sec:protosimulation}

In this section we describe in detail one simulation of a FIB, corresponding to the cartoon in \cref{fig:basicgeometry}.
In \cref{sec:parameterinteractions} we will consider the effects of changing the cell parameters and in \cref{sec:variationgeometry} we will consider the effects of changing the geometry and grid size.
\begin{figure}
    \centering
    \begin{tikzpicture}[scale=0.2]
       \path[draw,help lines] (1,1) grid +(41,41); 
       {\foreach \x in {1,...,18,23,24,...,42} \foreach \y in {1,...,42} \draw[healthy,fill=red] (\x,\y) circle (0.15);} 
       \foreach \x in {19,20,21,22} \foreach \y in {1,...,17,25,26,...,42} \draw[healthy,fill=red] (\x,\y) circle (0.15);
       \foreach \x in {19,20} \foreach \y in {18,...,24} \draw[fast,fill=green] (\x-0.3,\y-0.3) rectangle +(0.6,0.6);
       \draw[->] (10,0) node[anchor=north] {Fast-recovery} -- (18.5,17.5);
       \foreach \x in {21,22} \foreach \y in {18,...,24} \draw[slow,fill=blue] (\x,\y) circle (0.4);
       \draw[->] (30,0) node[anchor=north] {Slow-recovery} -- (22.5,17.5);
       \draw[fill=yellow!50,nearly transparent] (0,0) rectangle +(1.5,43);
       \node[anchor=east] at (0,21) {Initial Voltage Stimulus};
    \end{tikzpicture}
    \caption{The computational grid used in the simulation corresponding to the cartoon basic geometry in \cref{fig:basicgeometry} and cartoon sequence in \cref{fig:basicsequence}.
    Healthy cells are indicated by small red discs, fast-recovery cells by green squares, and slow-recovery cells by large blue discs.
    }
    \label{fig:simbasicgrid}
\end{figure}
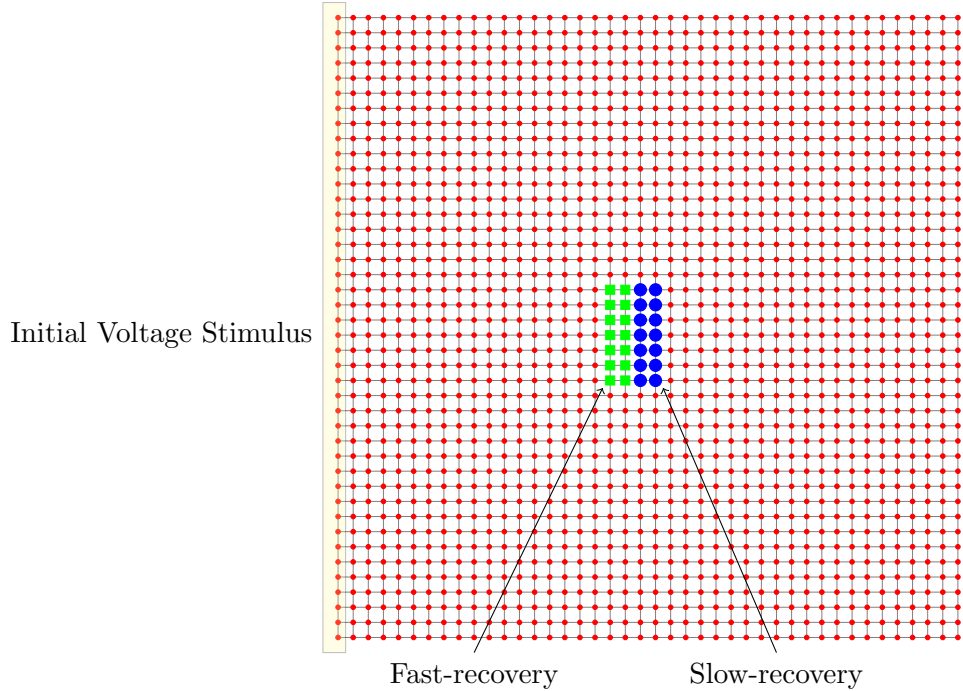
We use a 41 cell by 41 cell grid with a group of fast-recovery cells sized 2 wide by 7 tall, 
and an adjacent equally-sized group of slow-recovery cells, as illustrated in \cref{fig:simbasicgrid}.

There is no coupling outside the grid, so effectively the boundaries are insulating.
We use the parameters in \cref{tab:parameters}, with $\epsilon=0.0065$ for the fast-recovery cells and $b=5.3$ for the slow-recovery cells.
The horizontal and vertical conductances are $g_{\mathrm{horz}} = 0.027$ and $g_{\mathrm{vert}} = 0.015$.
The simulation begins with a voltage stimulus applied to the first column of cells. 
The fourth-order Runge-Kutta method is used to approximate \cref{eqn:FHNcoupled}, with time step $\delta t=0.1$ units. 

\begin{figure}
    \begin{tikzpicture} 
    \setlength{\figscale}{0.23\textwidth}
    \matrix[row sep=1, column sep=0, cells={scale=1}] at (0,0)
        {
        \node[anchor=south] at (0,0) {\includegraphics[width=\figscale,height=\figscale]{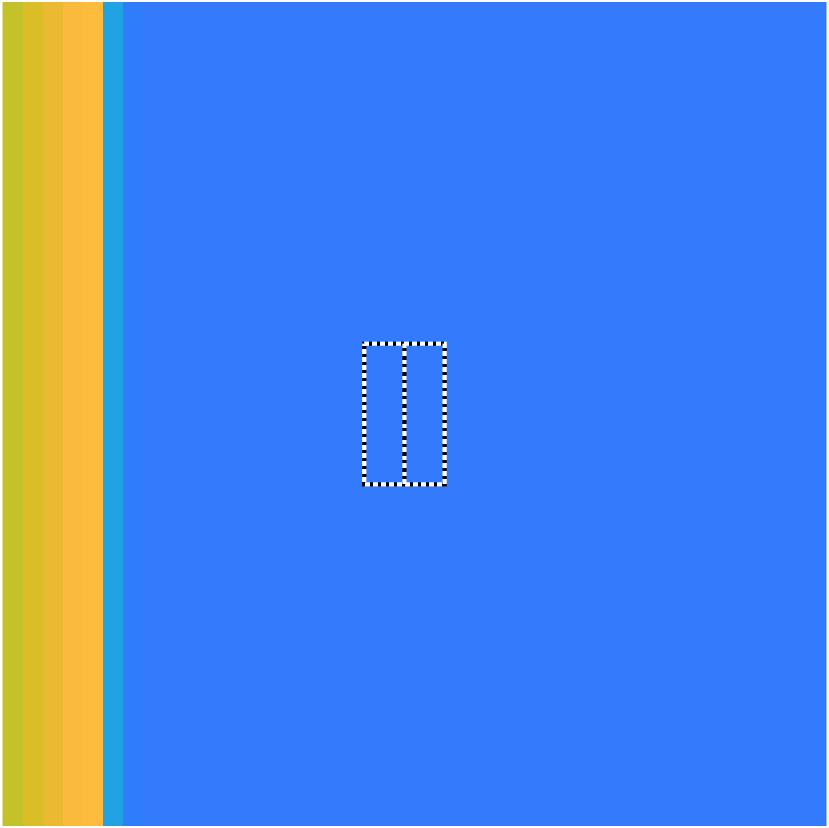}};
        \node[anchor=north] at (0,0) {$t=50$};
        &
        \node[anchor=south] at (0,0) {\includegraphics[width=\figscale,height=\figscale]{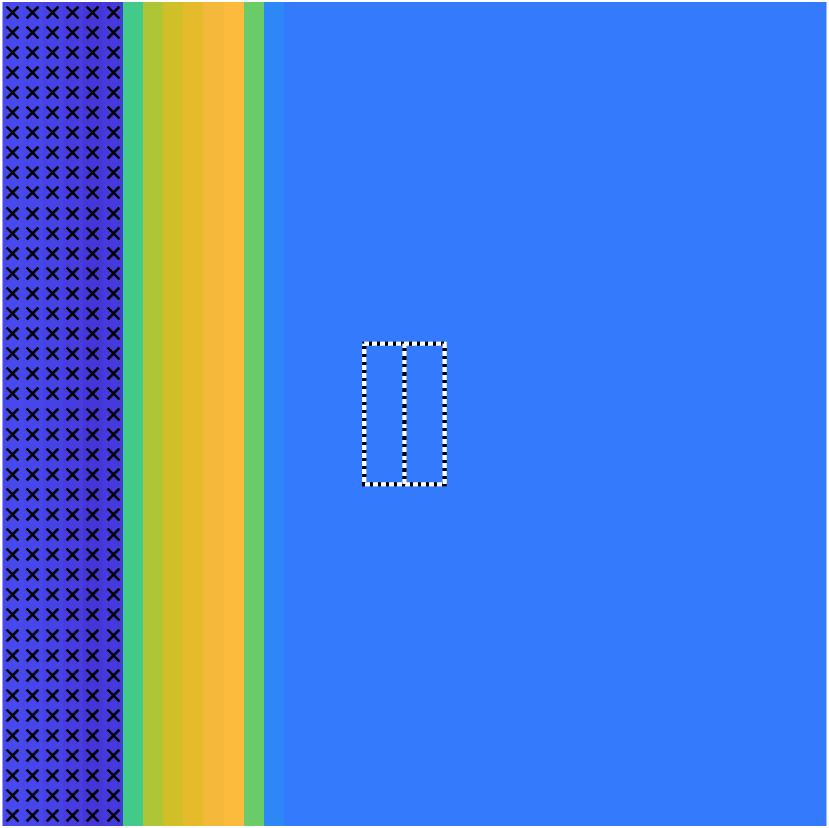}};
        \node[anchor=north] at (0,0) {$t=130$};
        &
        \node[anchor=south] at (0,0) {\includegraphics[width=\figscale,height=\figscale]{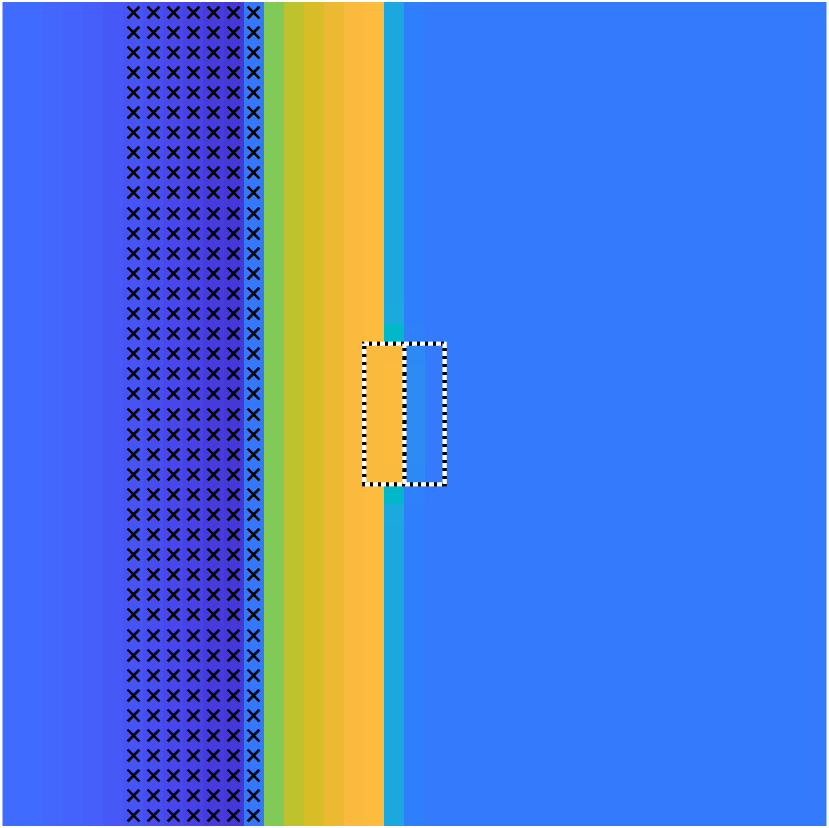}};
        \node[anchor=north] at (0,0) {$t=200$};
        &
        \node[anchor=south] at (0,0) {\includegraphics[width=\figscale,height=\figscale]{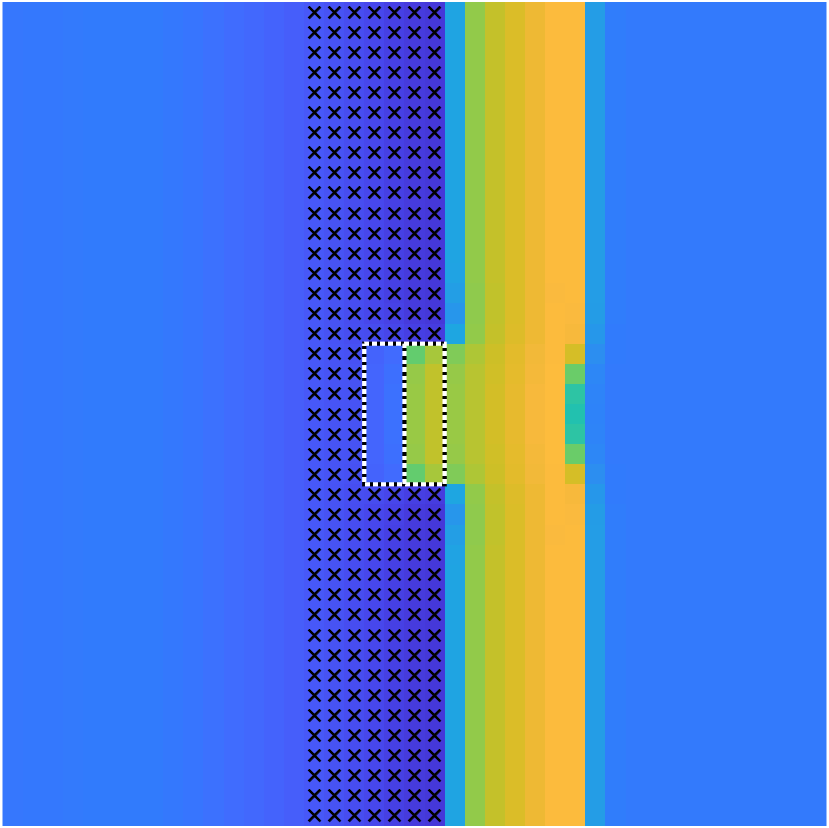}};
        \node[anchor=north] at (0,0) {$t=305$};
        \\
        \node[anchor=south] at (0,0) {\includegraphics[width=\figscale,height=\figscale]{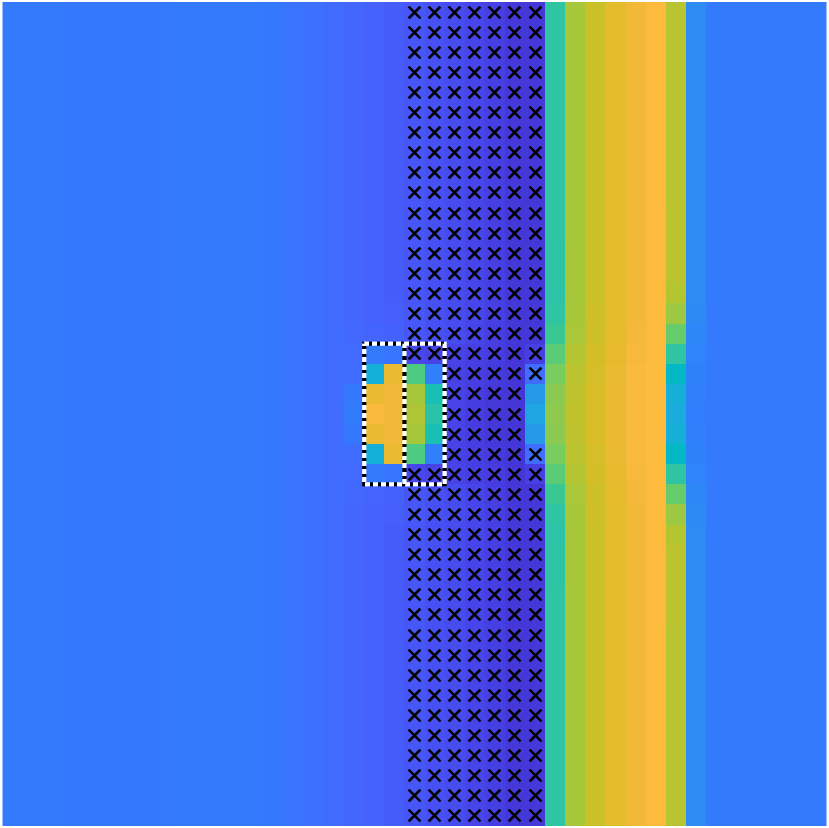}};
        \node[anchor=north] at (0,0) {$t=355$};
        &
        \node[anchor=south] at (0,0) {\includegraphics[width=\figscale,height=\figscale]{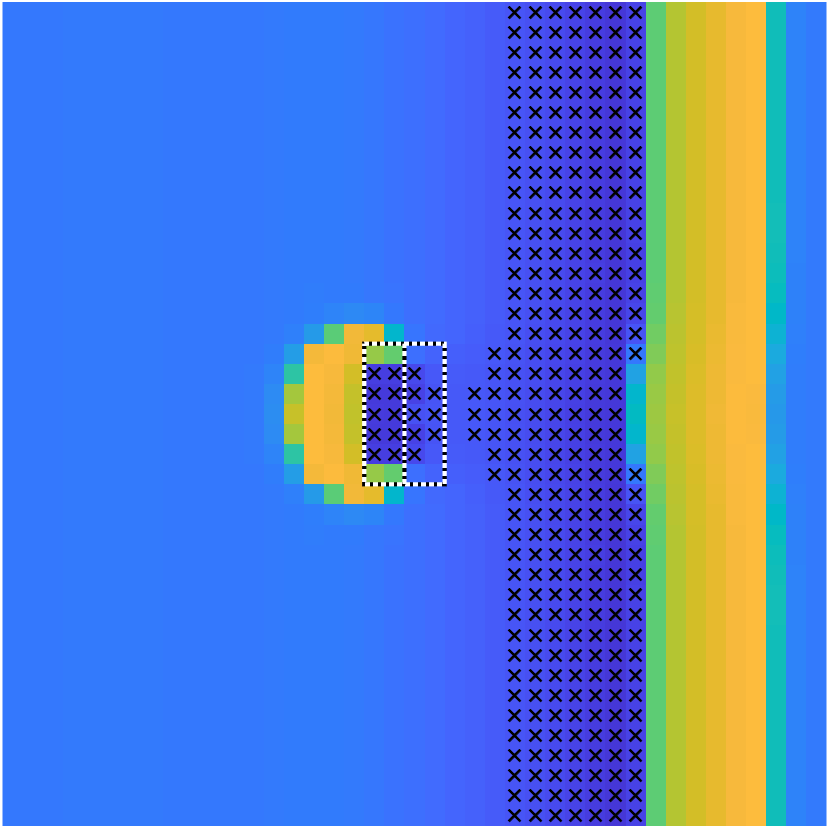}};
        \node[anchor=north] at (0,0) {$t=405$};
        &
        \node[anchor=south] at (0,0) {\includegraphics[width=\figscale,height=\figscale]{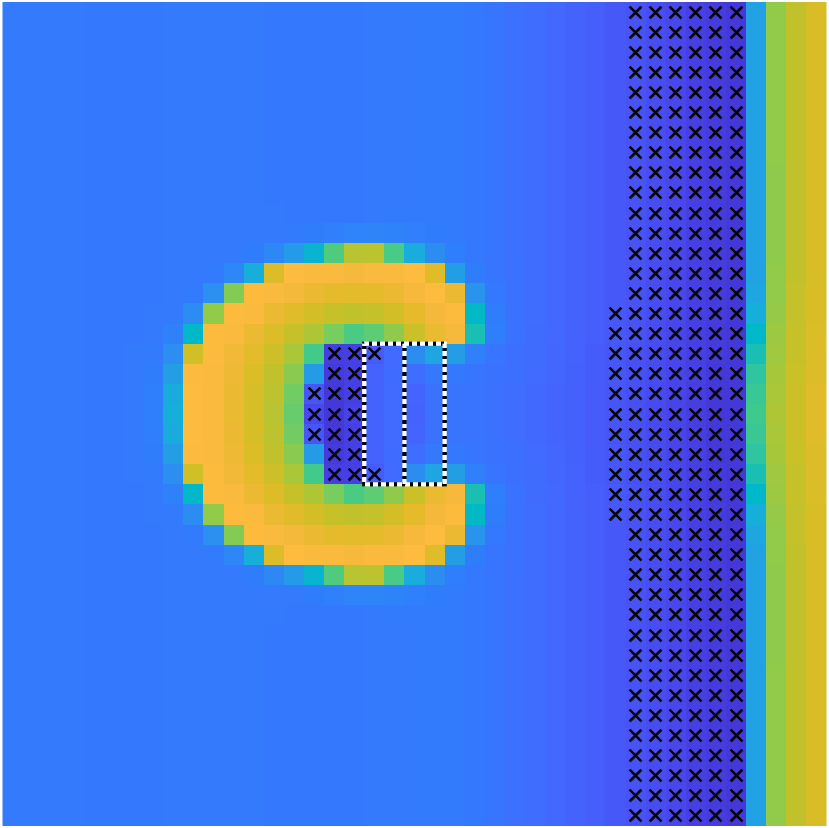}};
        \node[anchor=north] at (0,0) {$t=465$};
        &
        \node[anchor=south] at (0,0) {\includegraphics[width=\figscale,height=\figscale]{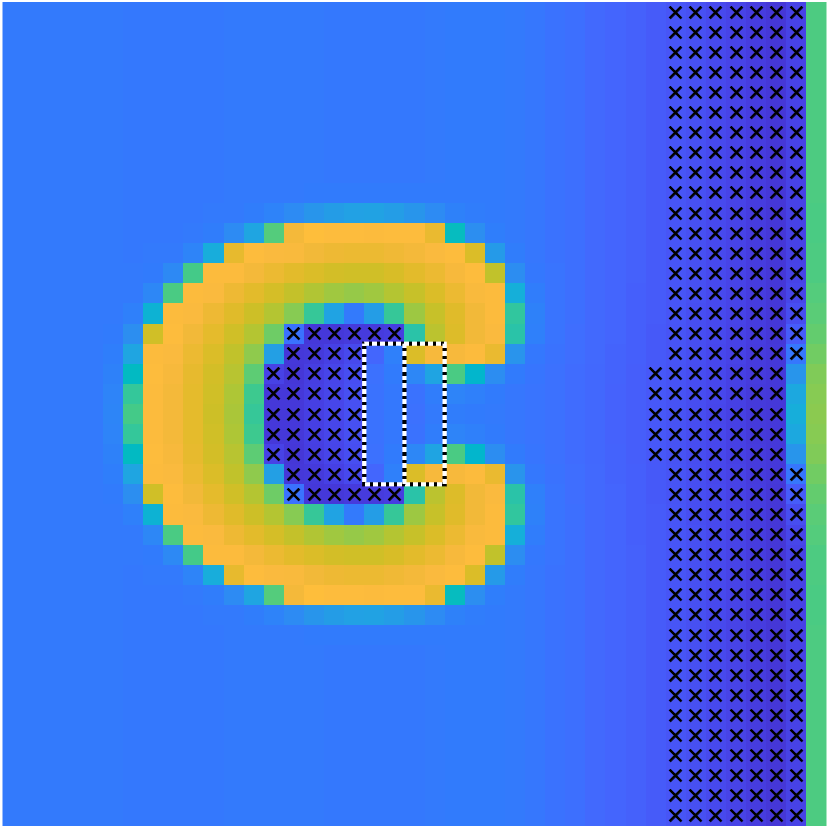}};
        \node[anchor=north] at (0,0) {$t=490$};
        \\
        \node[anchor=south] at (0,0) {\includegraphics[width=\figscale,height=\figscale]{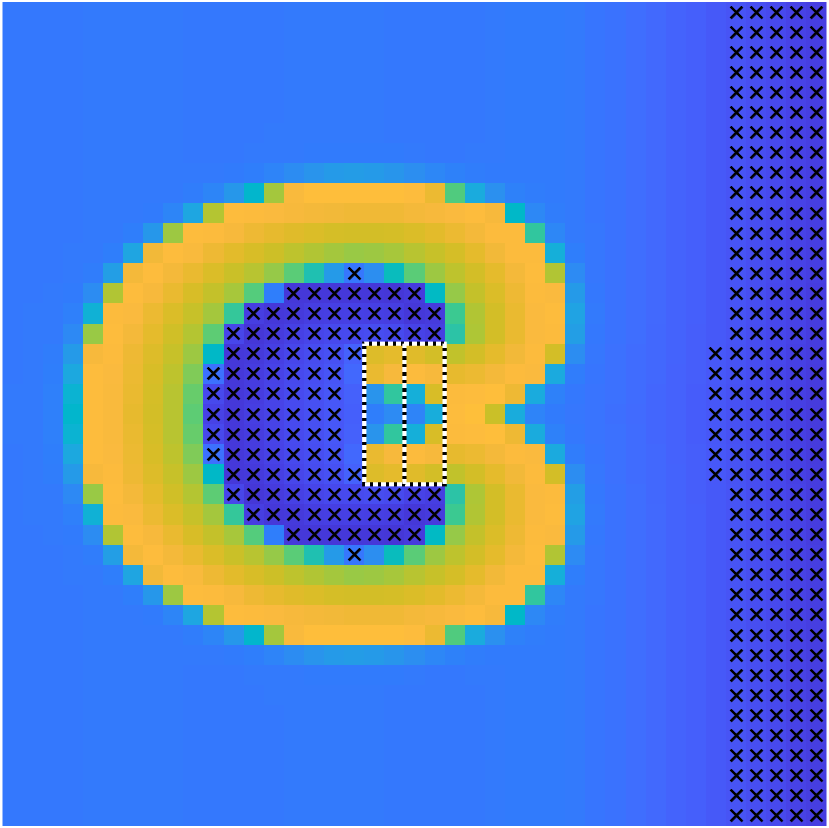}};
        \node[anchor=north] at (0,0) {$t=520$};
        &
        \node[anchor=south] at (0,0) {\includegraphics[width=\figscale,height=\figscale]{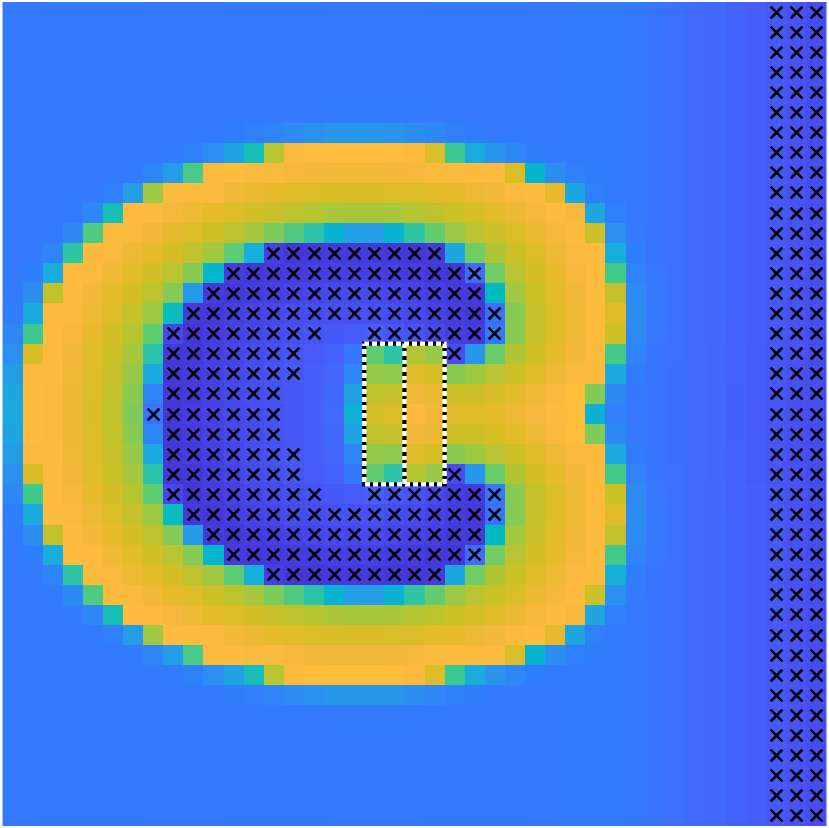}};
        \node[anchor=north] at (0,0) {$t=550$};
        &
        \node[anchor=south] at (0,0) {\includegraphics[width=\figscale,height=\figscale]{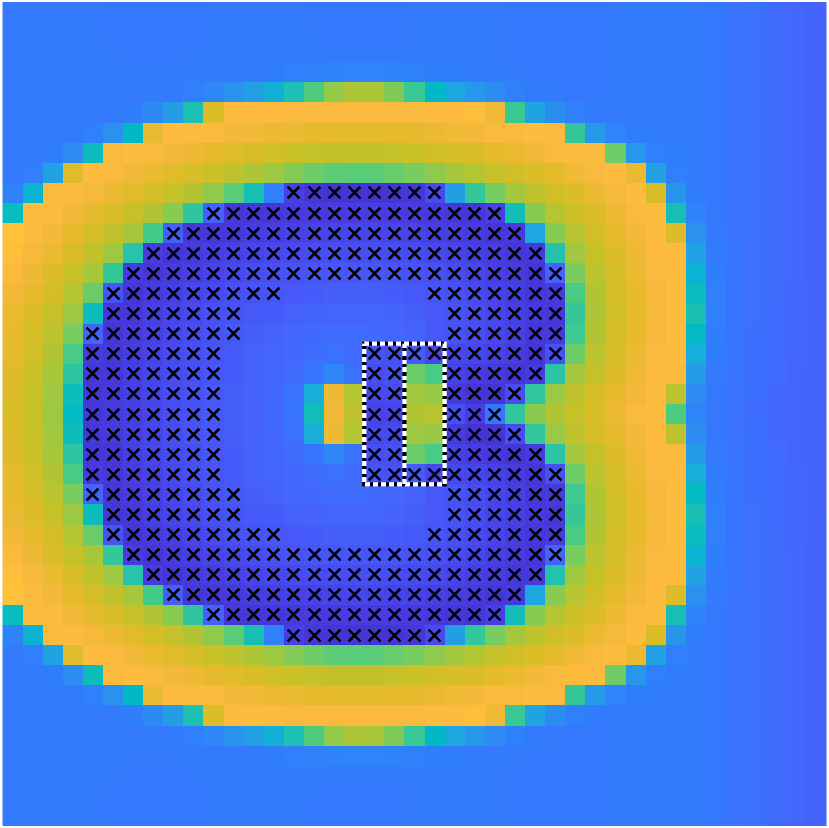}};
        \node[anchor=north] at (0,0) {$t=590$};
        &
        \node[anchor=south] at (0,0) {\includegraphics[width=\figscale,height=\figscale]{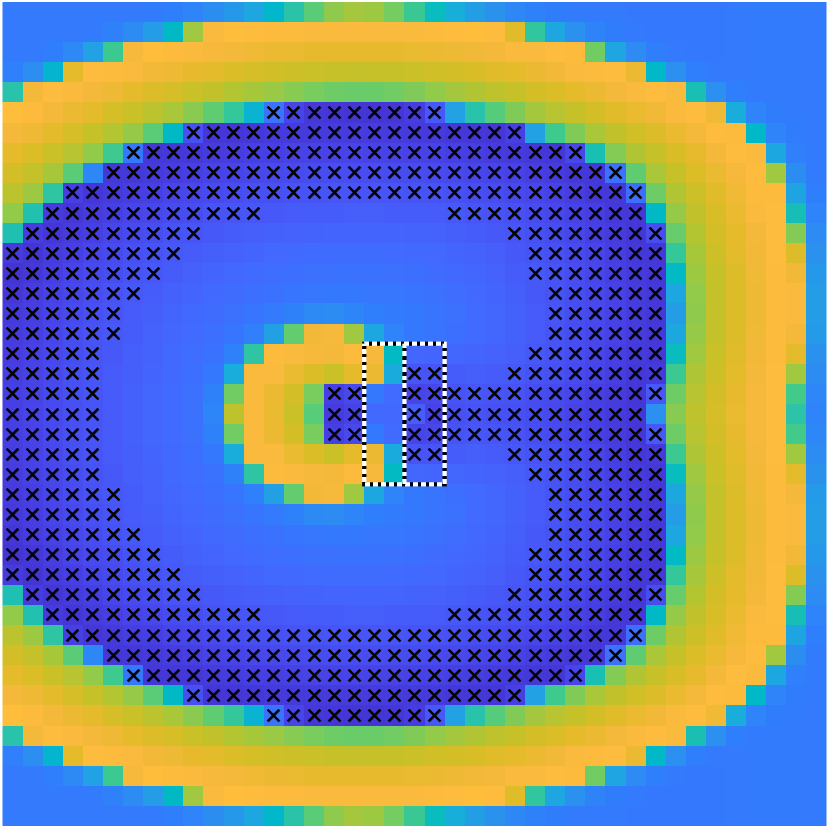}};
        \node[anchor=north] at (0,0) {$t=650$};
        \\
        \node[anchor=south] at (0,0) {\includegraphics[width=\figscale,height=\figscale]{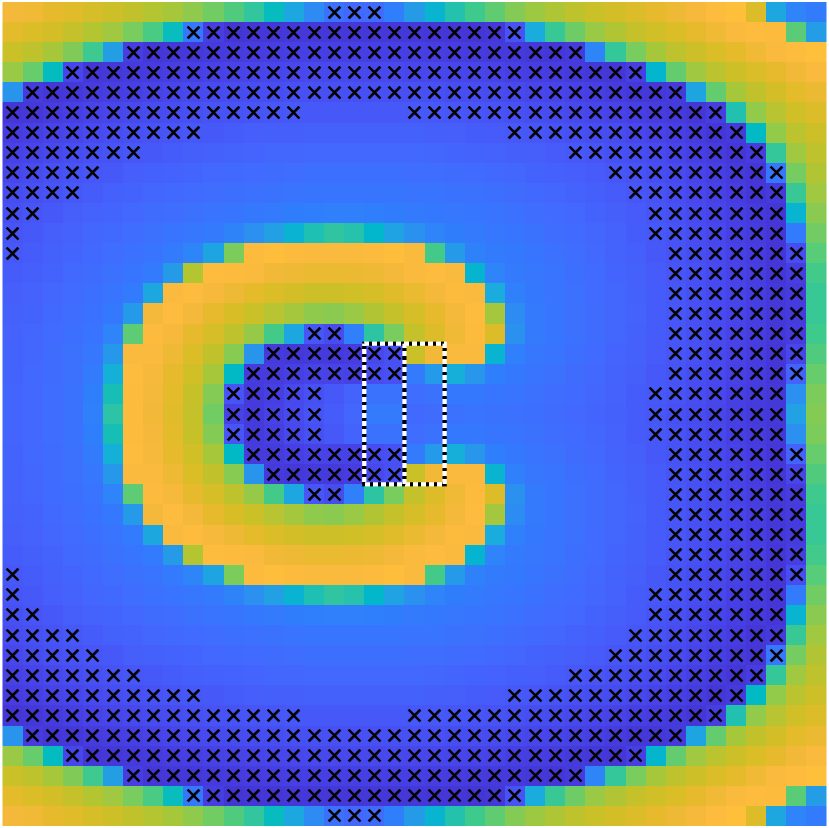}};
        \node[anchor=north] at (0,0) {$t=720$};
        &
        \node[anchor=south] at (0,0) {\includegraphics[width=\figscale,height=\figscale]{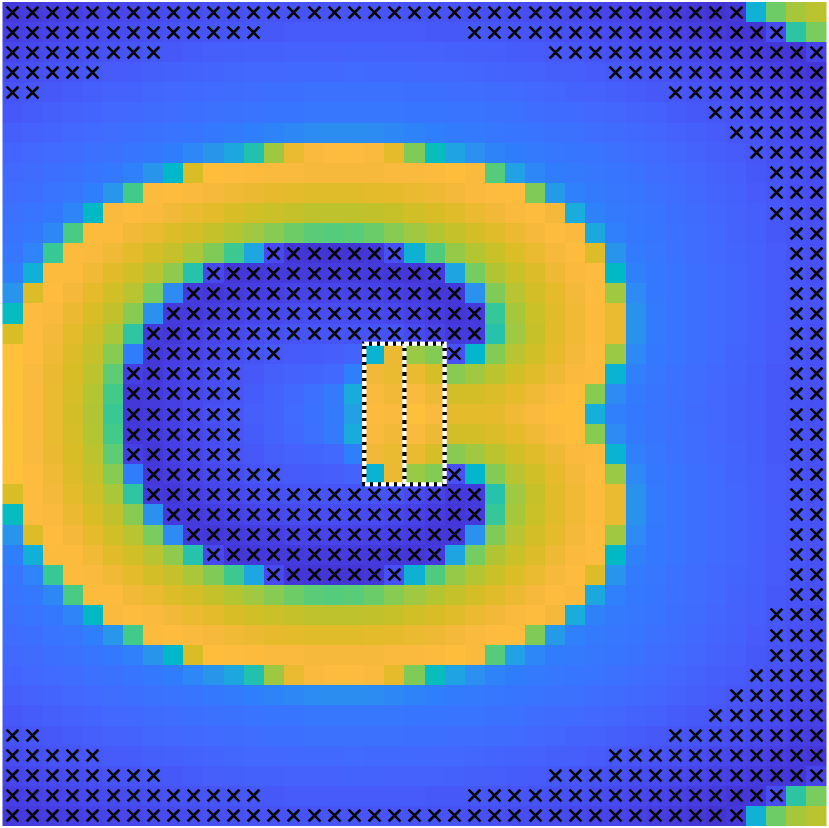}};
        \node[anchor=north] at (0,0) {$t=785$};
        &
        \node[anchor=south] at (0,0) {\includegraphics[width=\figscale,height=\figscale]{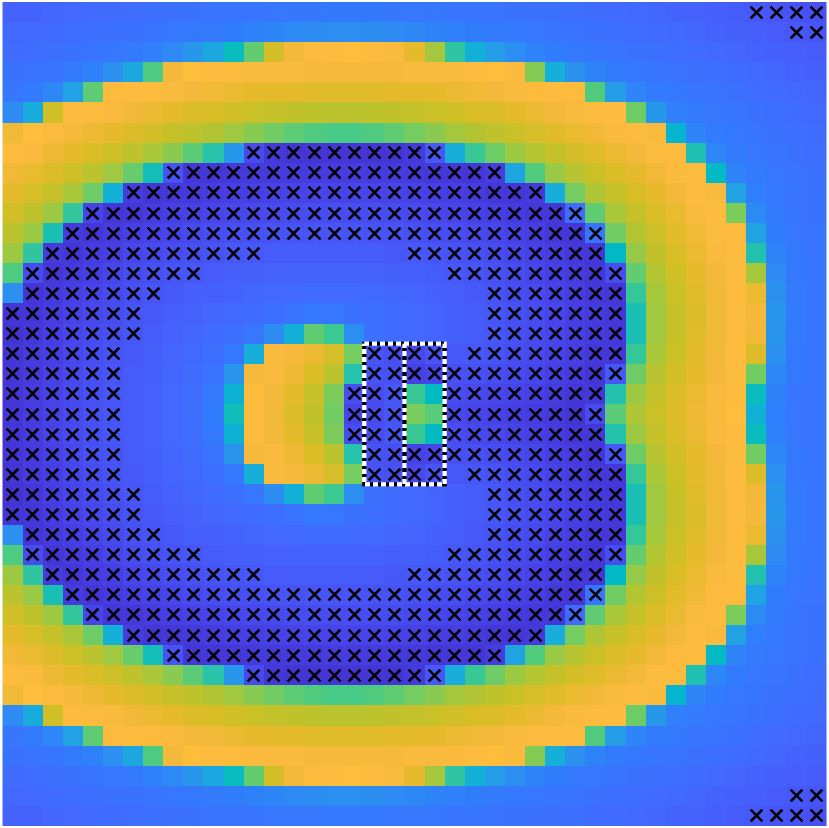}};
        \node[anchor=north] at (0,0) {$t=860$};
        &
        \node[anchor=south] at (0,0) {\includegraphics[width=\figscale,height=\figscale]{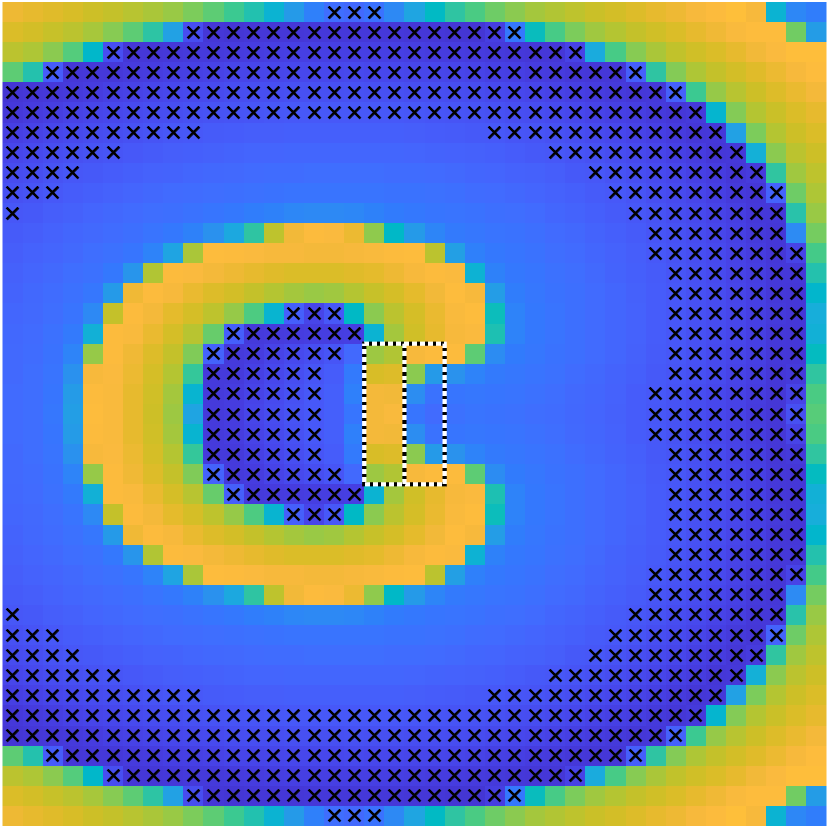}};
        \node[anchor=north] at (0,0) {$t=955$};
        \\
        };
    \end{tikzpicture}
    \caption{Snapshots from a simulation corresponding to \cref{fig:basicsequence}, using the computational grid in \cref{fig:simbasicgrid}, with the layers of unhealthy cells marked by dashed rectangles.
    Cells are colored by voltage, from blue (low) to yellow (high).
    Cells with $v<0$ and $w>0.0225$ are marked with $\times$ to indicate that they are in recovery and inhibited from firing.
    Note that the dynamics are not precisely periodic.  At $t = 490$ the wave re-stimulates the slow-recovery cells first, while at $t = 955$ the wave re-stimulates the fast-recovery layer before reaching the slow-recovery layer.
    }
    \label{fig:basicsnaps}
\end{figure}
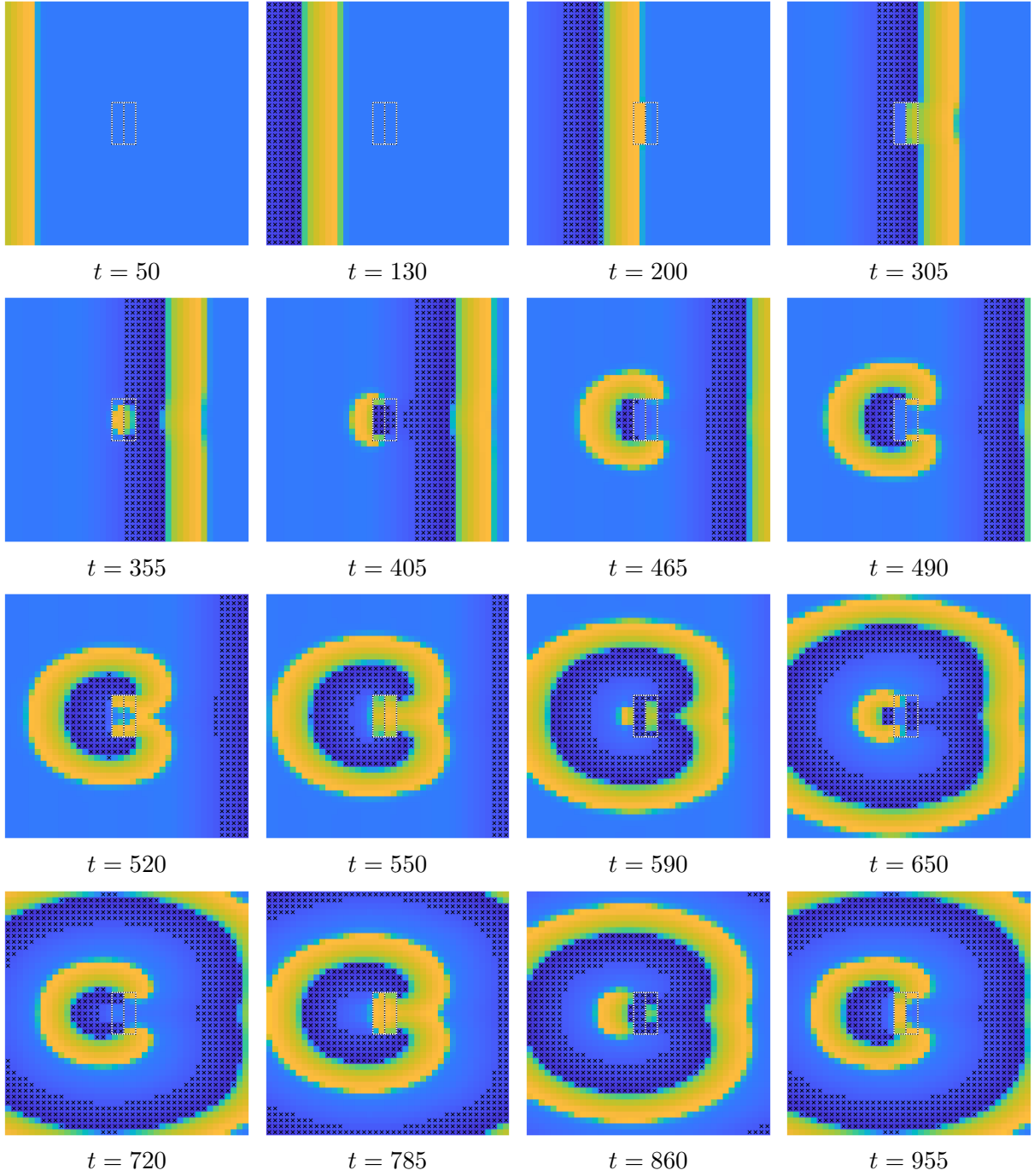

In \cref{fig:basicsnaps} we give snapshots from the simulation corresponding to the cartoon sequence of events in \cref{fig:basicsequence}.
A video is included in the Supplemental Material \cite{AF-supplement}.
The events occur as follows.
\begin{enumerate}
    \item A normal  impulse wave arrives, triggering all cells to fire. Once this pulse has passed the layers of unhealthy cells, the three types of cells are in different states:
    \label{item:ititialpulse}
    \begin{itemize}
        \item The slow-recovery cells still have high voltage.
        \item The fast-recovery cells have recovered and are ready to fire again.
        \item The healthy cells adjacent to the slow-recovery cells have not yet recovered and so cannot fire again.
    \end{itemize}
    The slow-recovery cells in the middle of their layer retain the highest voltage, because the other slow-recovery cells above and below them help insulate them from charge leakage.  
    \item The slow-recovery cells stimulate the fast-recovery cells by the plateau (phase 2) mechanism,  causing the fast-recovery cells to fire a second time, which triggers a second wave in the reverse direction.
    \label{itm:fastfireagain}
    \begin{itemize}
    \item The fast-recovery cells in the middle of their layer fire first, because the slow-recovery cells in the middle of their layer have retained the highest voltage.
        \item The healthy cells adjacent to the fast-recovery cells have now recovered, and thus fire.
    \end{itemize}
\end{enumerate}
In \cref{item:ititialpulse}, it is important that the slow-recovery cells retain enough voltage for long enough to cause the fast-recovery cells to fire again. 
In \cref{sec:parameterinteractions} we observe constraints on the parameter combinations that allow this to happen and in \cref{sec:localtiming} we validate these constraints using properties of the fast-recovery and slow-recovery cells in isolation.
Having adjacent slow-recovery cells and fast-recovery cells with compatible parameter combinations is the essential assumption of our model.

\begin{enumerate}
\setcounter{enumi}{2}
\item \label{itm:curlback} The second wave travels through the healthy cells. 
\begin{itemize}
        \item The broken ends of the wave propagate around the layers of fast-recovery and slow-recovery cells, which have finished firing but have not yet recovered, forming swirling-back curves.
        \item By the time the wave has propagated around to the middle of the right side of the slow-recovery cells layer, the slow-recovery cells are ready to fire again.
\end{itemize}    
\end{enumerate}

In \cref{itm:curlback}, the geometry required the wave through the healthy cells take enough time that the slow-recovery cells could recover and be ready to fire again.
In \cref{sec:variationgeometry} we observe the effects of moderate variations in the geometry and conclude that they affect the details of the phenomena but are not essential to the model.
In \cref{sec:exoticsimulations} we consider the effects of more exotic geometries.


\begin{enumerate}
\setcounter{enumi}{3}
    \item The slow-recovery cells fire a second time, which triggers the fast-recovery cells to fire a third time, initiating a third wave.
\end{enumerate}
The third wave is different than the second wave, and did not depend on compatible parameter combinations between the slow-recovery and fast-recovery cells.

\begin{enumerate}
\setcounter{enumi}{4}
    \item The third wave travel around in swirling-back curves and initiates a fourth wave.
\end{enumerate}
The fourth wave is different than the second and third waves.
From this point, the sequence of events becomes somewhat irregular. 
By $t=2500$, there had been 8 waves and at $t=20000$ waves were still being generated.

Various secondary timing and geometric effects also occur. 
For example, if the slow-recovery cells still have high voltage when the second wave reaches them, then they may be pushed back to higher voltage and never recover.
This state is non-physiologic in that such cells would eventually die.
In \cref{sec:parameterinteractions} we observe the parameters that produce this state and in \cref{sec:localnonphysslow} we validate these observations using a two-cell model coupling one fast-recovery cell and one slow-recovery cell.


\section{Interaction of the Fast-recovery and Slow-recovery Governing Parameters}
\label{sec:parameterinteractions}
\label{sec:local}

In this section we study how the slow-recovery cell parameter $b$ and the fast-recovery cell parameter $\epsilon$ jointly determine if wave formation is sustained in simulations using the computational grid in \cref{fig:simbasicgrid}.
We designate the slow-recovery cell that is vertically centered and adjacent to the fast-recovery cell block as the \emph{sensing cell} and monitor what it does.
We count its number of beats, measured by the number of times its voltage increases above $0.4$ and subsequently decreases below 0.
We classify its state at the end of the simulation as either recovered, still firing, or stuck at high voltage, based on its behavior in the last 300 time units of the 2500 time unit simulation.
If its maximum voltage stays below $0.4$, then it is labeled recovered.
If its minimum is above $0.4$, then it is labeled stuck high.
If the number of beats within those last 300 time units is greater than one, it is labeled still firing.
This classification scheme does not necessarily include all possible voltage traces, but our simulations have not produced a $(b,\epsilon)$ pair that does not fall under one of these three labels.

\begin{figure}
    \begin{tabular}{cc}
         Beat Count& End State  \\
         \includegraphics[width=0.49\textwidth]{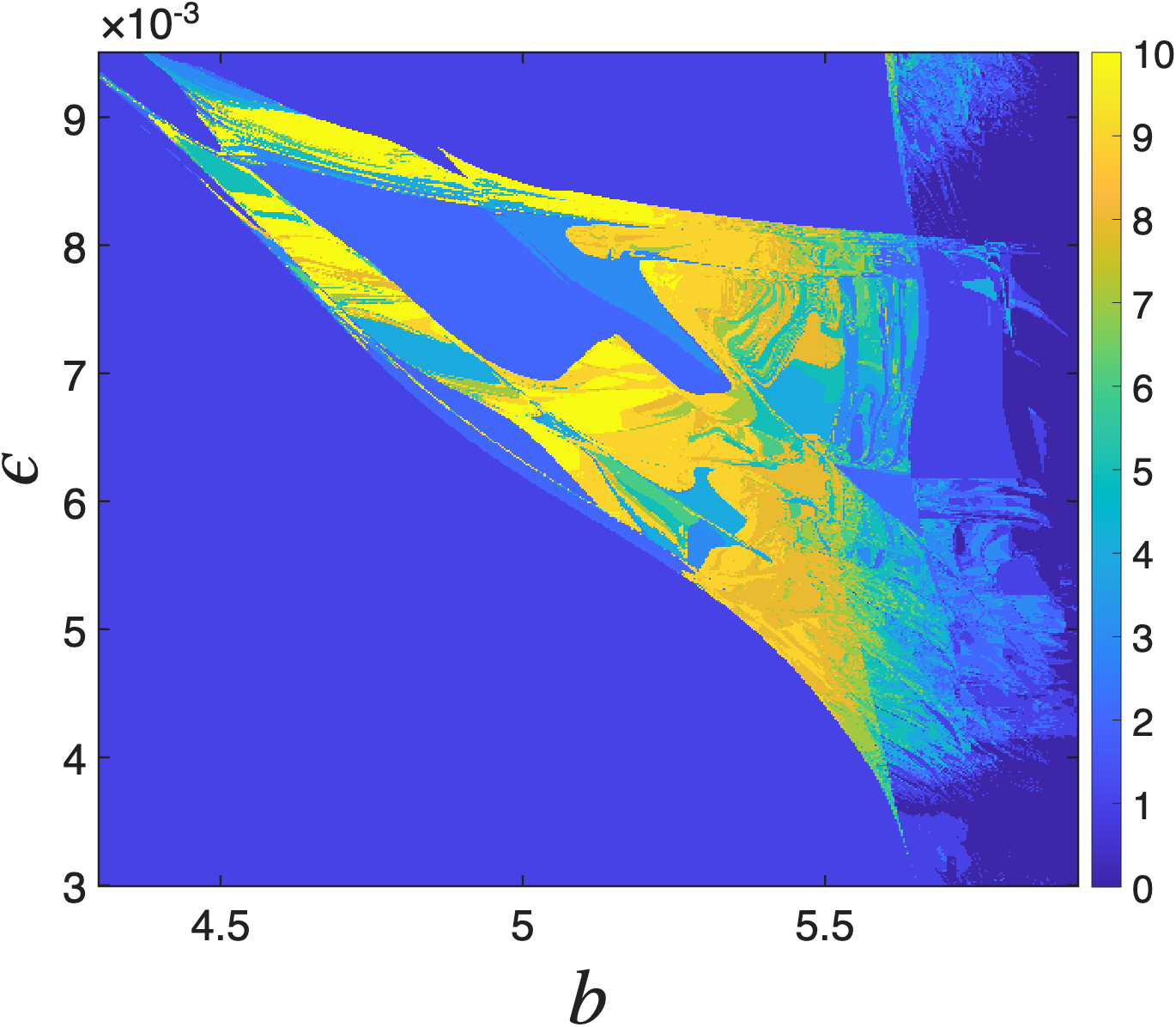} & \includegraphics[width=0.476\textwidth]{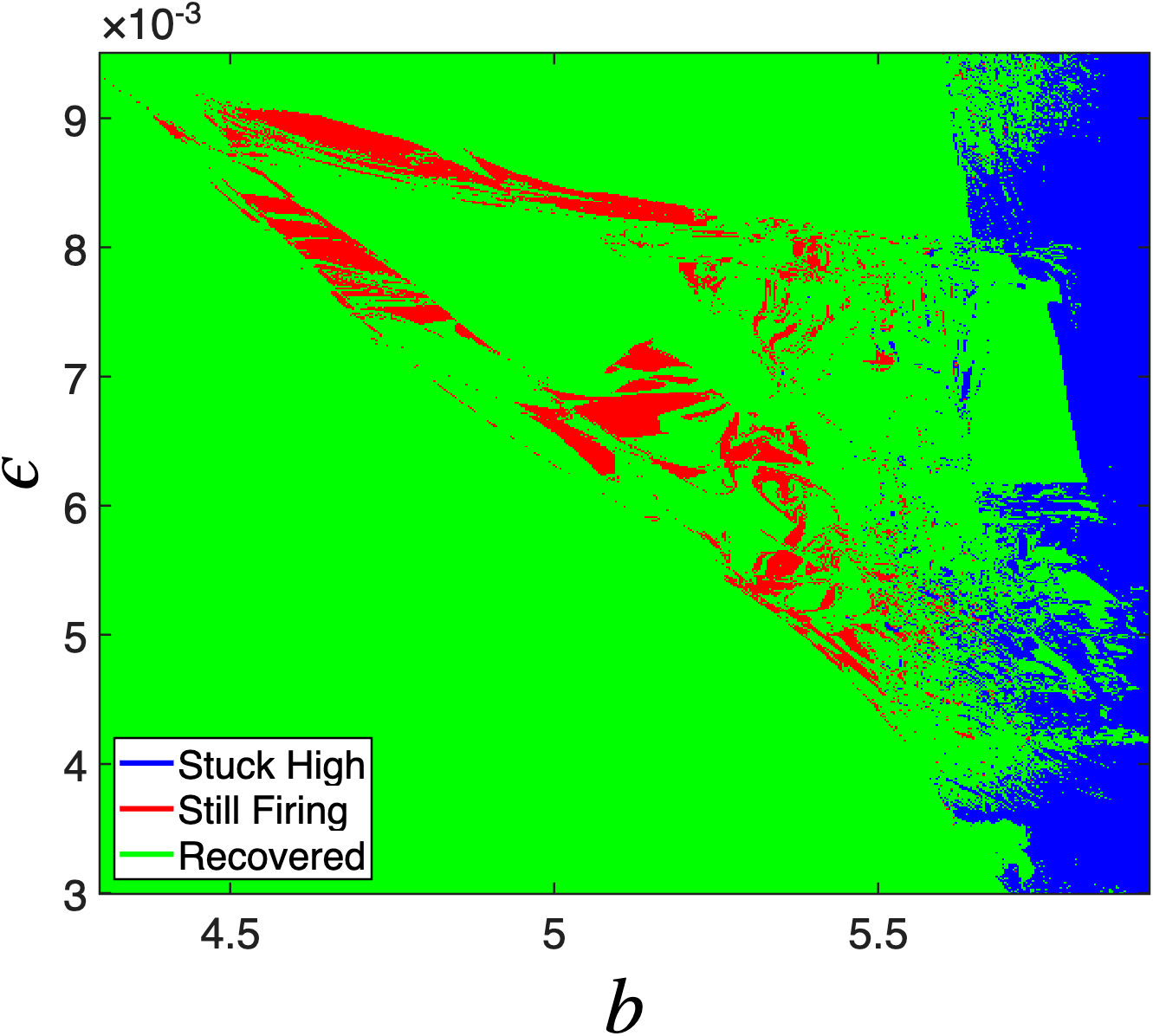}
    \end{tabular}
    \caption{ Beat count (left) and end state (right) of the sensing cell as a function of $(b, \epsilon)$ in a simulation on the grid of cells in \cref{fig:simbasicgrid} evolved for 2500 time units. 
         There is a well-defined
        diagonal line marking the onset of multiple firings; this is a local
        effect that can be explained in terms of two one-cell models. For $b \gtrsim 5.6$, 
        slow-recovery cells become stuck in a high-voltage state and never recover. 
    }
    \label{fig:OSC 7 FIB}
\end{figure}

In \cref{fig:OSC 7 FIB}, we show the number of beats and end state for physiologically relevant intervals of both $b$ and $\epsilon$.
Two notable features are apparent. 
First is a diagonal boundary, to the lower-left of which there is only a normal beat, but to the upper-right of which the sensing cell fires multiple times.
In \cref{sec:localtiming}, we show that this boundary is explained by the relative recovery times of the two types of unhealthy cells, which determine whether the slow-recovery cells can initiate the second wave.
Second is a vertical boundary around $b = 5.7$, to the right of which the slow-recovery cells enter a non-physiologic state of maintaining high voltage, as explained in \cref{sec:localnonphysslow}.


\subsection{Slow-recovery cells charge retention time interacting with fast-recovery cells recovery time}
\label{sec:localtiming}
In this section, we examine the key element required to create FIBs: the time differential between how long the slow-recovery cells retain voltage and how soon the fast-recovery cells are prepared to fire after the initial impulse has occurred. 
We revert to the one-cell model for this computation to ensure that no geometric effects are included.
In \cref{sec:variationgeometry} we will consider geometric effects, such as how the types of cells surrounding a slow-recovery cell affect how long it retains voltage. 

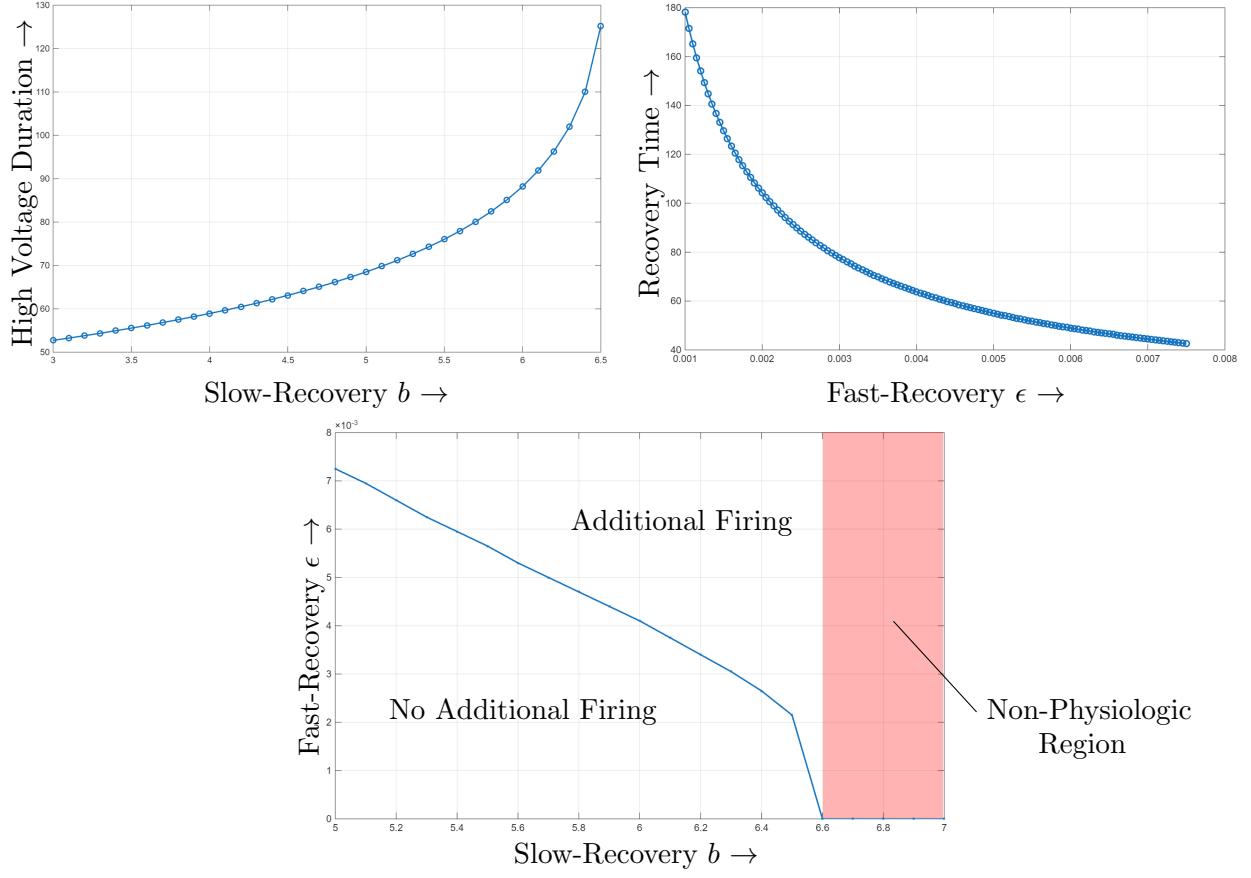
\begin{figure}
    \centering
    \begin{tikzpicture}
    \setlength{\figscale}{0.46\textwidth}
        \node[anchor=east] at (-.25,0) {\includegraphics[width=\figscale]{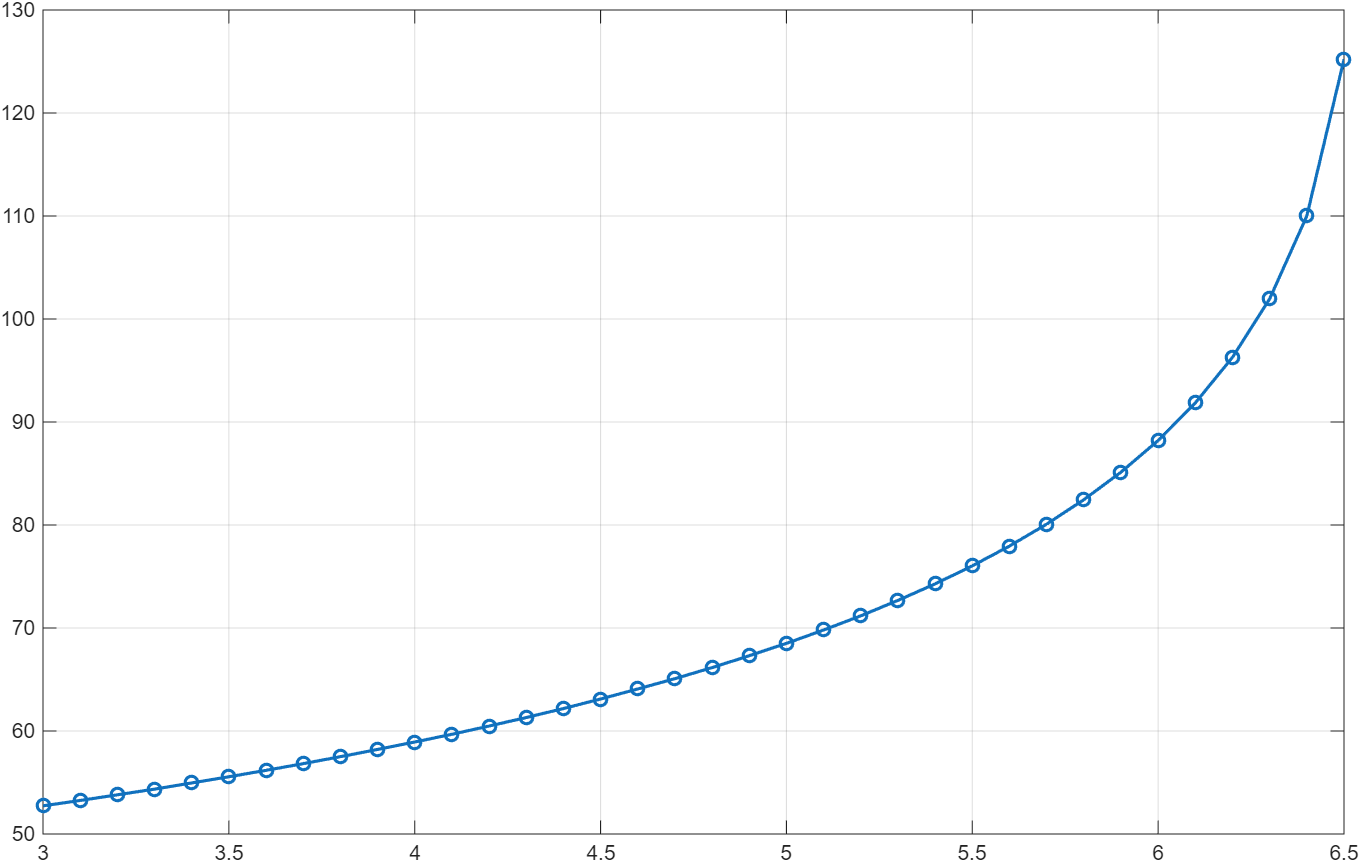}};
        \node[anchor=west] at (.25,0) {\includegraphics[width=\figscale]{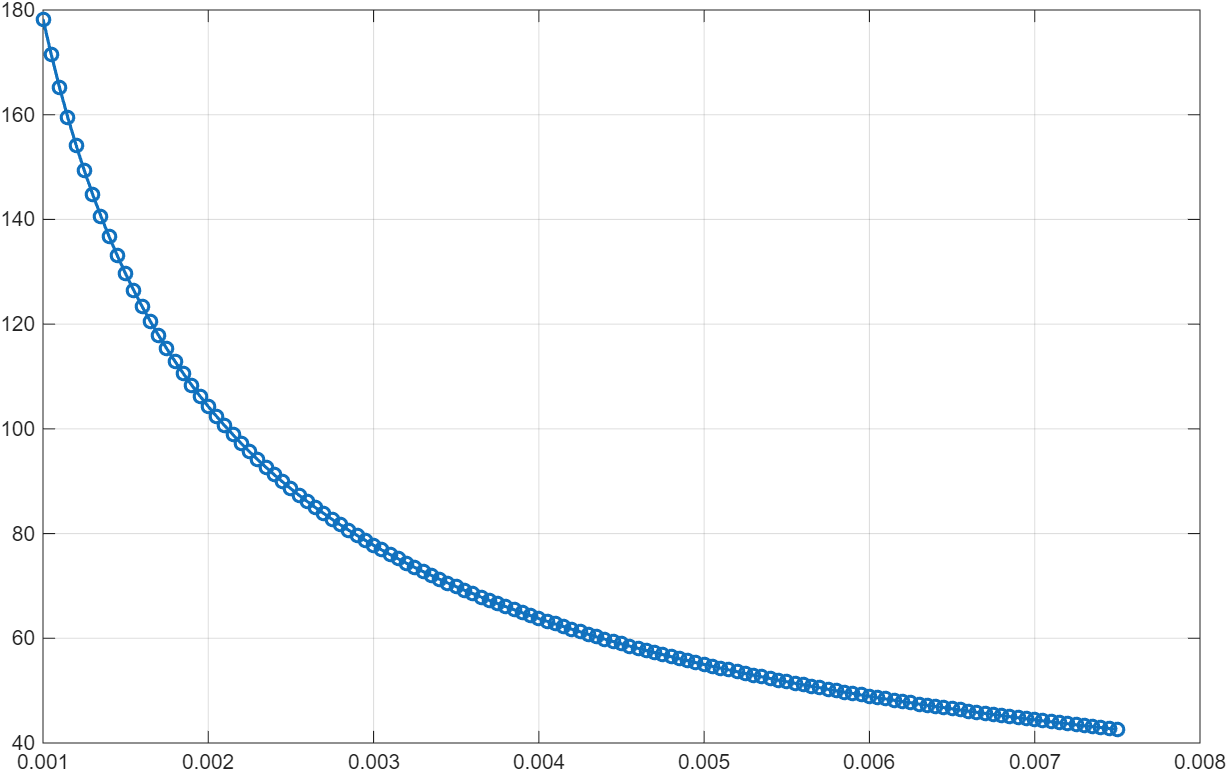}};
        \node[anchor=north] at (0,-3) {\includegraphics[width=.5\textwidth]{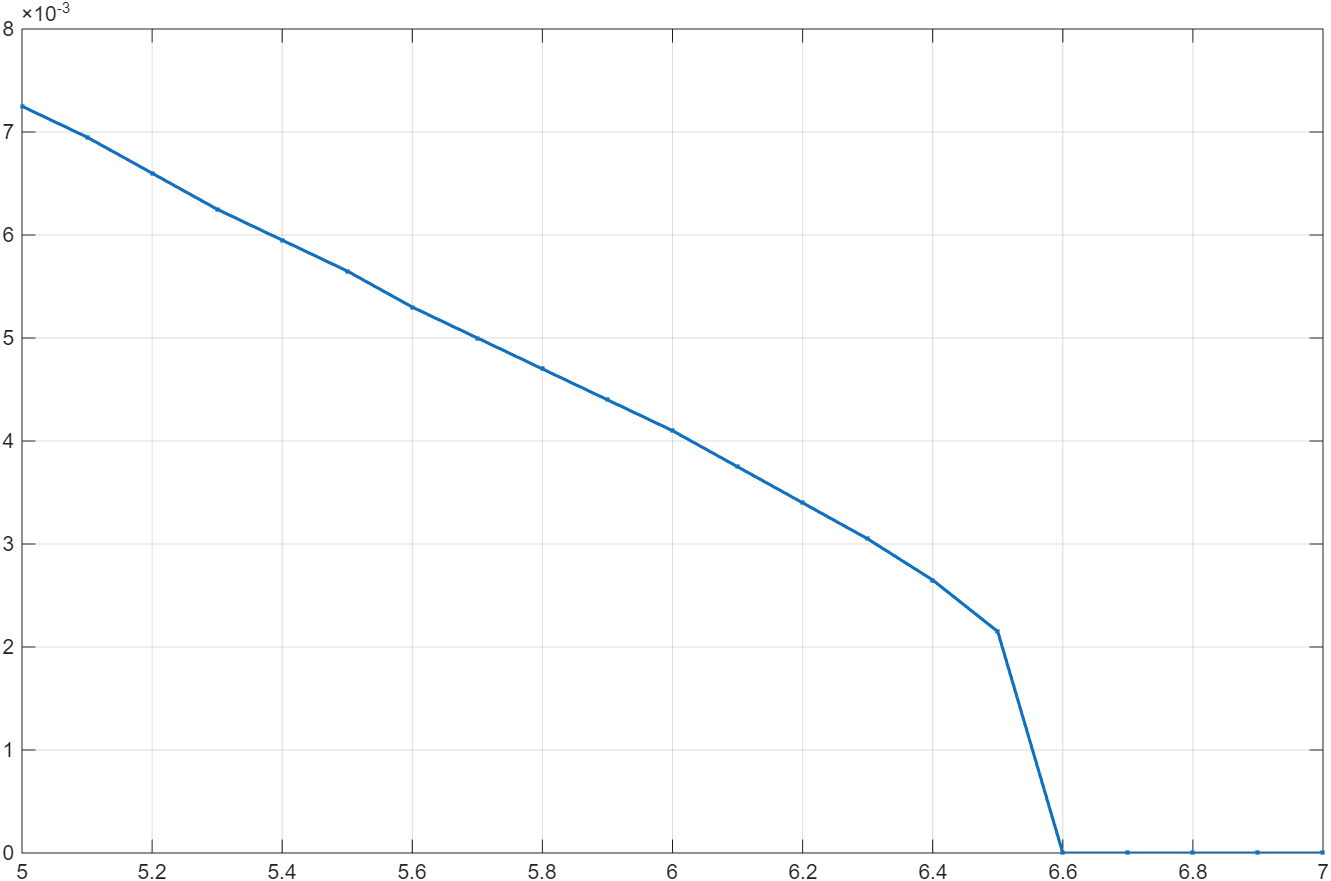}};

        \node[anchor=north] at (-4.1,-2.5) {Slow-Recovery $b$ $\rightarrow$};
        \node[anchor=north, rotate=90] at (-8.4,0) {High Voltage Duration $\rightarrow$};

        \node[anchor=north] at (4.1,-2.5) {Fast-Recovery $\epsilon$ $\rightarrow$};
        \node[anchor=north, rotate=90] at (-.1,0) {Recovery Time $\rightarrow$};
        
        \node at (.6,-4.5) {Additional Firing};
        \node at (-1.5,-7) {No Additional Firing};
        \fill[red, opacity=0.3] (2.47,-8.4) rectangle (4.05,-3.3);
        \draw[-] (4.5,-7) node[anchor=west] {Non-Physiologic} -- (3.4,-5.8);
        \node at (5.9,-7.45) {Region};
        \node[anchor=north] at (0,-8.6) {Slow-Recovery $b$ $\rightarrow$};
        \node[anchor=north, rotate=90] at (-4.6,-6) {Fast-Recovery $\epsilon$ $\rightarrow$};
    \end{tikzpicture}
    \caption{In the first plot, we show how the duration of high voltage increases as the $b$ value of the slow-recovery cell increases. In the second plot, we show how recovery time decreases as the $\epsilon$ value of the fast-recovery cell increases. 
    In the third plot, we conduct the inverse mapping computation with a 25 time unit buffer to determine what combinations of $b$ and $\epsilon$ would result in an additional firing of the fast-recovery cell. When the $b$ value of the slow-recovery cell exceeds 6.5, it would never recover (and so die), so the third plot assigns these non-physiologic cases with an $\epsilon$ value of 0.}
    \label{fig:timing interraction}
\end{figure}

We measure the duration of high voltage for a single slow-recovery cell by applying a single impulse and finding the difference between the times when the cell voltage rises above $0$ and later falls below $0.5$. 
The value of $0.5$ is somewhat arbitrary, but occurs during the rapid voltage descent phase ((3) in \cref{fig:pp_trace_3}), so the results are not sensitive to this choice.
We find that for $b$ values greater than $6.5$, the slow-recovery cell never repolarizes, and so its voltage never crosses the falling threshold.
The first plot in \cref{fig:timing interraction} shows the duration as a function of $b$.

We measure the recovery time for a single fast-recovery cell in a similar manner, 
though now we calculate the difference between the times when the cell's recovery variable rises above $0$ and when the recovery variable later falls below $0.1$. 
This falling recovery threshold of $0.1$ is selected experimentally to mark the point at which the fast-recovery cell has recovered from the initial impulse sufficiently and is able to be excited by nearby cells.
The second plot in \cref{fig:timing interraction} shows this time as a function of $\epsilon$.


The slow-recovery cells must maintain high voltage for some period beyond the recovery time of the fast-recovery cell in order for enough charge to transfer to trigger the fast-recovery cell.
Using 25 time units as this period, in the third plot in \cref{fig:timing interraction} we show the minimal $\epsilon$ to achieve a trigger, as a function of $b$. 
Above this curve, slow-recovery cells sustain a high voltage long enough to re-activate adjacent fast-recovery cells.
Note that this calculated curve matches well with the diagonal curve in \cref{fig:OSC 7 FIB} that delineates regions with and without multiple firings. 


\subsection{Bifurcations and Stable States of Coupled Fast- and Slow-Recovery Cells}
\label{sec:localnonphysslow}

\begin{figure}
\centering
\begin{tikzpicture}
        \setlength{\figscale}{0.8\textwidth}
        \node[anchor=south west] at (0,0) {{\includegraphics[width=\figscale]{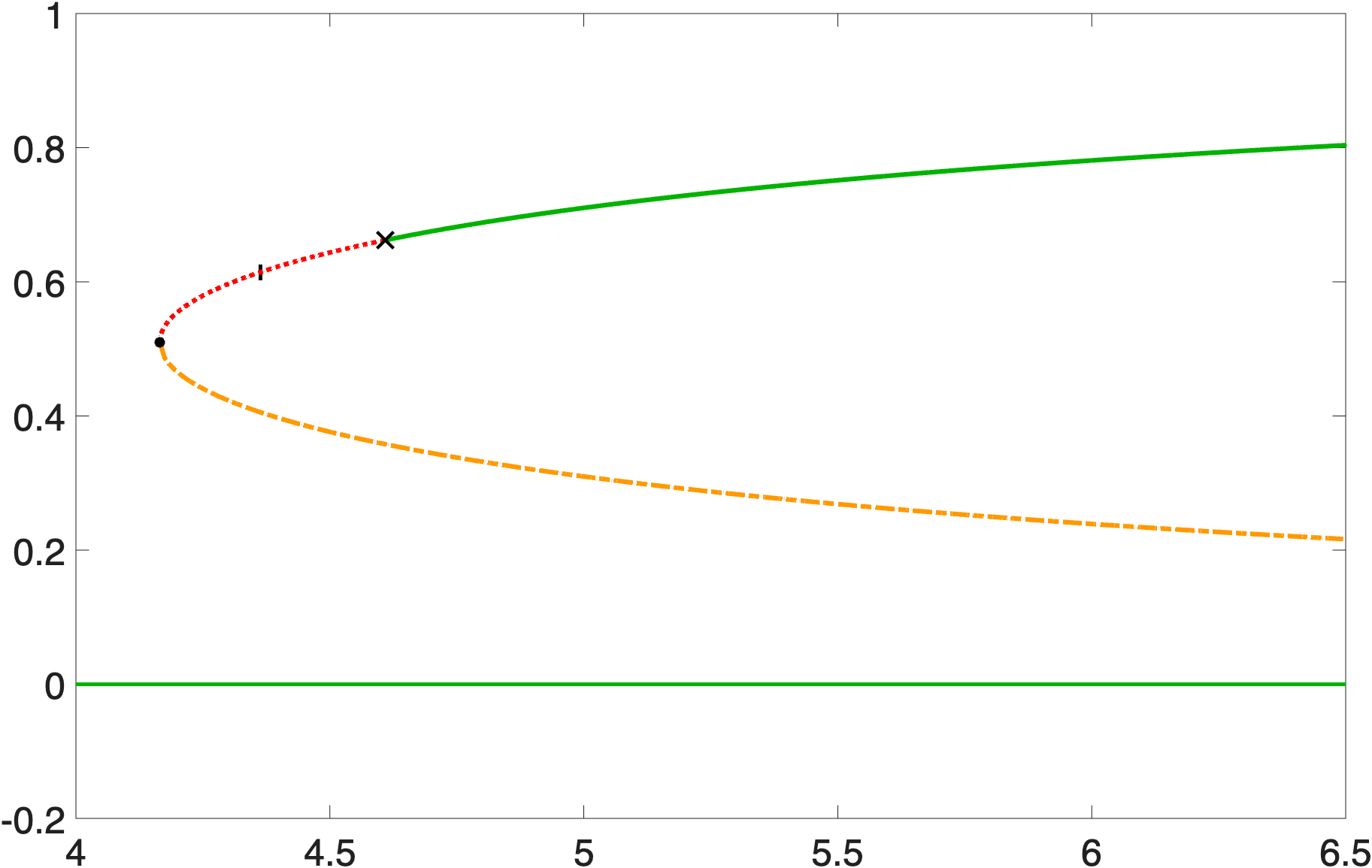}}};
        \node[rotate=90,anchor=south] at (0,5) {voltage $v$};
        \node[anchor=north] at (6.7,0) {$b$};
        \draw[->] (5,1.4) node[anchor=west] {stable equilibrium} -- (4.2,1.8);
        \draw[->] (3,5) node[anchor=west] {saddle-node bifurcation} -- (1.8,5.2);
        \draw[->] (5,7.5) node[anchor=west] {Hopf bifurcation} -- (3.85,6.4);
        \draw[->] (10,6) node[anchor=north] {stable equilibrium} -- (9,6.6);
        \draw[->] (2.5,7) node[anchor=south] {node$\rightarrow$focus} -- (2.6,6.0);
\end{tikzpicture}
\caption{
Bifurcations for a single-cell FHN model \cref{eqn:FHN}. 
At $b \approx 4.1649$,  a saddle-node bifurcation creates a saddle (lower branch)  and unstable node (upper branch). 
At $b \approx 4.3629$ the unstable node becomes an unstable focus.
A Hopf bifurcations occurs at $b \approx 4.6097$, making the upper branch into a stable focus.}
\label{fig:onecell_bif}
\end{figure}
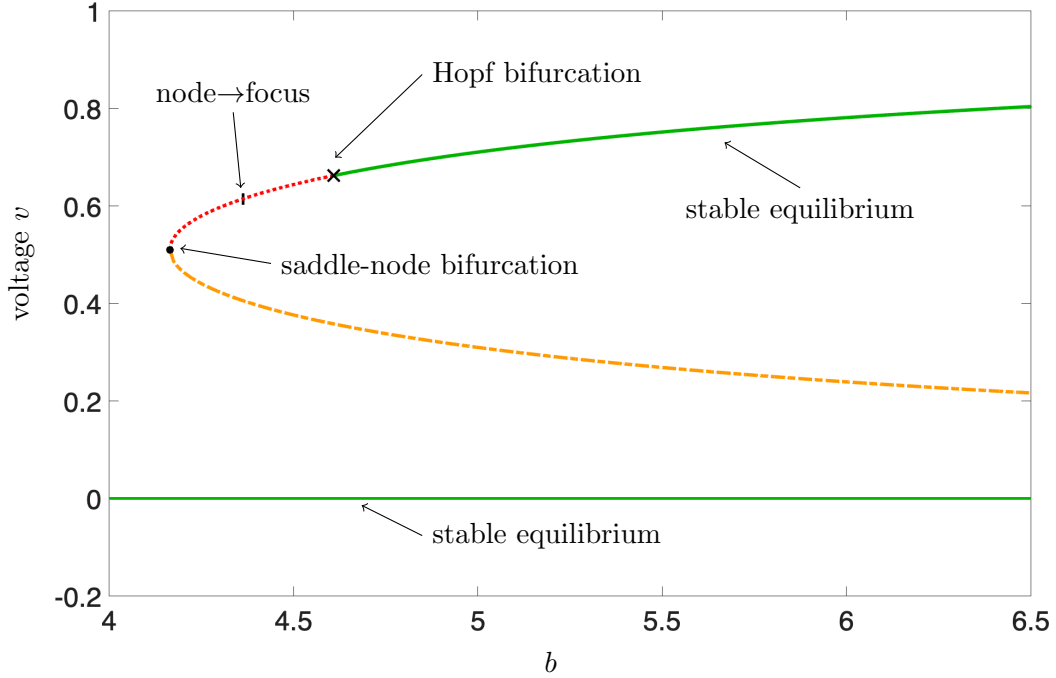

In \cref{fig:onecell_bif} we show the bifurcation diagram for the FHN model \cref{eqn:FHN} of one cell. 
All constants are those for a healthy cell, except $b$ is varied (it is the bifurcation parameter). 
The diagram shows that $(v,w) = (0,0)$ is an asymptotically stable equilibrium for all values of $b$. 
At $b \approx 4.1649313$, a saddle point (orange, dashed line) and an  unstable node (red, dotted line) are created via a saddle-node bifurcation. 
This curve is obtained by a direct algebraic solution.
At $b \approx 4.3629588$ the unstable node becomes an unstable focus. 
(A pair of positive real eigenvalues collide and become a complex pair of eigenvalues with positive real part.) 
The unstable focus is transformed into a stable focus by a Hopf bifurcation at $b \approx 4.6096588$. 

The creation of the stable equilibrium past $b \approx 4.61$ would theoretically mean that an isolated cell that is activated would get stuck in a high-voltage (depolarized) equilibrium. 
However, heart cells are never isolated; they are always electrically coupled to other heart cells. 
A one-cell bifurcation diagram, equivalent to \cref{fig:onecell_bif}, appeared in \cite{H-S-V-H:2021}. They noted that in heart tissue, charge leakage to surrounding cells generally prevents cells from becoming stuck in the high-voltage state.

\begin{figure}
    \centering
\begin{tikzpicture}
        \setlength{\figscale}{0.8\textwidth}
        \node[anchor=south west] at (0,0) {{\includegraphics[width=\figscale]{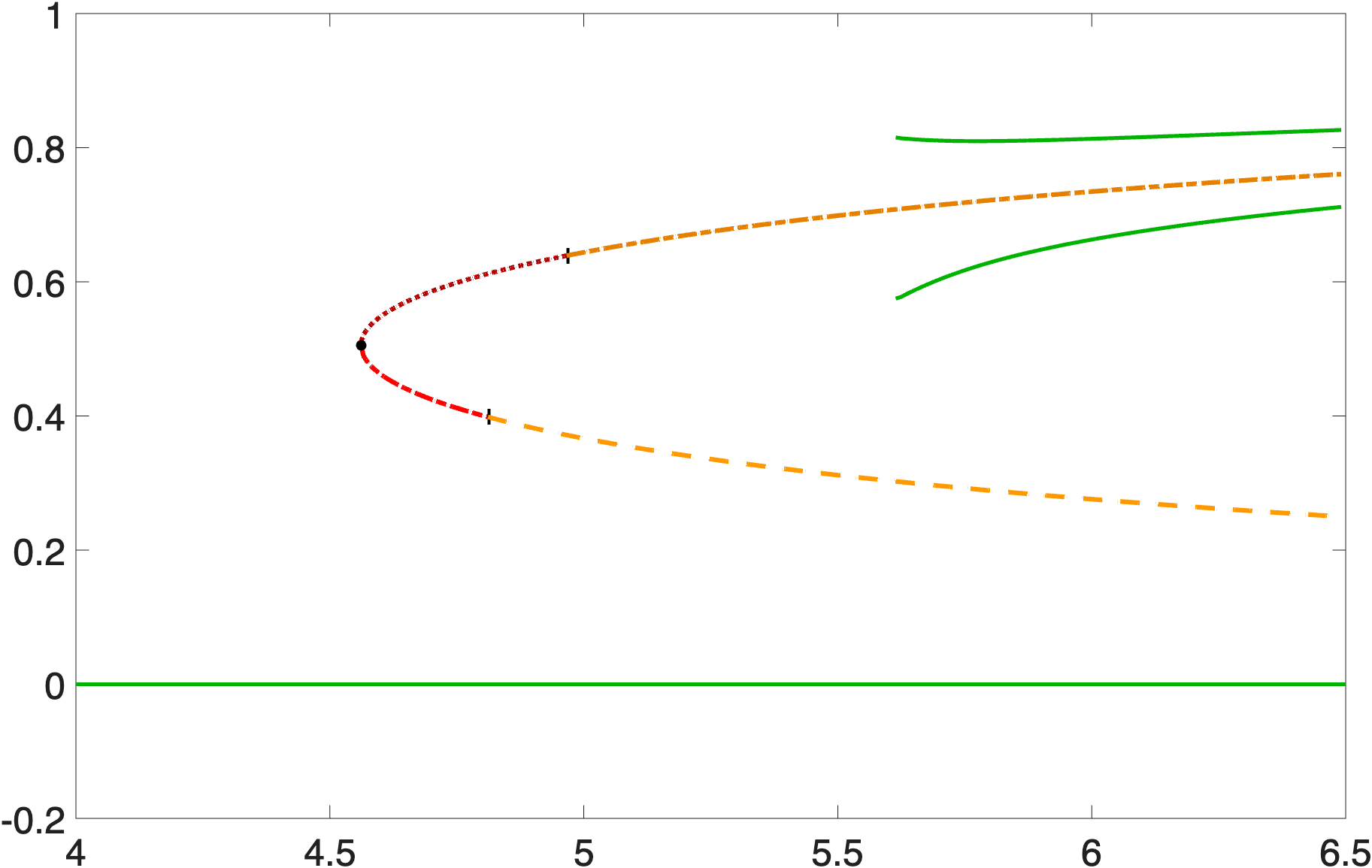}}};
        \node[rotate=90,anchor=south] at (0,5) {$v$ of slow-recovery cell};
        \node[anchor=north] at (6,0) {$b$ of slow-recovery cell};
        \draw[->] (5,1.4) node[anchor=west] {stable equilibrium} -- (4,1.8);
        \draw[->] (3.0,3.5) node[anchor=north] {saddle-node bifurcation} -- (3.5,5.1);
        \draw[->] (6.3,5.4) node[anchor=north] {Hopf bifurcations} -- (5.6,5.9);
        \draw[->] (6.0,4.8) -- (5.0,4.5);
        \draw[->] (11,5.6) node[anchor=north,align=center] {max and min of\\ stable periodic orbit} -- +(0,.5);
        \draw[<->] (11,6.33) -- +(0,0.74);
        \draw[->] (6.6,7.6) node[anchor=south] {reverse period-doubling bifurcation} -- +(2.0,-0.45);
        \draw[->] (6.6,7.6) -- +(2.0,-2.0);
\end{tikzpicture}
    \caption{
    Bifurcation diagram for a two cell system with one slow-recovery cell
    electrically coupled to one fast-recovery cell as in \cref{eqn:FHNcoupled}. 
    Unstable equilibria are generated by a saddle-node bifurcation
    at $b_1 \approx 4.56117$. The upper branch initially has 4 unstable eigenvalues and the lower branch initially has 3 unstable eigenvalues. 
    These equilibria undergo local bifurcations, but remain unstable.
    At $b_1 \approx 5.61397$ a reverse period doubling
    bifurcation creates a stable periodic orbit in which the slow-recovery cell never recovers, 
    but is trapped at an elevated voltage, which is triggering the fast-recovery cell to periodically fire.
    For $b_1$ below this bifurcation in the two-cell model
    the slow-recovery cell can induce multiple firings of the fast-recovery cell, but itself still recover. }
    \label{fig:twocell_bif}
\end{figure}
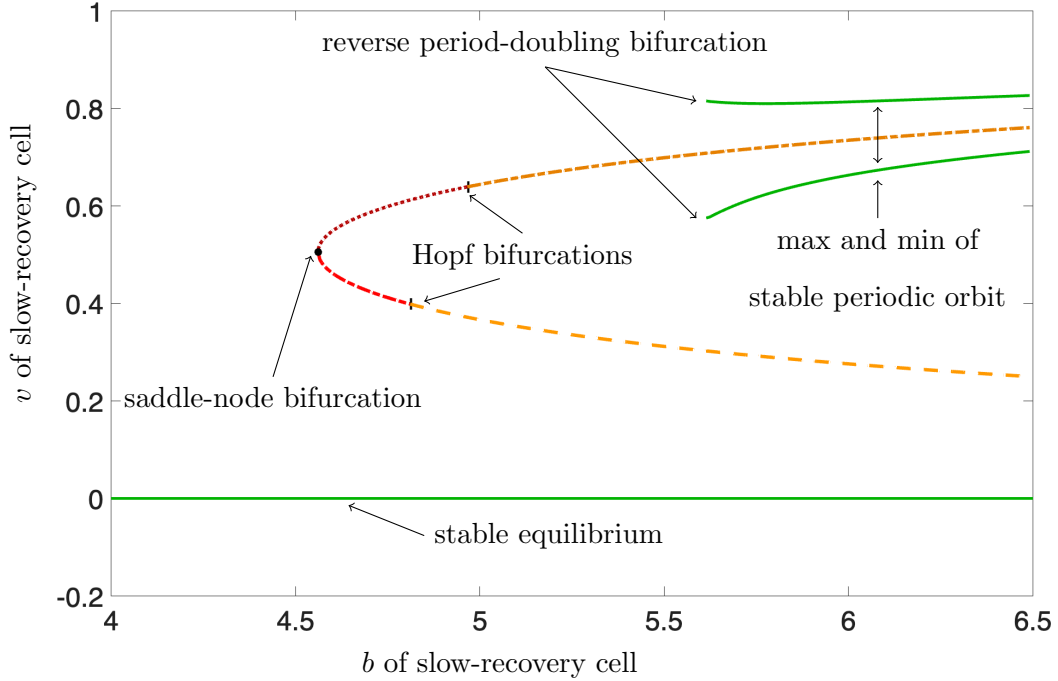

In order to understand the dynamics more fully, we next consider a two-cell model where cell 1 is a slow-recovery cell and cell 2 is a fast-recovery cell.  
The slow-recovery cell parameters were set to $\sigma_1 = 0.6$, $\epsilon_1 = .004$, $a_1 = 0.01$, and $c_1 = 1.0$; the fast-recovery cell parameters were set to $\sigma_2 = 1.6$, $\epsilon_2 = .008$, $a_2 = 0.02$, $b_2 = 2.9$, and $c_2 = 1.0$; and the conductance between them was set to $g = .028$.
In \cref{fig:twocell_bif}, we present a numerically generated bifurcation diagram for this model using the slow-recovery cell parameter $b_1$ as the bifurcation parameter and the voltage of the slow-recovery cell $v_1$ response variable. 
The only stable equilibrium is the zero solution. 
Unstable equilibria are generated by a saddle-node bifurcation at $b_1 \approx 4.56117$. 
The upper branch initially has 4 unstable eigenvalues, but two of them become stable via a Hopf bifurcation at $b_1 \approx 4.963$.
The lower branch initially has 3 unstable eigenvalues and two of them become stable at $b_1 \approx 4.811$ again via a Hopf bifurcation. 
In all cases, these equilibrium points remain unstable, so they would not play a vital role in the dynamics that are observed.

\begin{figure}
    \centering
    \includegraphics[width=0.7\linewidth]{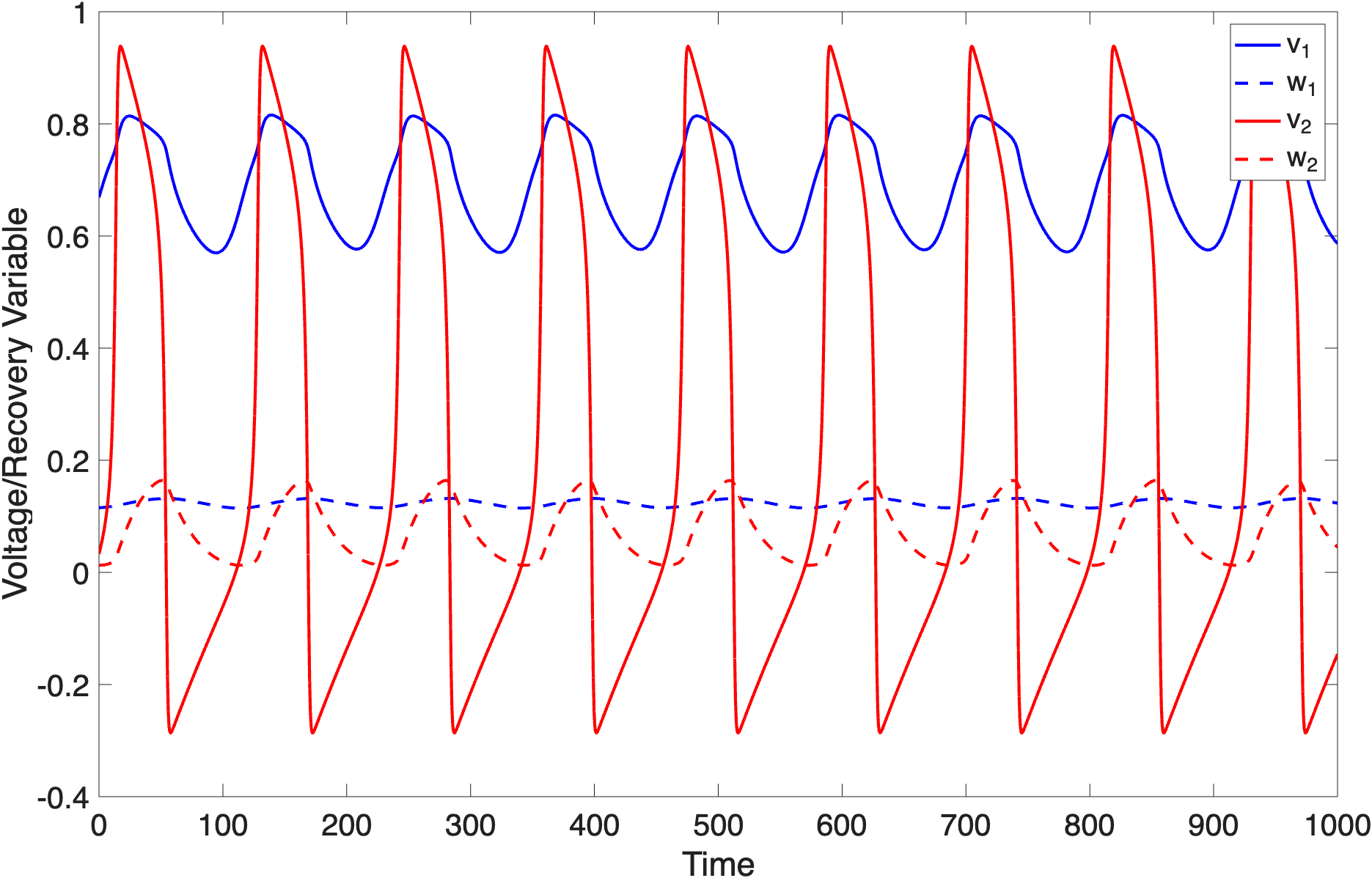}
    \caption{Time series for a stable periodic solution in the two-cell model 
    with $b_1=5.61397$, 
    which is just to the right of the reverse period doubling bifurcation in  \cref{fig:twocell_bif}. Here $v_1, w_1$ (blue) are the potential and recovery variable of the slow-recovery cell and $v_2, w_2$ (red) represent the fast-recovery cell.
    The slow-recovery cell never recovers from the  depolarized state, 
    while the fast-recovery cell repeatedly fires.}
    \label{fig:stable_periodic}
\end{figure}

At $b_1 \approx 5.61397069328$, a reverse period-doubling bifurcation creates a stable periodic orbit in which the slow-recovery cell is trapped in a periodic state at elevated voltage, which triggers the fast-recovery cell to periodically fire (see \cref{fig:stable_periodic}).
Such a scenario would kill the slow-recovery cell (due to \ce{Ca^{2+}} overload \cite{C-Z-K-H-M-Metal:2005}), so we consider it non-physiologic. 
In our large scale simulations, for $b \gtrsim 5.75$ we do not observe multiple firings because the slow-recovery cells, including the ``sensing cell", become stuck in the high-voltage state and never recover. 
This explains the dark regions on the far right of the panels in \cref{fig:OSC 7 FIB}. 

\begin{figure}
    \centering
    \includegraphics[width=0.7\linewidth]{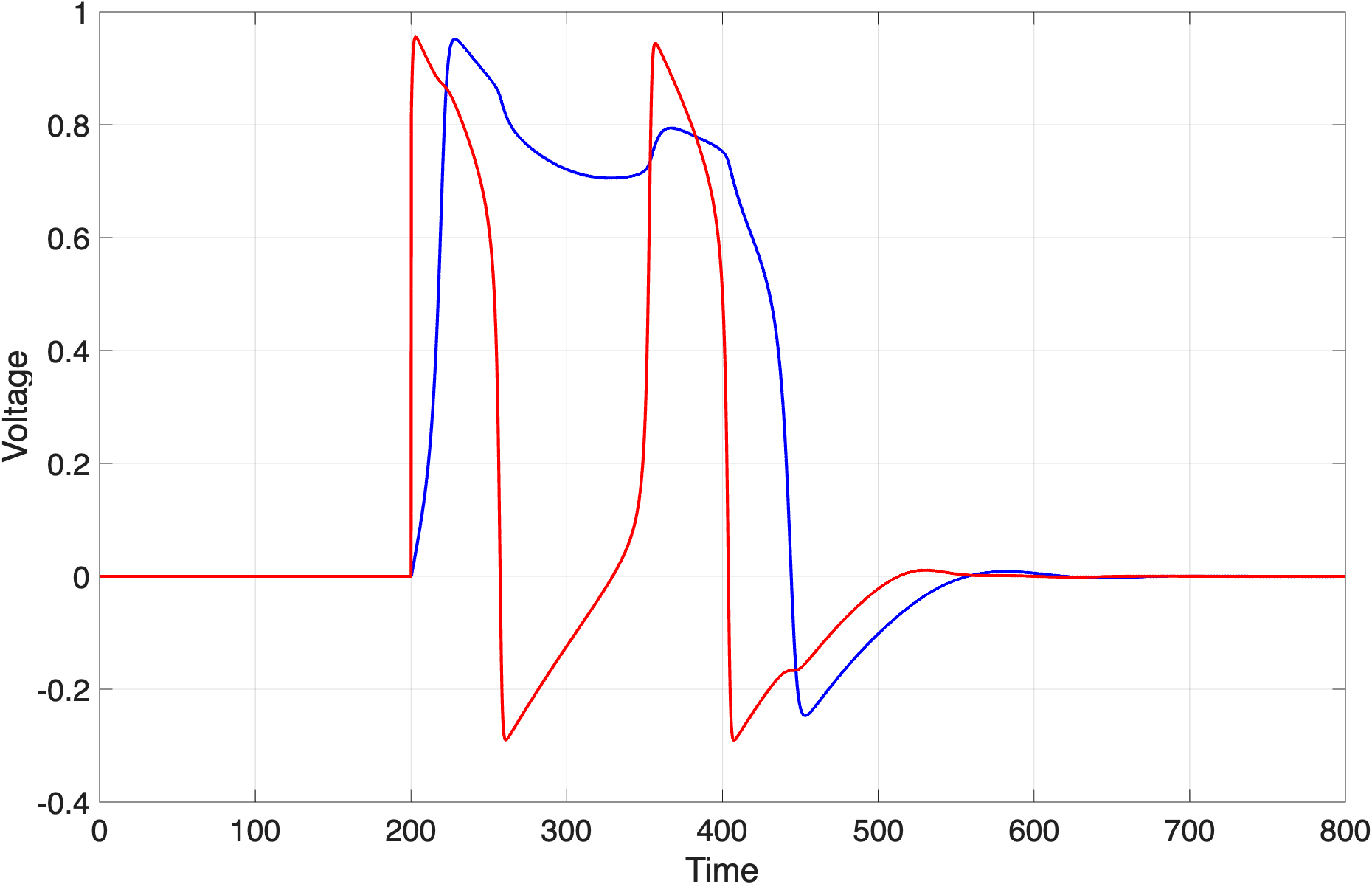}
    \begin{tabular}[b]{cc}
         \includegraphics[width=0.12\linewidth]{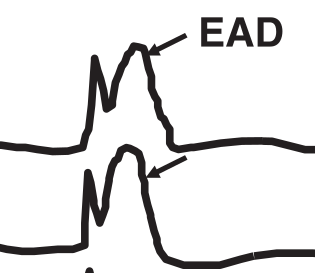} & \cite[Fig.~3]{O-H-K-L-L-W-C-K:2007}\\
         \includegraphics[width=0.12\linewidth]{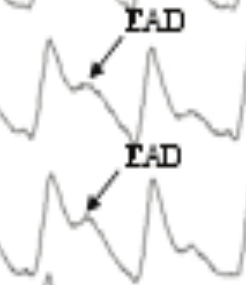} &\cite[Fig.~4]{P-C-F-B-W-K:2018}\\
         \includegraphics[width=0.12\linewidth]{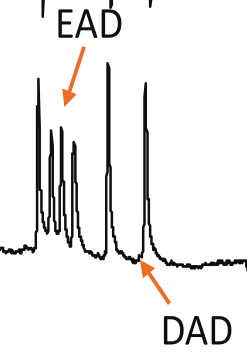} & \cite[Fig.~5]{P-C-F-B-W-K:2018}
    \end{tabular}
    \caption{{\bf Left:} Time series for a two-cell FHN model of a slow-recovery cell and a fast-recovery cell, with $b_1=5.5$, 
    which is to the left of the reverse period doubling bifurcation in  \cref{fig:twocell_bif}. Here $v_1$ (blue) is the voltage of the slow-recovery cell and $v_2$ (red) is the voltage of the fast-recovery cell. 
    A short pulse is applied to the fast-recovery cell, triggering an initial firing in both cells.
    The slow-recovery cell stays at high voltage long enough to trigger the fast-recovery cell to fire a second time, which causes a voltage increase in the slow-recovery cell during its plateau phase, before it eventually recovers.
    {\bf Right:} Snippets of early afterdepolarizations (EADs) from the literature.
    }
    \label{fig:twocell_pulse}
\end{figure}

For $b_1$ below $\approx 5.75$ in the two-cell model, the slow-recovery cell can induce an additional firing of the fast-recovery cell,  but itself still recover. We observe this in \cref{fig:twocell_pulse} where we have used the same two-cell model as in \cref{fig:stable_periodic}, but with $b_1 = 5.5$, before the onset of a stable periodic orbit.
This is the situation where we observe the generation of fibrillation-like behavior in our large scale simulations.

We note that the voltage of the slow-recovery cell pictured in \cref{fig:twocell_pulse} is reminiscent of some instances of early afterdepolarization (EAD), a clinical condition in which a heart cell has an additional voltage spike or spikes (depolarizations) while already at elevated voltage (during phase 2 or 3) \cite{O-H-K-L-L-W-C-K:2007,ANT-BUR:2011}.
In \cite{P-C-F-B-W-K:2018}, microelectrodes were used to detect EADs at the cellular level and these were observed to trigger atrial fibrilation in laboratory preparations of rat heart tissue.
Single-cell models with more variables (representing ion channels) can produce EAD's  by using differences in channel rates, effectively accentuating the slow and fast internal dynamics of the models \cite{T-S-Y-W-G-Q:2009,B-J-M-P-S:2023}. 
We conjecture that there are multiple mechanisms to produce EADs, and the coupled-cell dynamics we observe is one of them.

In summary, while the FHN model for a single cell does not adequately capture features of the dynamics in heart tissue, in a two-cell FHN model we can find clear geometric explanations of the local dynamics exhibited at an interface between fast- and slow-recovery heart cells.





\section{Robustness of the FIB Mechanism} 
\label{sec:variationgeometry}
\label{sec:global}

Our prototypical simulation in \cref{sec:protosimulation} used the geometry in \cref{fig:simbasicgrid} and specific values for the cell parameters $b$ and $\epsilon$ and the conductance coefficients.
In \cref{sec:local}, we considered the effect of changing $b$ and $\epsilon$.
In this section, we consider the effects of changing the conductance and of moderate changes to the geometry, with the goal of demonstrating that the FIB mechanism is robust.
In \cref{sec:exoticsimulations}, we will consider more exotic changes to the geometry.

Our prototypical simulation used conductance coefficients $g_{\mathrm{horz}} = 0.027$ and $g_{\mathrm{vert}} = 0.015$, giving an anisotropy ratio of 1.8. 
We found that when we change these coefficients but keep their ratio, sustained refirings can still be obtained by simply adjusting the $\epsilon$ and $b$ values.
For instance, if we double both coefficients, then voltage travels more quickly between cells, so we need to increase the $b$ for slow-recovery cells to ensure that they retain their voltage long enough to excite the fast-recovery cells.    
We also found that the qualitative behavior is preserved up to ratios of about 4, but at higher ratios the coherence of the wave fronts breaks down after the first few refirings.

\begin{figure}
    \centering
    \begin{tabular}{cc}
         Beat Count& End State  \\
         \includegraphics[width=0.485\textwidth]{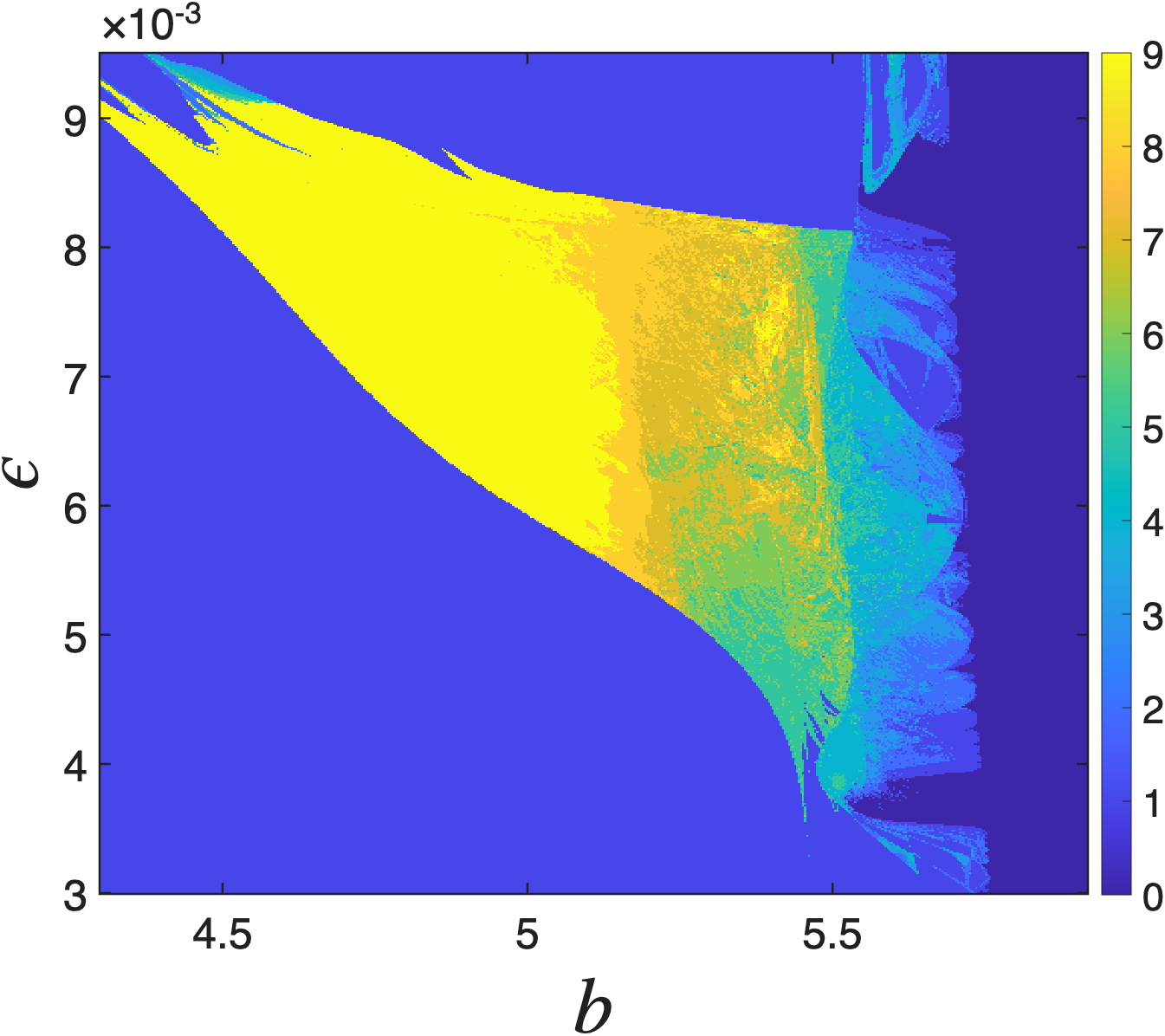}
         & \includegraphics[width=0.48\textwidth]{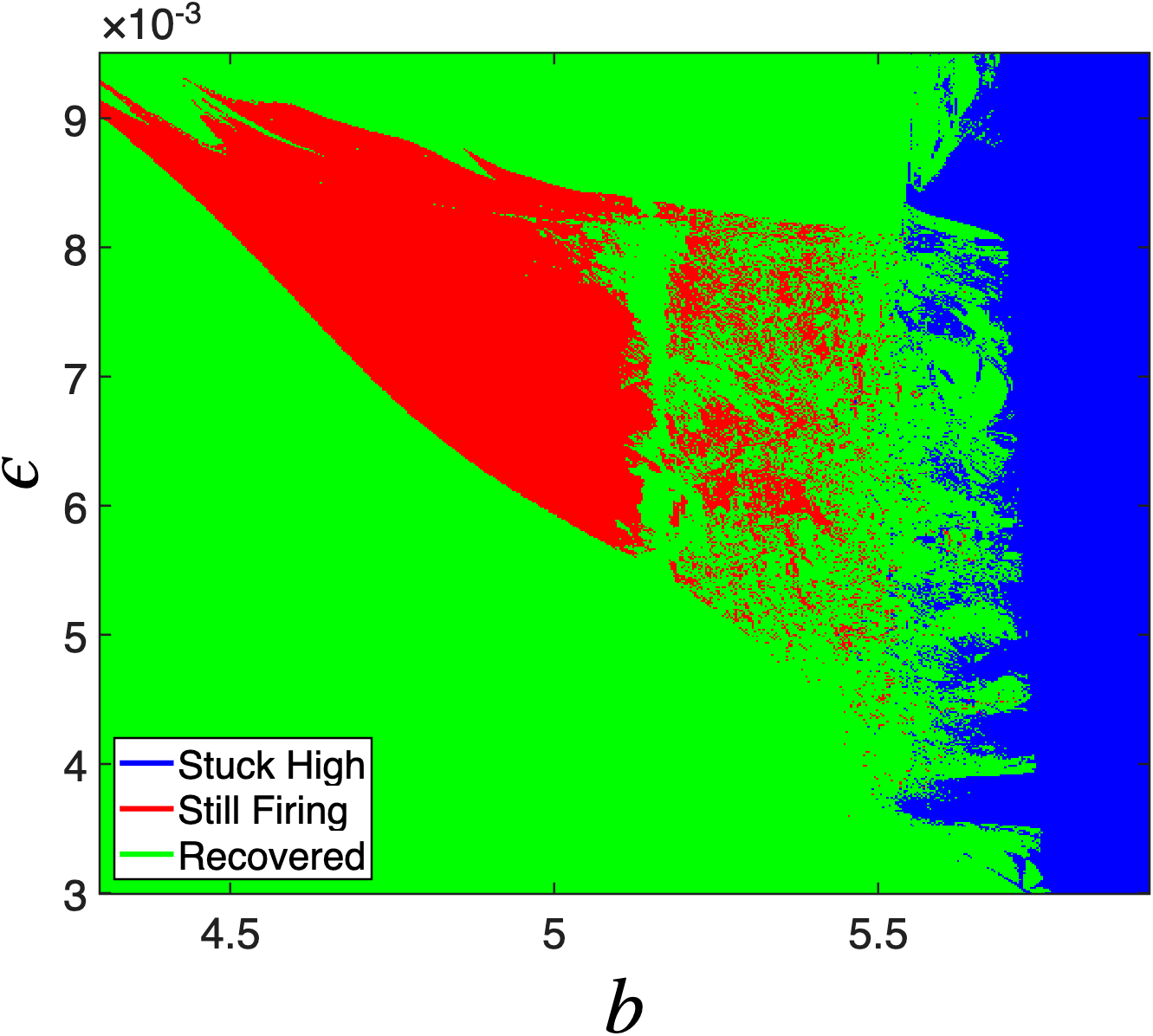}
    \end{tabular}
    \caption{Beat count (left) and end state (right) of the sensing cell as a function of $(b, \epsilon)$ in a simulation with a $4\times 27$ FIB.
    Comparing with \cref{fig:OSC 7 FIB}, which used a $4\times 7$ FIB, we see that the general shape is preserved, with a diagonal boundary for the self-sustaining behavior and a vertical boundary for the slow-recovery cells becoming stuck at elevated voltage. 
    However, the larger FIB produces the self-sustaining behavior more robustly.} 
    \label{fig:FIBlengths}
\end{figure}

Our prototypical simulation had a FIB (region of unhealthy cells) of dimension $4\times 7$. 
In \cref{fig:FIBlengths}, we show the number of beats and end state for a $4\times 27$ FIB, in the same format as \cref{fig:OSC 7 FIB}.
The qualitative behavior is the same, but \cref{fig:OSC 7 FIB} has holes and fractal-like features within the region where multiple firing occurs, whereas in \cref{fig:FIBlengths} the set of parameters with multiple firings is large and solid.
Thus, smaller FIBs might yield intermittent irregular dynamics, while longer FIBS yield persistent fibrillation-like behaviors.
Such behavior is consistent with worsening of atrial fibrillation in patients as more cells become unhealthy.
We observe two distinct mechanisms for the additional robustness of longer FIBs.
First, the additional slow-recovery cells help each other retain voltage by reducing the leakage of charge to non-slow-recovery neighbors.
Second, the travel time of the wave from the fast-recovery layer to the back of the slow-recovery layer increases, which reduces mistiming effects that spontaneously end the refiring.



Increasing the thickness of the slow-recovery cell layer also helps these cells retain voltage by reducing charge leakage to non-slow-recovery neighbors. For layers 4-5+ cells thick, this can lead to non-recovery of the cells in the center of the slow-recovery layer.

Increasing the thickness of the fast-recovery cell layer leads to more robust FIBs, by providing an additional mechanism of refiring.
Some of the fast-recovery cells at the top and bottom recover before the wave passes around to the back of the FIB, and so can refire without needing the slow-recovery layer.

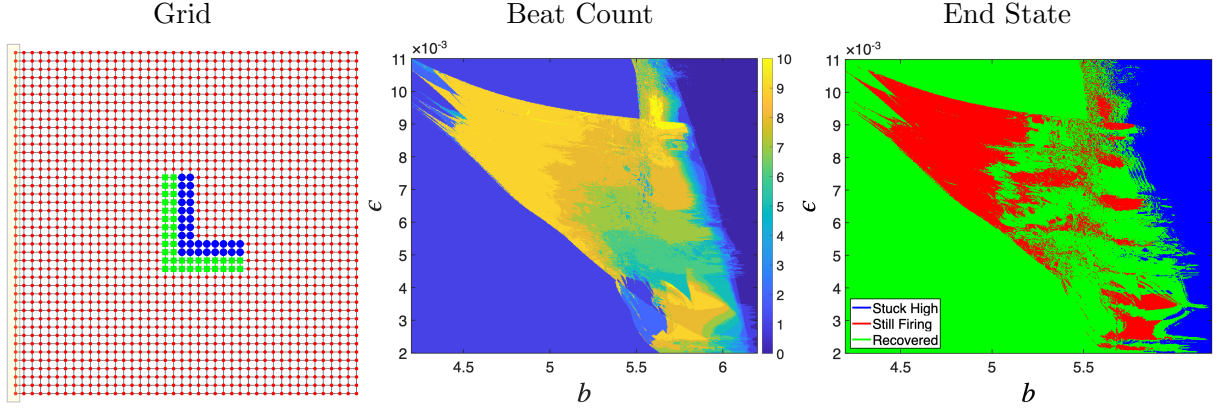
\begin{figure}
    \centering
    \begin{tabular}{ccc}
         Grid &Beat Count& End State  \\
        \begin{tikzpicture}[scale=0.11]
       \path[draw,help lines] (1,1) grid +(41,41); 
       {\foreach \x in {1,...,18,23,24,...,42} \foreach \y in {1,...,42} \draw[healthy,fill=red] (\x,\y) circle (0.15);} 
       \foreach \x in {19,20,21,22} \foreach \y in {1,...,17,25,26,...,42} \draw[healthy,fill=red] (\x,\y) circle (0.15);
       \foreach \x in {19,20} \foreach \y in {16,...,27} \draw[fast,fill=green] (\x-0.3,\y-0.3) rectangle +(0.6,0.6);
       \foreach \x in {21,...,28} \foreach \y in {16,...,17} \draw[fast,fill=green] (\x-0.3,\y-0.3) rectangle +(0.6,0.6);
       \foreach \x in {21,22} \foreach \y in {18,...,27} \draw[slow,fill=blue] (\x,\y) circle (0.4);
       \foreach \x in {23,...,28} \foreach \y in {18,...,19} \draw[slow,fill=blue] (\x,\y) circle (0.4);
       \draw[fill=yellow!50,nearly transparent] (0,0) rectangle +(1.5,43);
    \end{tikzpicture}
         &\includegraphics[width=0.34\textwidth]{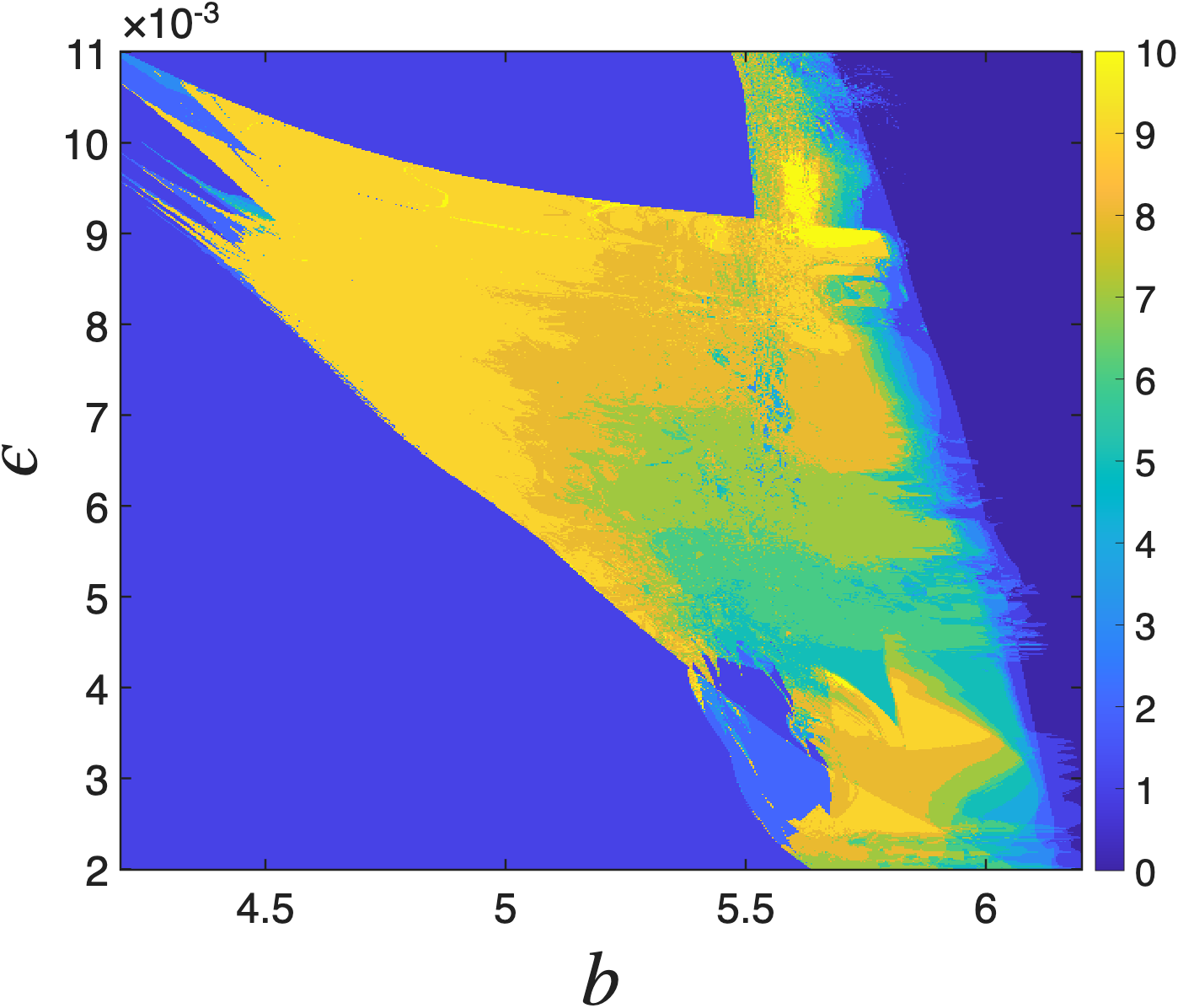}
         & \includegraphics[width=0.33\textwidth]{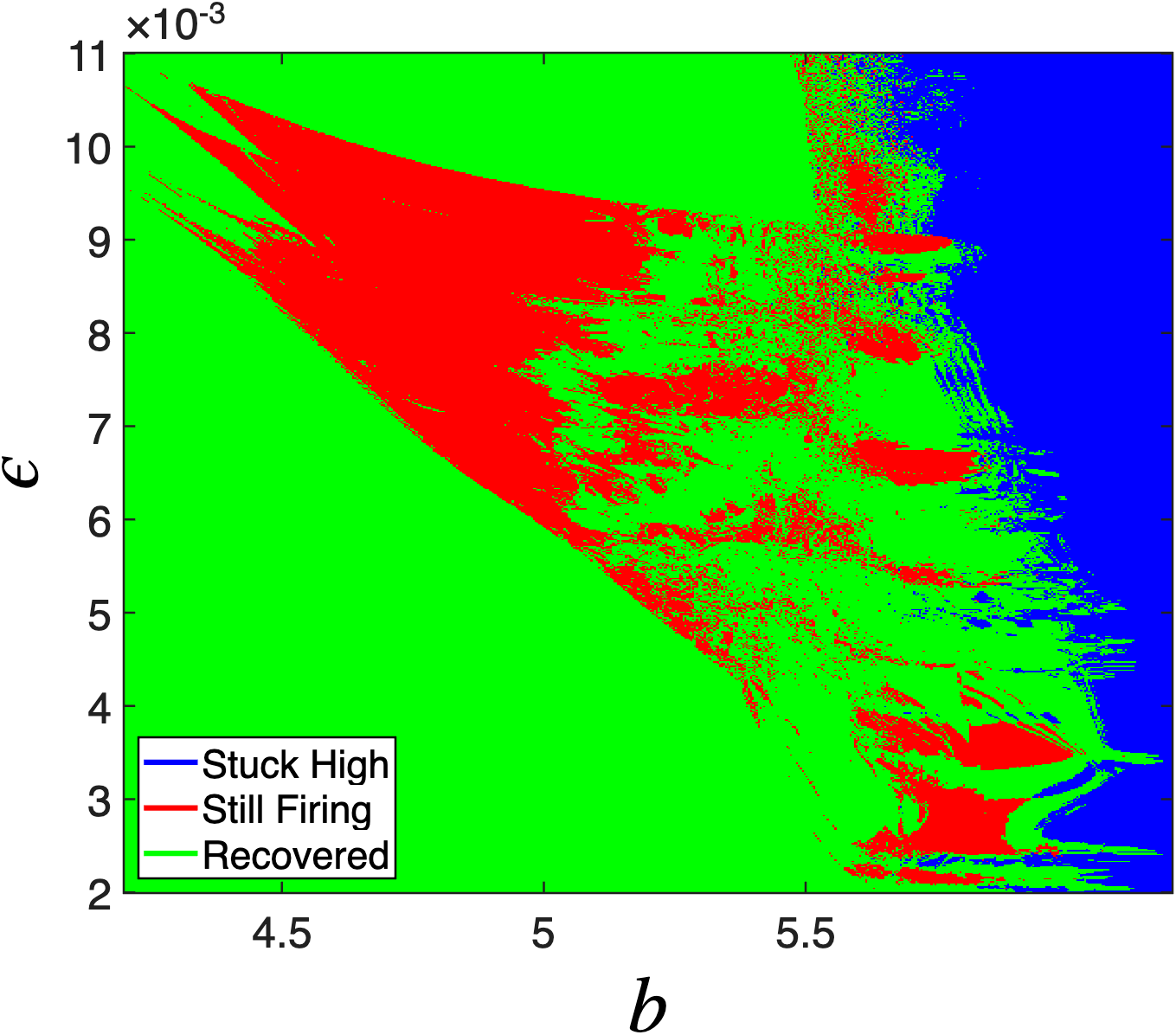}
    \end{tabular}
    \caption{Computational grid (left), beat count (center), and end state (right) as a function of $(b, \epsilon)$ in a simulation with an L-shaped FIB.  
    Comparing with \cref{fig:OSC 7 FIB,fig:FIBlengths}, the diagonal and vertical boundaries remain, but have shifted and the vertical boundary is slanted; note that the current figure includes larger $b$ and $\epsilon$ than the previous figures.
    }
    \label{fig:LFIB}
\end{figure}

    
We broke the vertical symmetry of the FIB by adding a horizontal region of unhealthy cells at the bottom to form an L-shaped FIB.
In \cref{fig:LFIB}, we show the computational grid and the number of beats and end state for a L-shaped FIB.
The L-shaped FIB is also more robust than the prototypical model.
It benefits from being effectively taller, which helps with voltage retention and increased travel time, and from effectively having a wider fast-recovery layer, which gives an additional refiring mechanism.
The wave shapes are more erratic, to the point that snapshots are not helpful.
A video is included in the Supplemental Material \cite{AF-supplement}.


Thus we have determined that the exact geometric arrangement of the FIB is not a limiting factor for sustained firings. So long as we have a continuous boundary between cell types, we can locate parameters in the model that result in sustained firings.
Overall, larger FIB regions lead to more fibrillatory-like behavior once the threshold of
local refiring is crossed.


\section{Spirals, Phantom Fibrillatory Centers, and Multiple Wave(let)s}
\label{sec:exoticsimulations}

Up to this point, we have taken a conservative approach to modeling, in order to demonstrate a simple mechanism that produces the characteristic features of atrial fibrillation.
We only considered very simple geometries (rectangular or L-shaped) and kept cell parameters and conductance coefficients physiologically plausible.
These conservative models were sufficient to produce disorganized waves of firing that persist is a self-sustaining manner, that are irregular, and that include double swirling-back curves.

In this section we loosen up a bit and include more exotic geometries, artificial initial conditions, and an additional cell type (scars).
With these, we can produce spirals, phantom fibrillatory centers, and multiple wave(let)s interacting.
Although in several cases the simulation results look like results from other models or observations from real atria, we cannot claim these looser models match reality without validation in real atria.
{\em If} the mechanisms we observe are real, then there are unfortunate implications for targeted ablation as a treatment for atrial fibrillation:
\begin{itemize}
    \item The occurrence of a spiral pattern does not imply any sort of structure at the apparent center of the spiral. 
    This fact was already known mathematically and explains the ineffectiveness of ablating spiral wave centers \cite{B-V-H-S-R-S-S-S:2025}. 
    \item The occurrence of a fibrillatory center, which repeatedly generate new waves, also does not imply any sort of structure at the apparent center. 
    Such phantom centers can be generated by a nearby FIB or by repeated external stimulation (pacing).
   
    \item When multiple FIBs are present near each other, they can turn each other on and off.
    Thus, some FIBS can be hidden.    
\end{itemize}
The common theme is that ablation targeting the interesting place might only kill healthy cells or might hit only one of many actual FIBs.



\subsection{Spirals}
\label{sec:spirals}

If we consider the limit of a two-dimensional grid of cells following the FitzHugh-Nagumo equations as the grid spacing goes to zero, we obtain a partial differential equation of the reaction-diffusion type. 
The existence of spiral solutions to reaction-diffusion partial differential equations in two dimensions has been known for some time, with studies on the stability of such spirals \cite{SAN-SCH:2000,SAN-SCH:2006,SAN-SCH:2007}, and the drifting and meandering motion of their centers \cite{SAN-SCH:2001,SAN-SCH:2007}.
Thus, it is expected that spiral solutions can arise here, even in a homogeneous grid of healthy cells.

\begin{figure}
    \centering
    \begin{tikzpicture} 
    \setlength{\figscale}{0.24\textwidth}
    \matrix[row sep=1, column sep=0, cells={scale=1}] at (0,0)
        {
        \node[anchor=south] at (0,0) {\includegraphics[width=\figscale,height=\figscale]{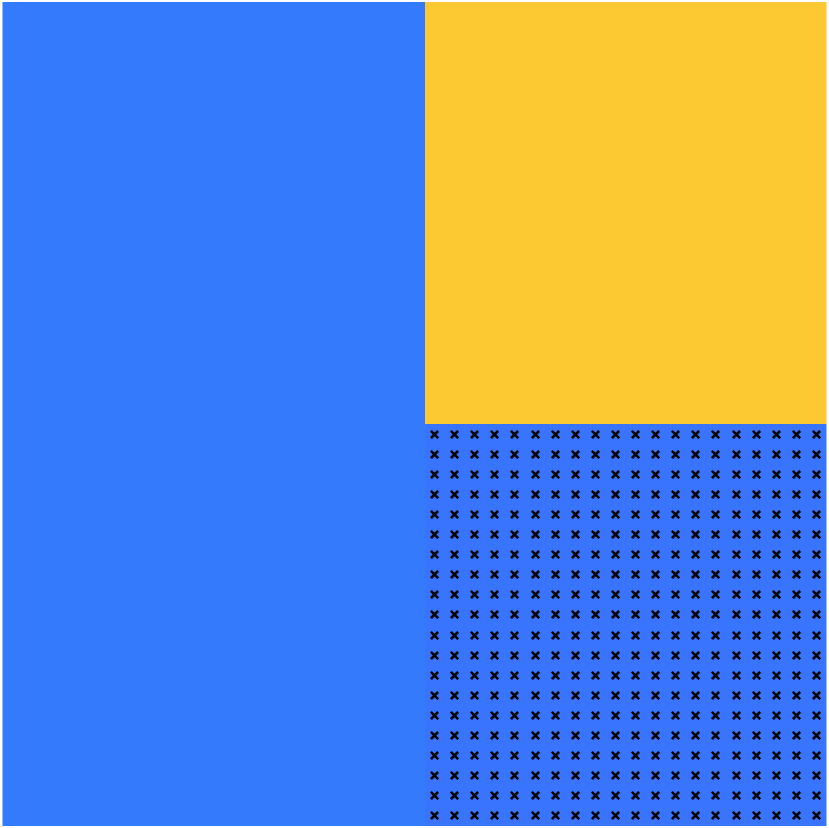}};
        \node[anchor=north] at (0,0) {$t=0$};
        &
        \node[anchor=south] at (0,0){\includegraphics[width=\figscale,height=\figscale]{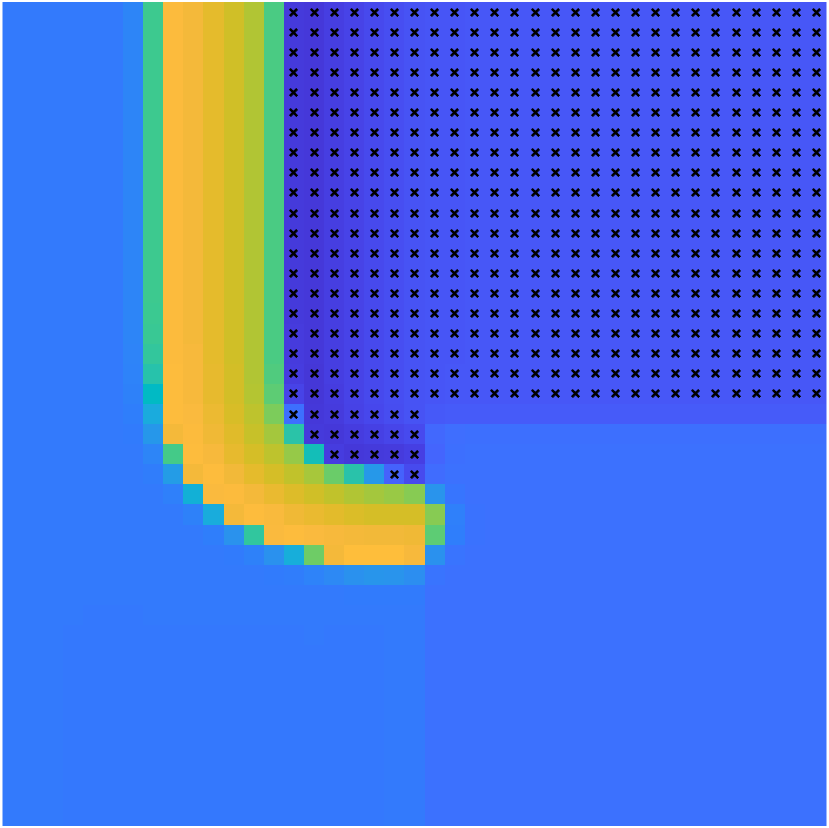}};
        \node[anchor=north] at (0,0) {$t=150$};
        &
        \node[anchor=south] at (0,0) {\includegraphics[width=\figscale,height=\figscale]{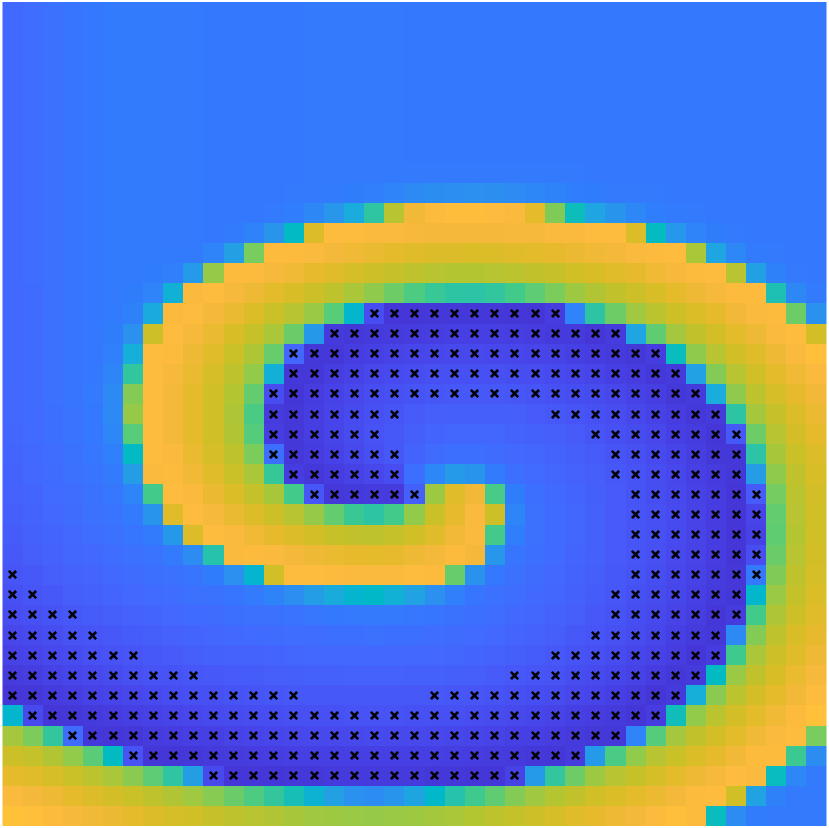}};
        \node[anchor=north] at (0,0) {$t=400$};
        &
        \node[anchor=south] at (0,0) {\includegraphics[width=\figscale,height=\figscale]{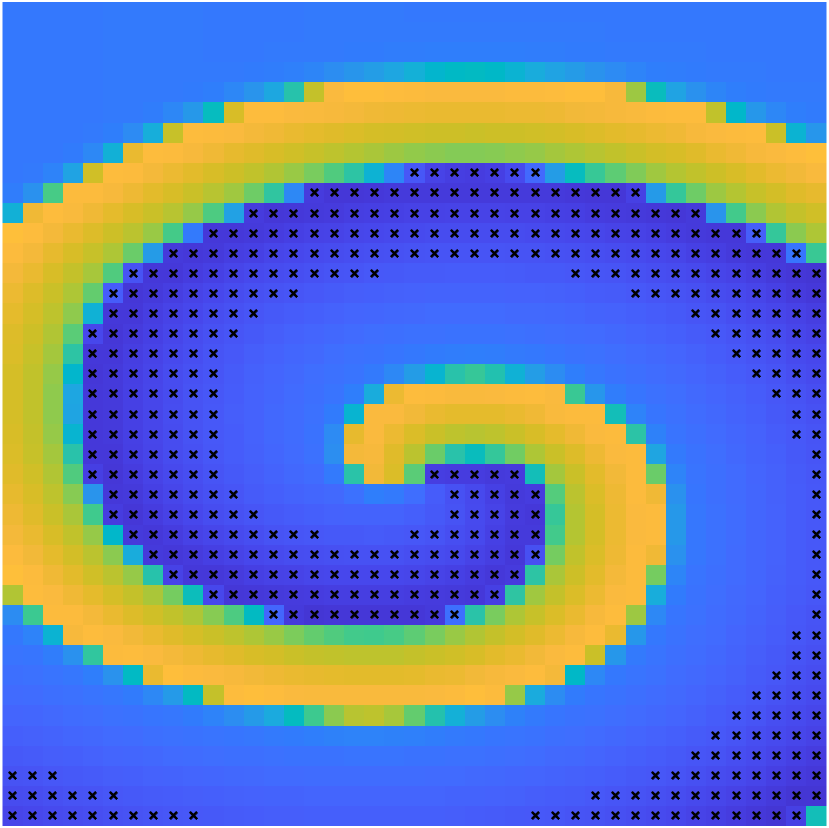}};
        \node[anchor=north] at (0,0) {$t=500$};
        \\
        };
    \end{tikzpicture}
    \caption{Snapshots of a simulation where setting the initial condition produces a spiral wave in healthy cells. 
    At $t=0$, cells in the upper-right quadrant have $v=1$ and $w=0$, cells in the lower-right quadrant have $v=0$ and $w=1/4$, and the remaining cells have $v=w=0$. We can observe that the wave initiated in the upper-right quadrant travels around the recovering tissue, forming a single broken end. Once the lower-right quadrant  recovers, the wave front travels into it and begins a spiral-like motion with an oscillating center.}
    \label{fig:healthyspiral}
\end{figure}
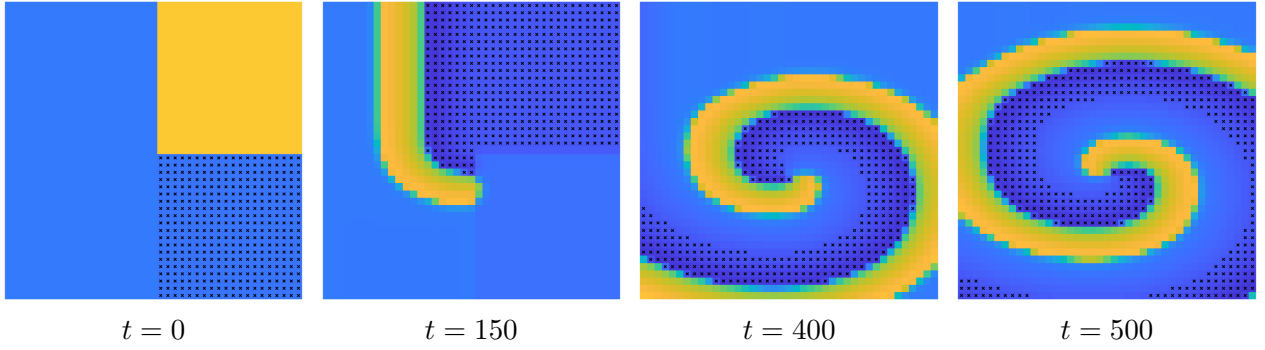

The direct cause of a spiral is known to be a {\em phase change point} \cite{ZYK-WIN:1987}, which appears as a broken end to the wavefront. 
As the main wavefront moves forward, the broken end moves sideway and curls back, forming the spiral shape.
In \cref{fig:healthyspiral}, we illustrate how setting initial conditions can lead to a persistent spiral.
In this simulation we see the creation of a single broken end which proceeds to curve back through the healthy tissue, forming a stable spiral.
A video is included in the Supplemental Material \cite{AF-supplement}.
Known mechanisms to create broken ends include mechanical perturbations (in experiments), a second stimulation during a narrow vulnerable time window, discontinuities in the wave travel speed with a second stimulation, obstacles with sharp corners, and wave breakage due to blocks caused by jumps in diffusion coefficients \cite{ZYKOV:2018}.


In our prototype simulation in \cref{fig:basicsnaps}, the FIB created two broken ends, which create double swirling-back curves.
By blocking one of these broken ends, the FIB can generate a spiral.
When functional heart cells (myocytes) die, they are replaced by scar cells (myofibroblasts), which are conductive but not excitable.
In isolation, a scar cell model is simply $dv/dt=0$, so when coupled to other cells it is modeled as $dv_i/dt=\sigma_i\sum_jg_{ij}(v_j-v_i)$; we use $\sigma=0.6$ to match our slow-recovery cells.
We adjust our prototypical model by adding a line of scar cells extending from the top of the FIB to the left boundary, and show snapshots in \cref{fig:scar-line}.
A video is included in the Supplemental Material \cite{AF-supplement}.
\begin{figure}
    \centering
    \begin{tikzpicture}
    \setlength{\figscale}{0.3\textwidth}
    \matrix[row sep=1, column sep=0, cells={scale=1}] at (0,0)
        {
        \node[anchor=south] at (0,0) {\includegraphics[width=\figscale]{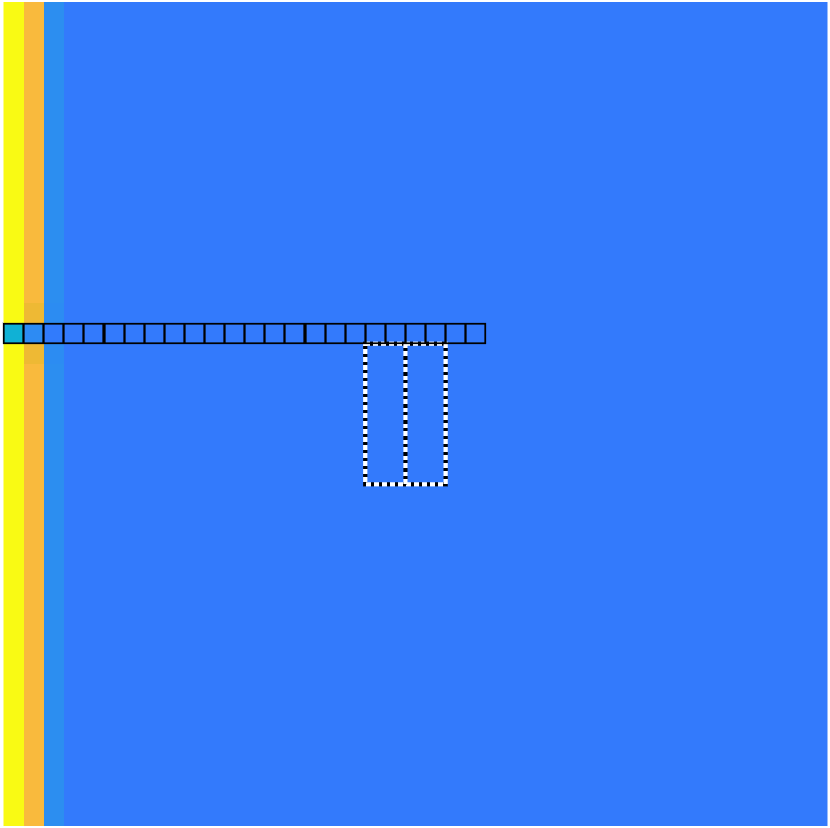}};
        \node[anchor=north] at (0,0) {$t=10$};
        &
        \node[anchor=south] at (0,0) {\includegraphics[width=\figscale]{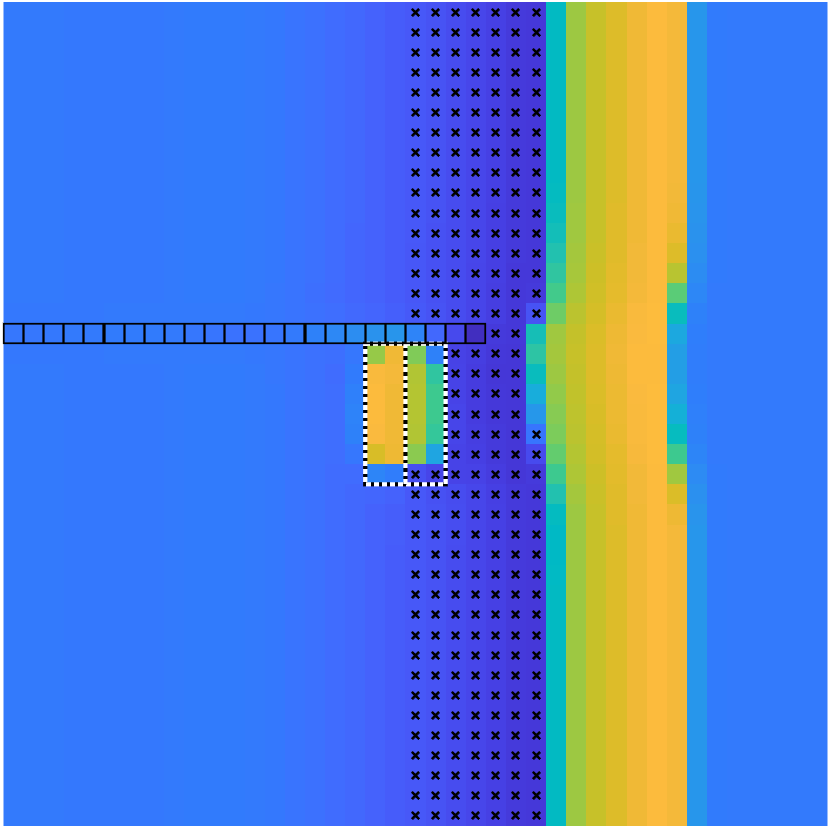}};
        \node[anchor=north] at (0,0) {$t=350$};
        &
        \node[anchor=south] at (0,0) {\includegraphics[width=\figscale]{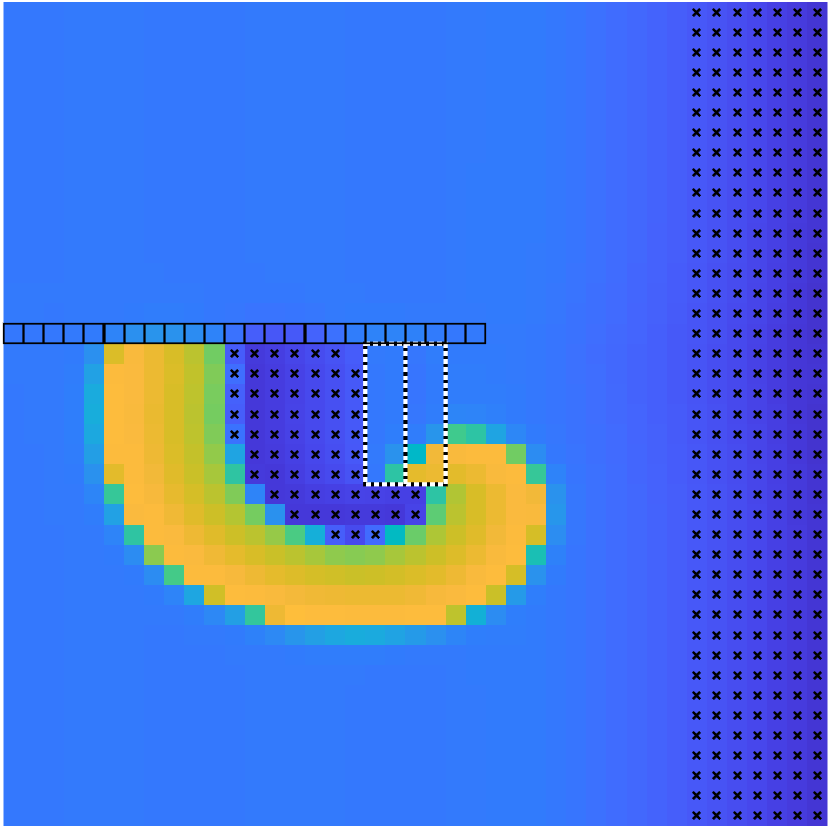}};
        \node[anchor=north] at (0,0) {$t=500$};
        \\
        \node[anchor=south] at (0,0) {\includegraphics[width=\figscale]{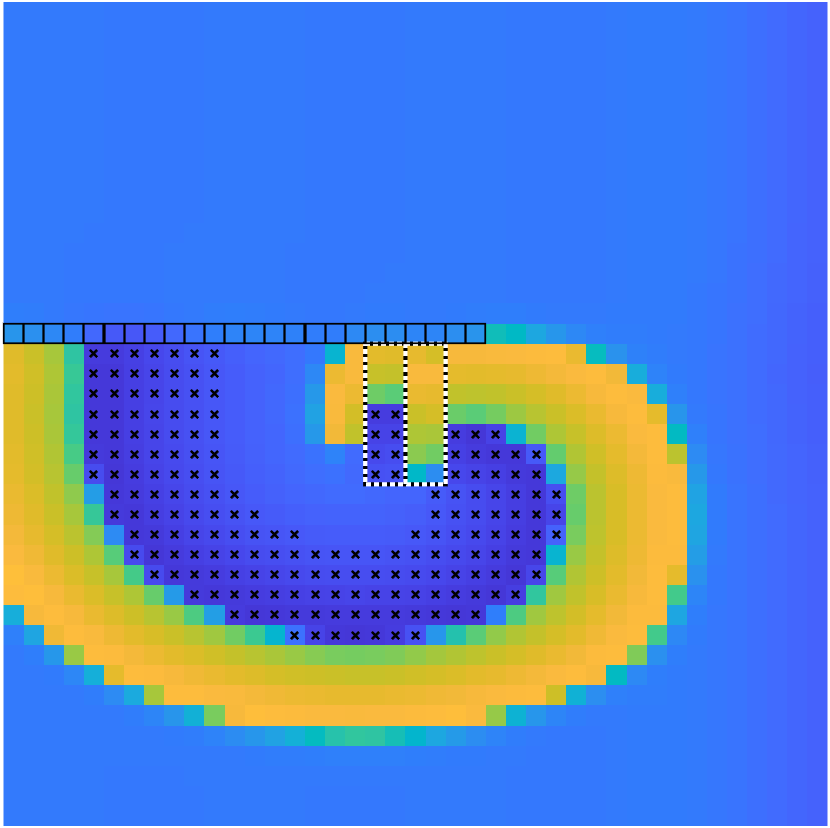}};
        \node[anchor=north] at (0,0) {$t=580$};
        &
        \node[anchor=south] at (0,0) {\includegraphics[width=\figscale]{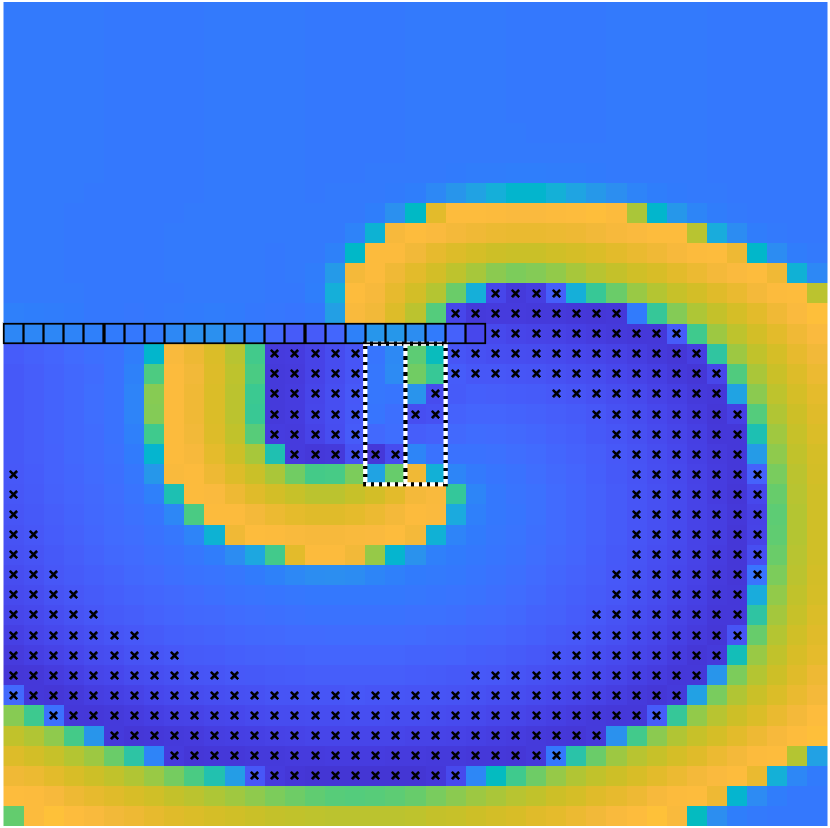}};
        \node[anchor=north] at (0,0) {$t=690$};
        &
        \node[anchor=south] at (0,0) {\includegraphics[width=\figscale]{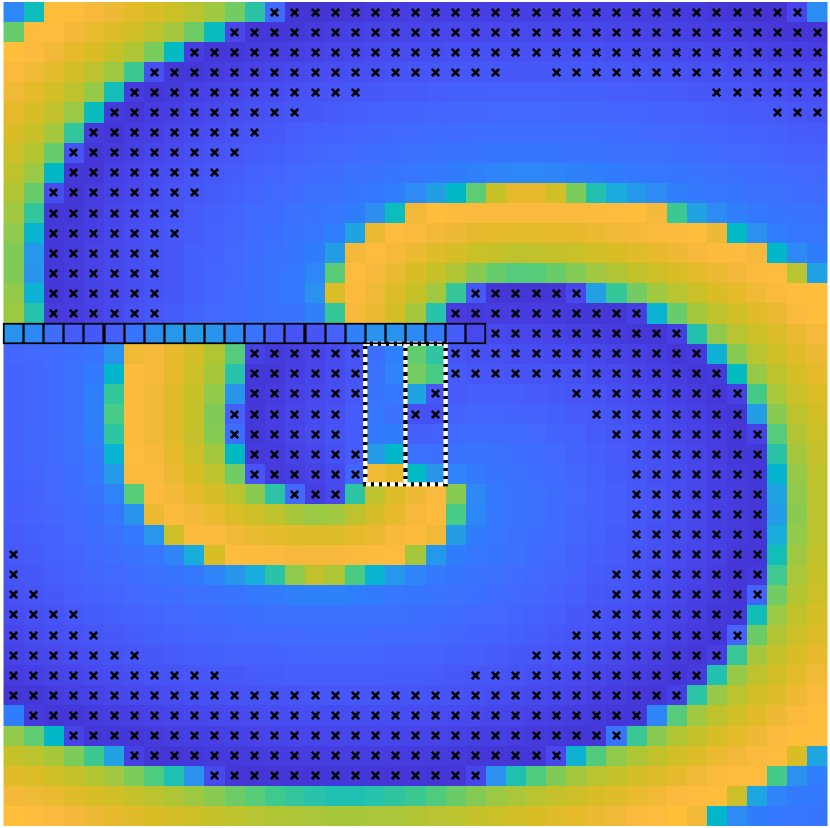}};
        \node[anchor=north] at (0,0) {$t=930$};
        \\
        };
    \end{tikzpicture}
    \caption{Snapshots of a simulation in which a line of scar cells (indicated by black borders) is included  directly above the prototypical FIB.
    The FIB fires in its usual self-sustained manner, and the scar shapes the emitted waves into multiple locally spiral-like waves.}
    \label{fig:scar-line}
\end{figure}
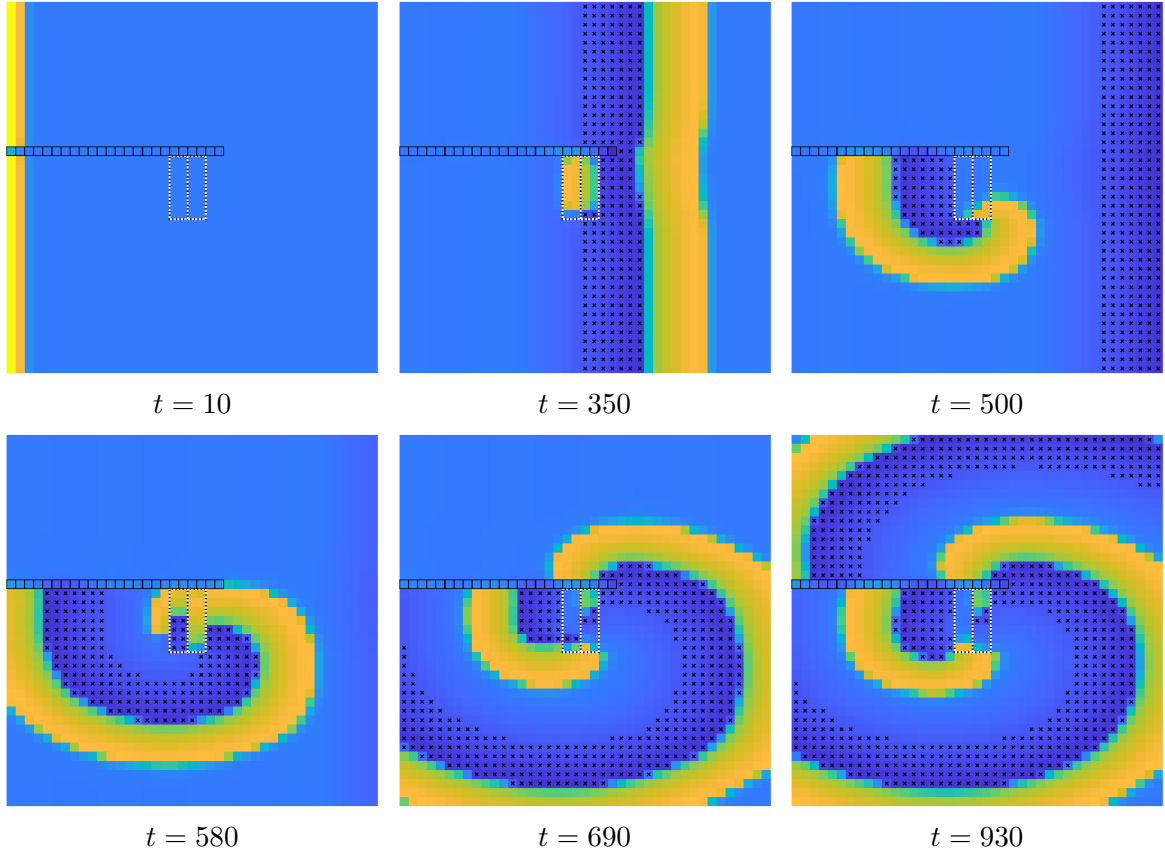
The single broken end creates a spiral, which rotates around and within the FIB and becomes self-sustaining.
When the wave emitted by this spiral center intercepts the scar tissue, the initial wavefront separates and creates a new broken end which travels along the line of scar tissue.




In \cite{A-A-K-Z-M-Betal:2019}, the authors produced spiral-like waves in an automaton model by sending a double-impulse diagonally toward a boundary between cells of different recovery times.
To produce a similar phenomenon, we embed a $6\times 81$ FIB in a  $81 \times 81$ computational grid and apply an impulse from the upper-left corner of the grid. 
In \cref{fig:vertical_strip_fib}, we show snapshots of the simulation.
A video is included in the Supplemental Material \cite{AF-supplement}.
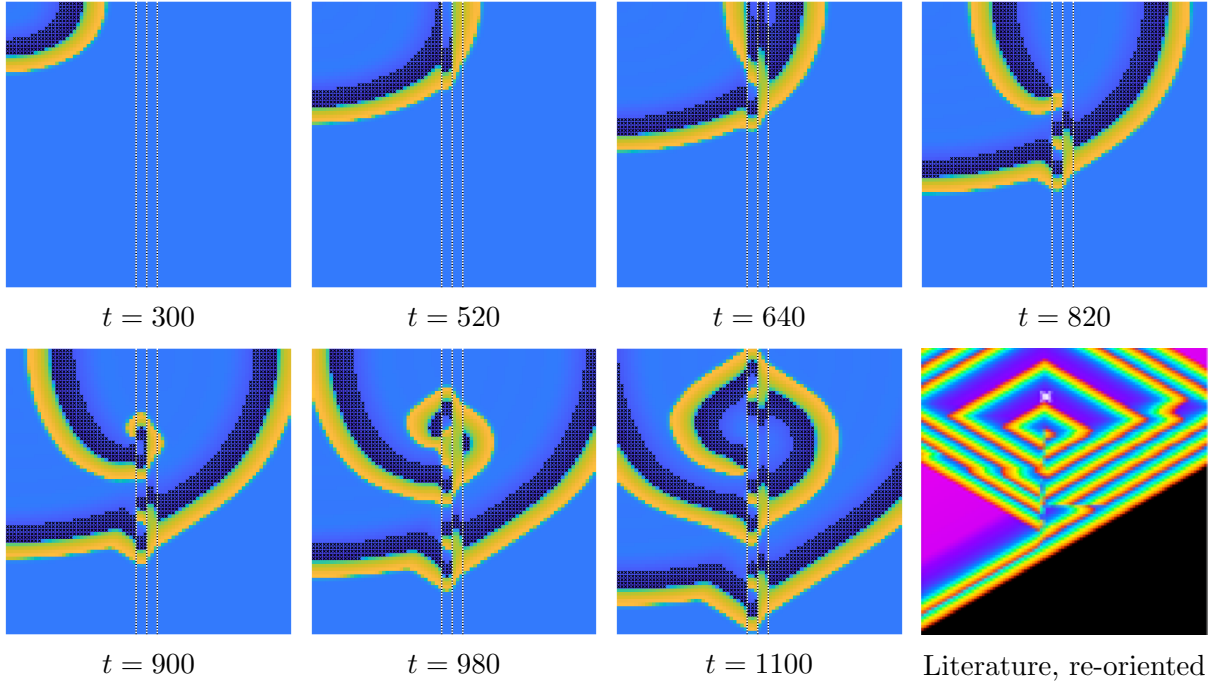
\begin{figure}
    \centering
    \begin{tikzpicture}
    \setlength{\figscale}{0.23\textwidth}
     \matrix[row sep=1, column sep=0, cells={scale=1}] at (0,0)
        {
        \node[anchor=south] at (0,0) {\includegraphics[width=\figscale]{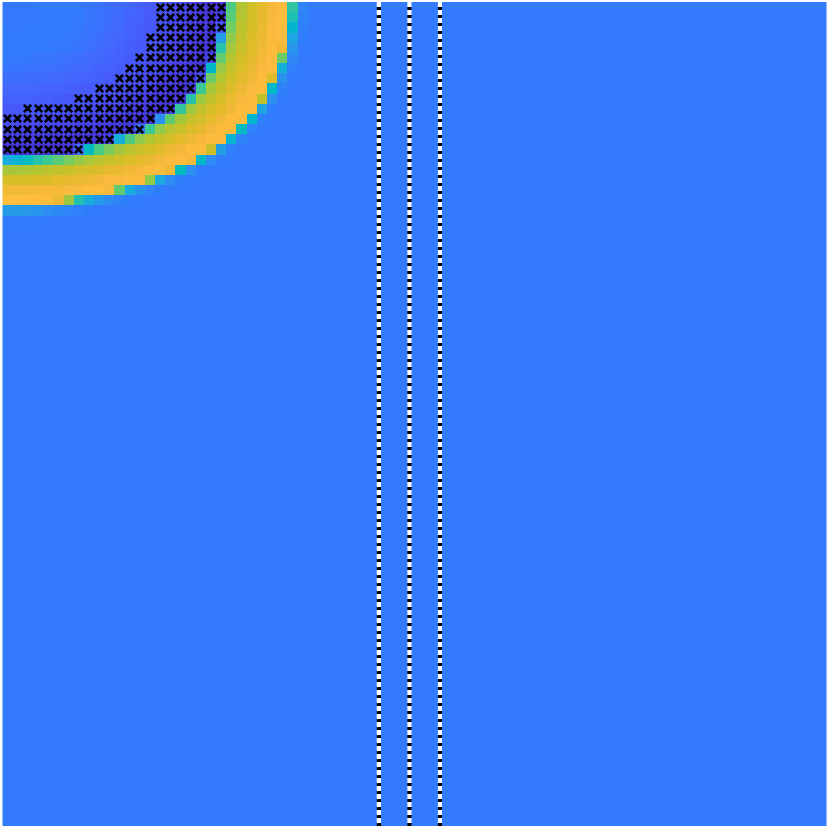}};
        \node[anchor=north] at (0,0) {$t=300$};
        &
        \node[anchor=south] at (0,0) {\includegraphics[width=\figscale]{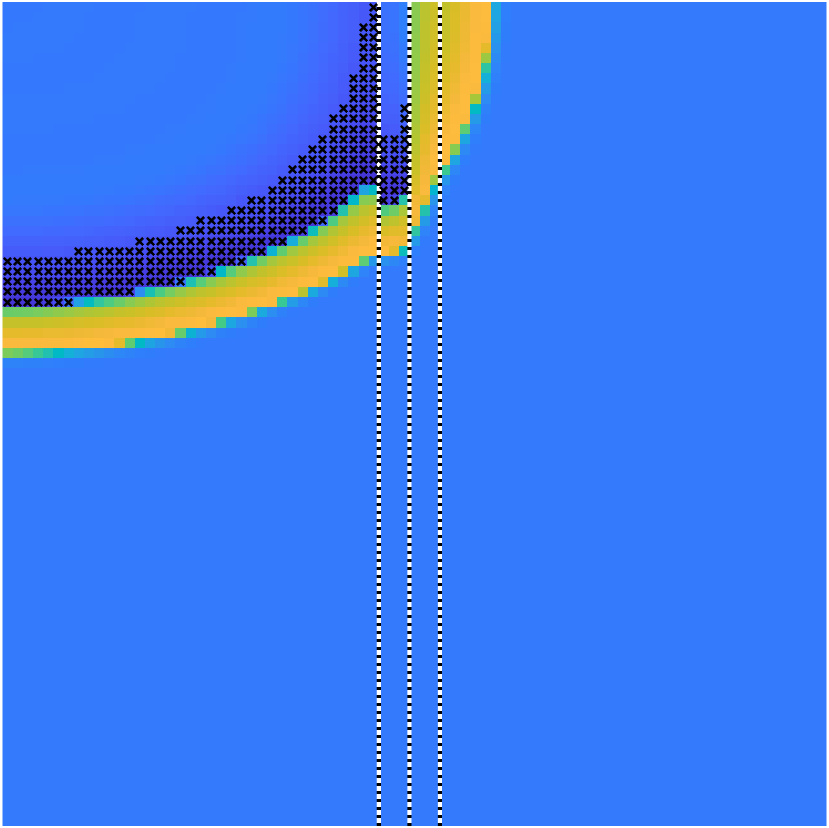}};
        \node[anchor=north] at (0,0) {$t=520$};
        &
        \node[anchor=south] at (0,0) {\includegraphics[width=\figscale]{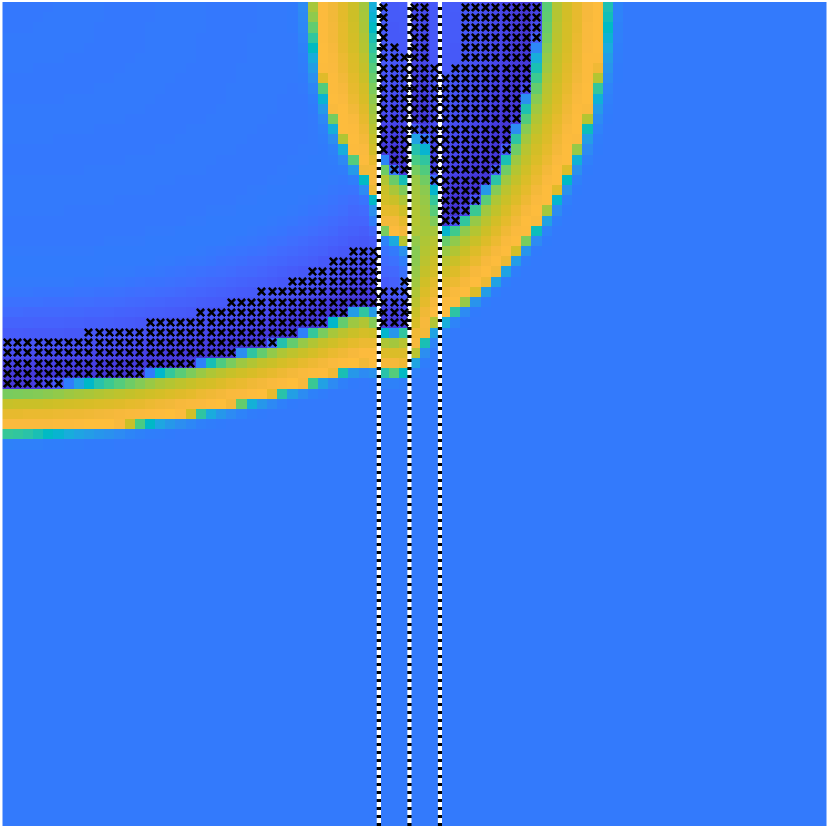}};
        \node[anchor=north] at (0,0) {$t=640$};
        &
        \node[anchor=south] at (0,0) {\includegraphics[width=\figscale]{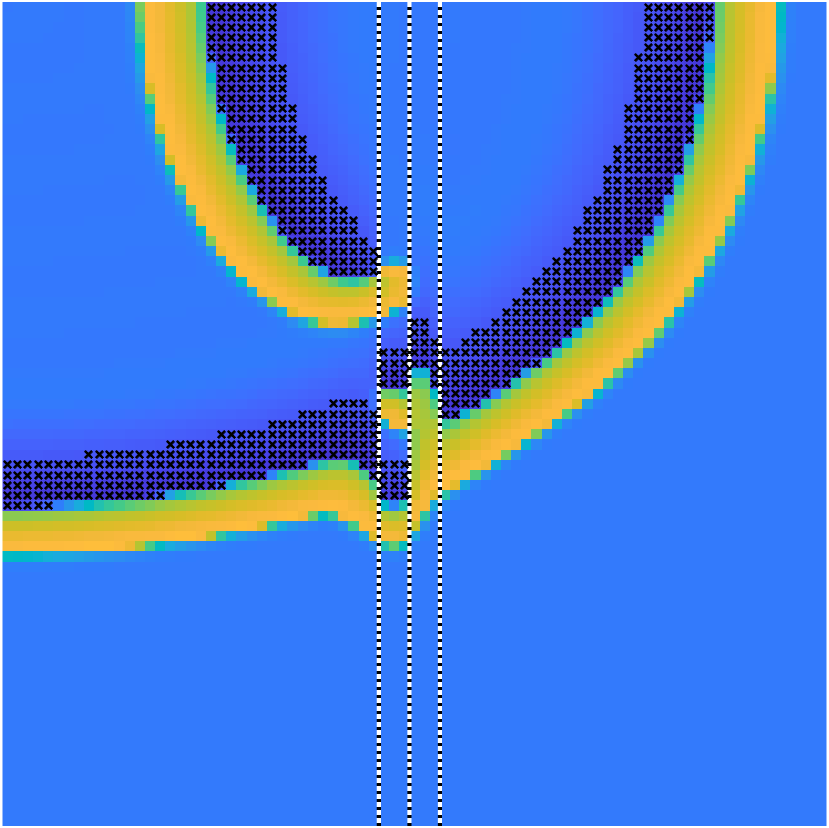}};
        \node[anchor=north] at (0,0) {$t=820$};
        \\ 
        \node[anchor=south] at (0,0) {\includegraphics[width=\figscale]{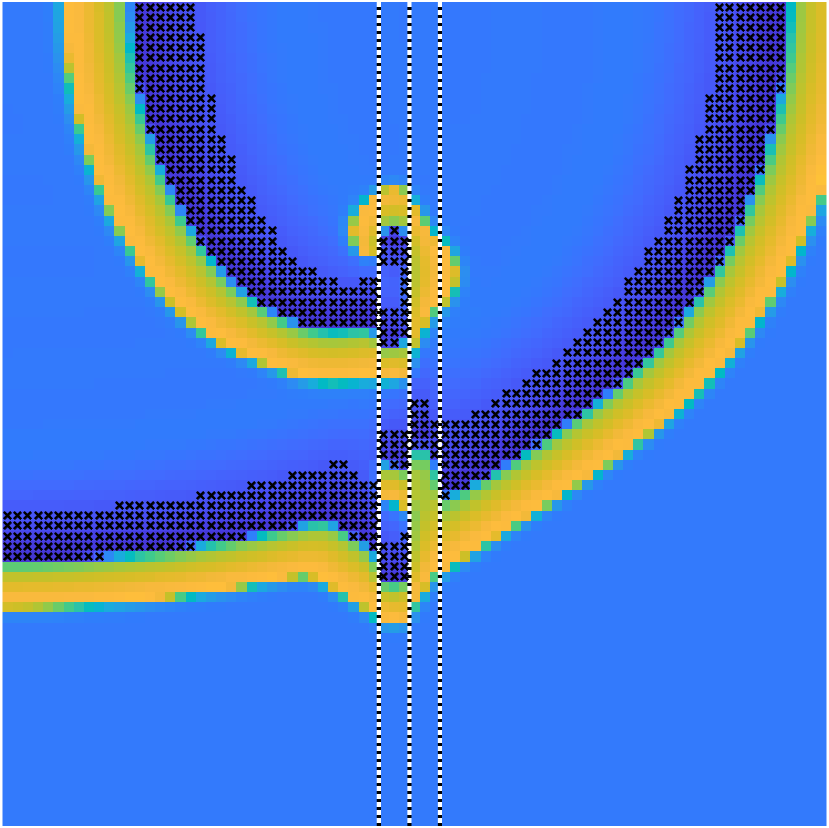}};
        \node[anchor=north] at (0,0) {$t=900$};
        &
        \node[anchor=south] at (0,0) {\includegraphics[width=\figscale]{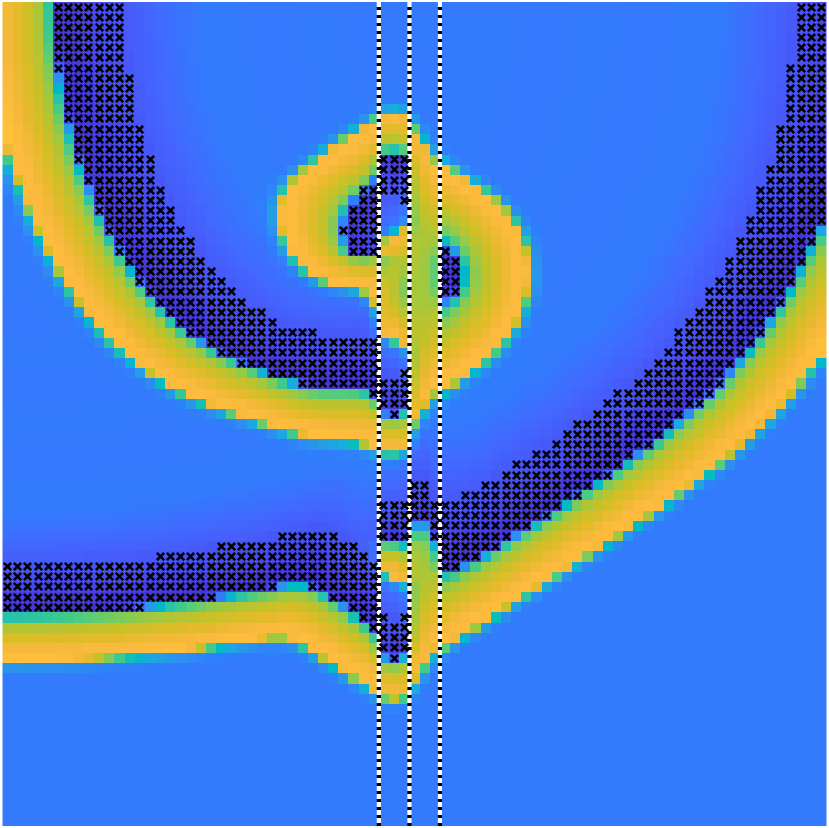}};
        \node[anchor=north] at (0,0) {$t=980$};
        &
        \node[anchor=south] at (0,0) {\includegraphics[width=\figscale]{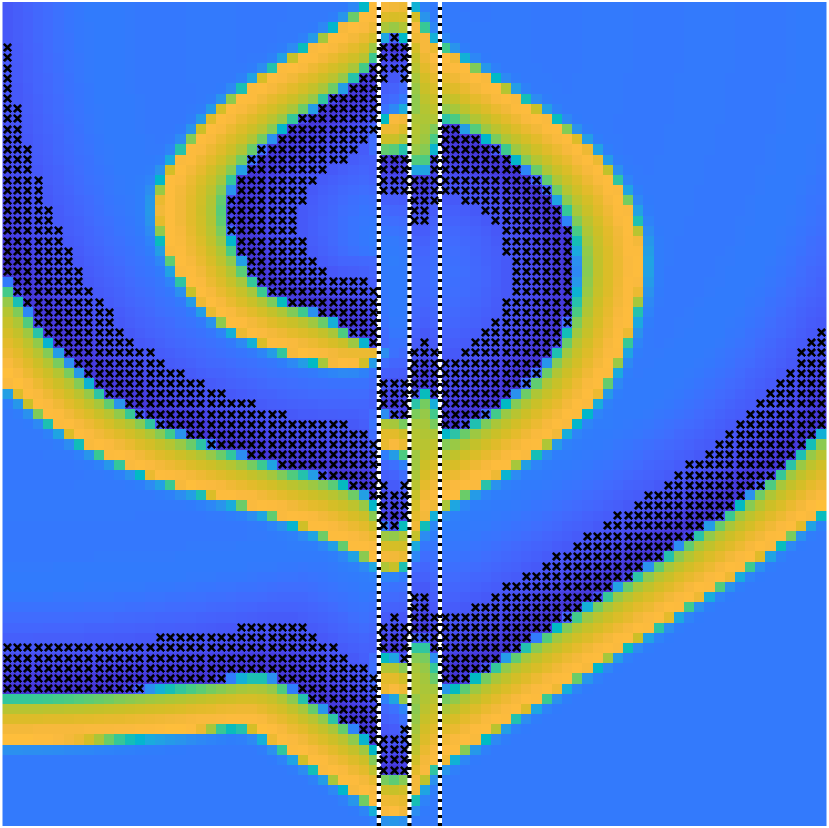}};
        \node[anchor=north] at (0,0) {$t=1100$};
        &
         \node[anchor=south] at (0,0) {\includegraphics[width=\figscale,height=\figscale]{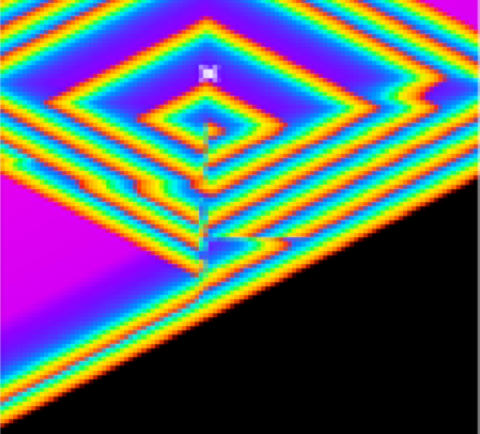}};
        \node[anchor=north] at (0,0) {Literature, re-oriented};
        \\
        };
    \end{tikzpicture}
    \caption{Snapshots of a simulation with a $6\times 81$ FIB in a  $81 \times 81$ computational grid, with fast-recovery cells with $\epsilon = 0.007$ and slow-recovery cells with $b = 4.9$, and an impulse sent from the upper-left corner. 
    Refirings occur both due to our usual FIB mechanism of slow-recovery cells retaining charge long enough that neighboring fast-recovery cells recover ($t=640$) and due to differences in wave speed creating broken ends ($t=820$).
    The final image is from \cite[Figure 5 part M]{A-A-K-Z-M-Betal:2019}, re-oriented to match our geometry.
    }
    \label{fig:vertical_strip_fib}
\end{figure}
At some times we see our usual FIB mechanism of slow-recovery cells retaining charge long enough that neighboring fast-recovery cells recover, thus initiating new waves propagating down and left.
At other times we see a different mechanism.
Since healthy cells have $\sigma=1$ and fast-recovery cells have $\sigma=1.6$, the voltage in fast-recovery cells changes faster and waves travels more quickly through them.
As a wave travels down through both healthy and fast-recovery cells simultaneously, the speed difference causes the wavefront to separate and create a broken end, which curls back and creates a locally spiral-like pattern.


\subsection{Phantom Fibrillatory Centers}

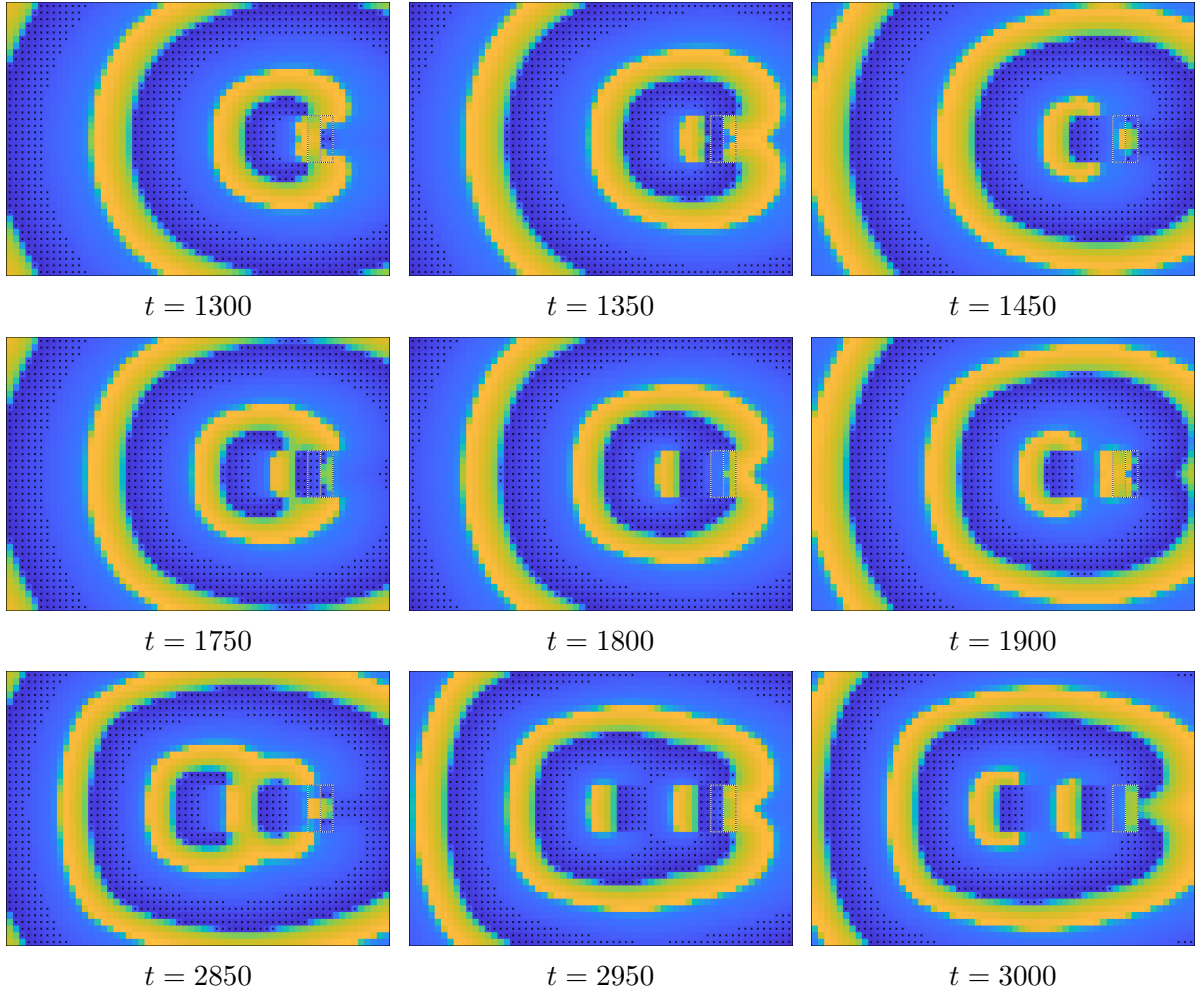
\begin{figure}
    \centering
    \begin{tikzpicture} 
    \setlength{\figscale}{0.22\textwidth}
    \matrix[row sep=1, column sep=0, cells={scale=1}] at (0,0)
        {
        \node[anchor=south] at (0,0) {\includegraphics[width=1.4\figscale,height=\figscale]{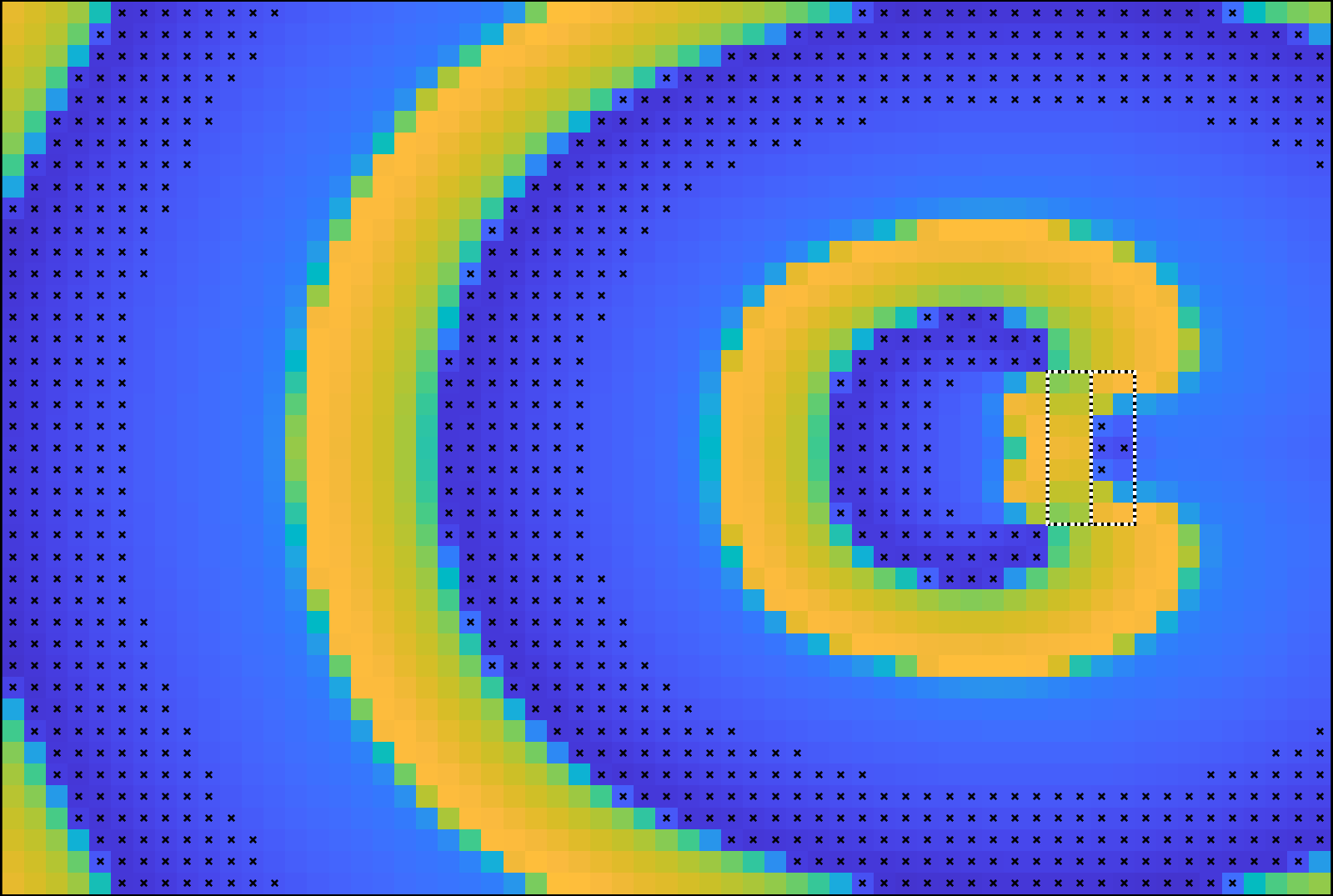}};
        \node[anchor=north] at (0,0) {$t=1300$};
        &
        \node[anchor=south] at (0,0) {\includegraphics[width=1.4\figscale,height=\figscale]{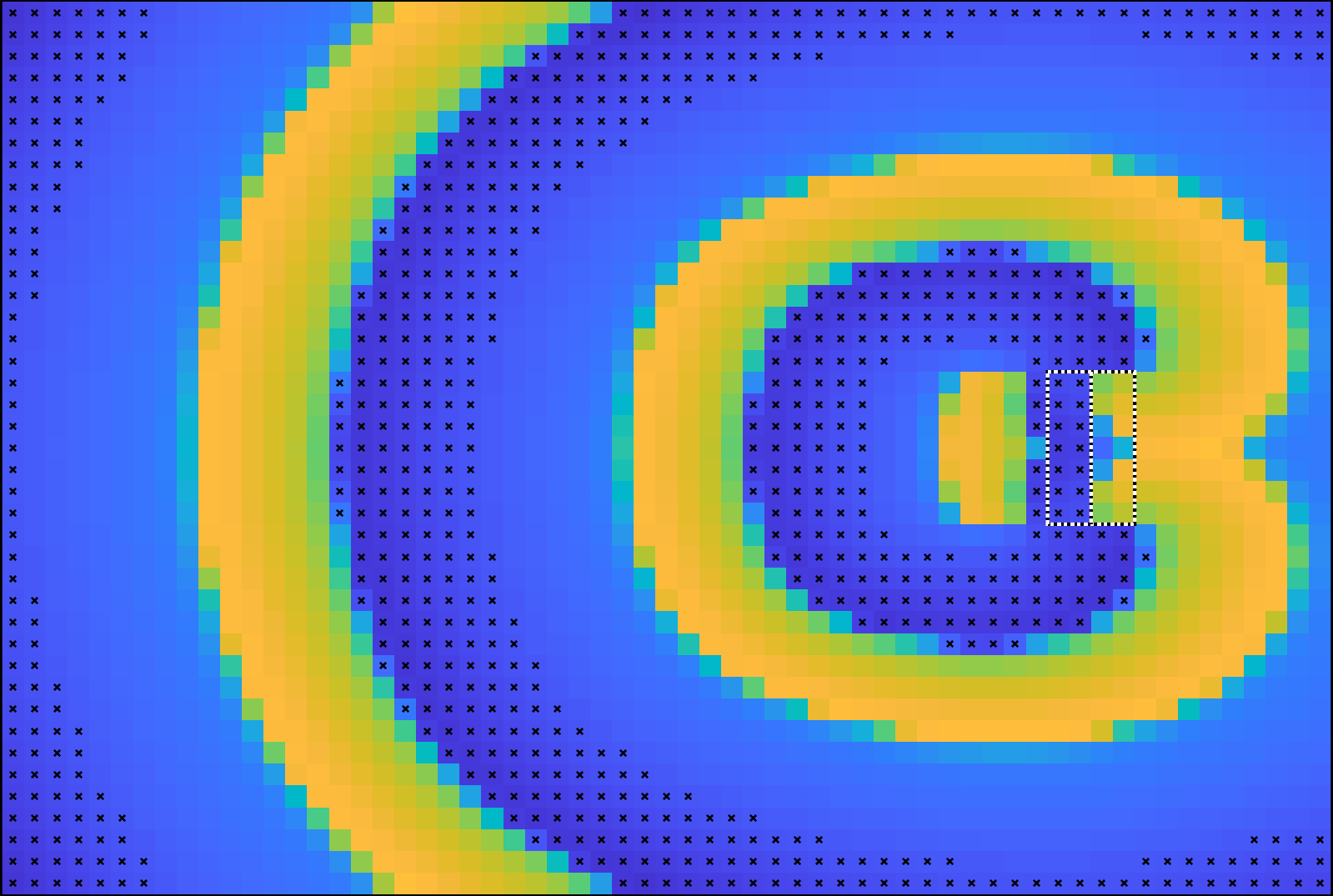}};
        \node[anchor=north] at (0,0) {$t=1350$};
        &
        \node[anchor=south] at (0,0) {\includegraphics[width=1.4\figscale,height=\figscale]{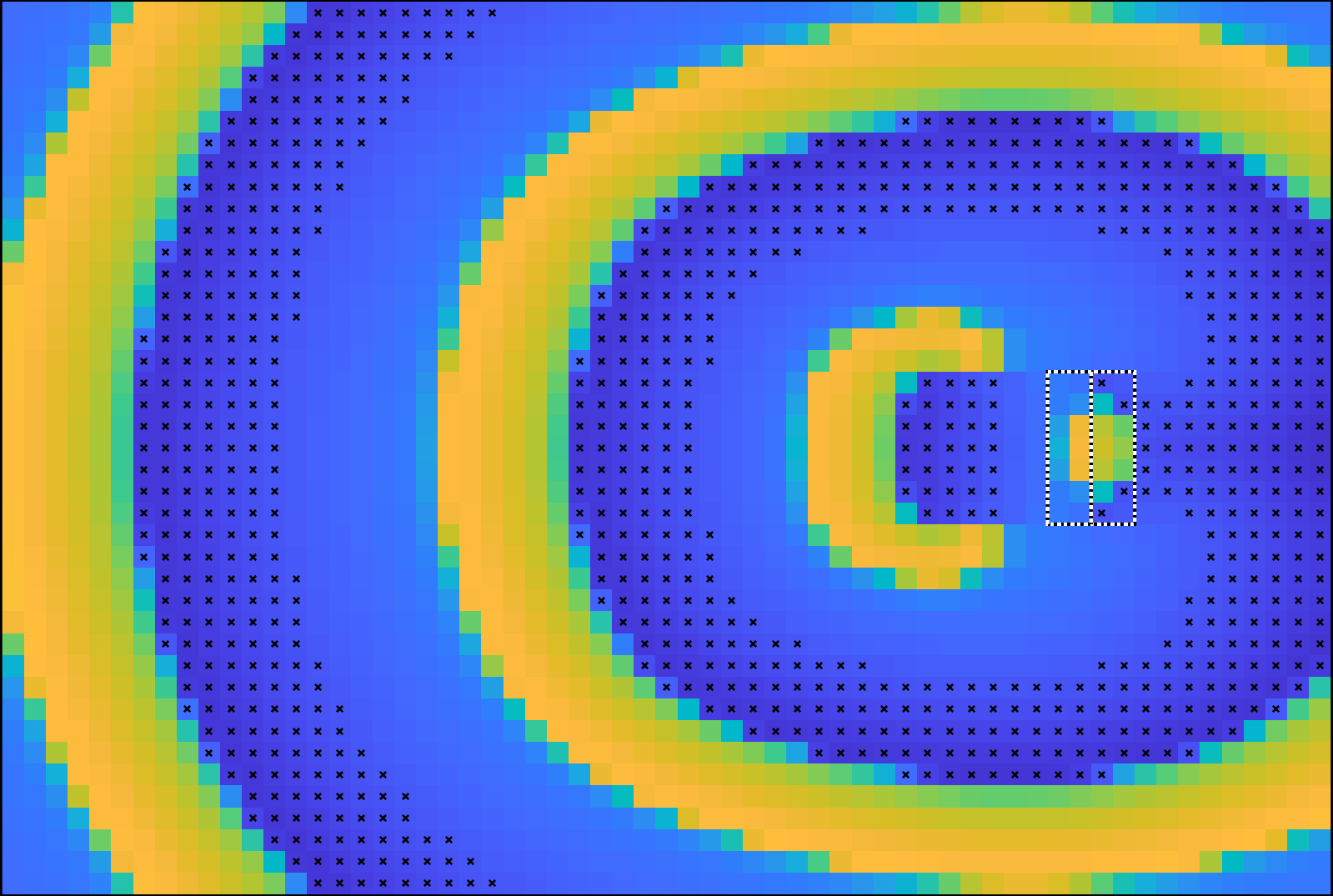}};
        \node[anchor=north] at (0,0) {$t=1450$};
    \\
        \node[anchor=south] at (0,0) {\includegraphics[width=1.4\figscale,height=\figscale]{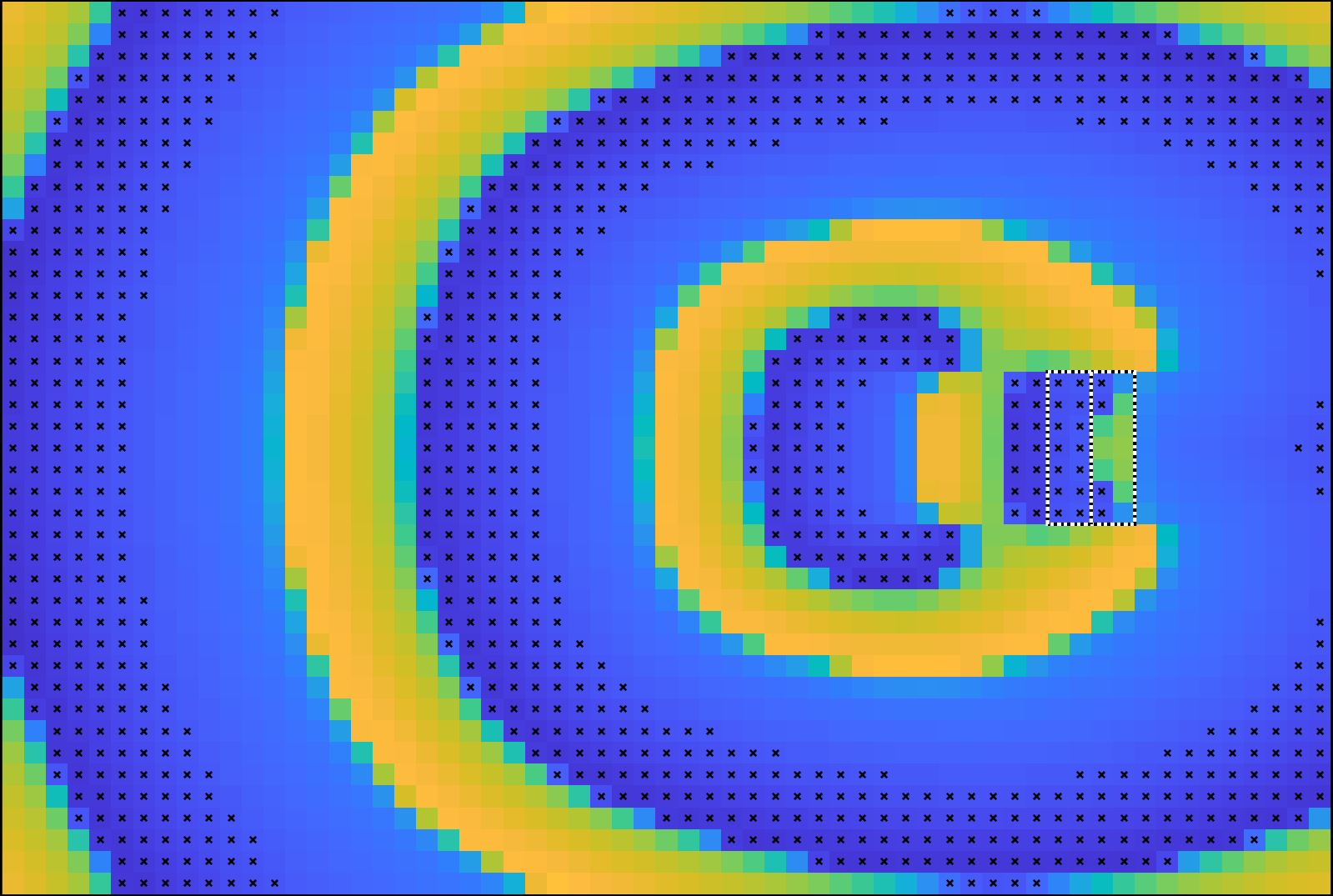}};
        \node[anchor=north] at (0,0) {$t=1750$};
        &
        \node[anchor=south] at (0,0) {\includegraphics[width=1.4\figscale,height=\figscale]{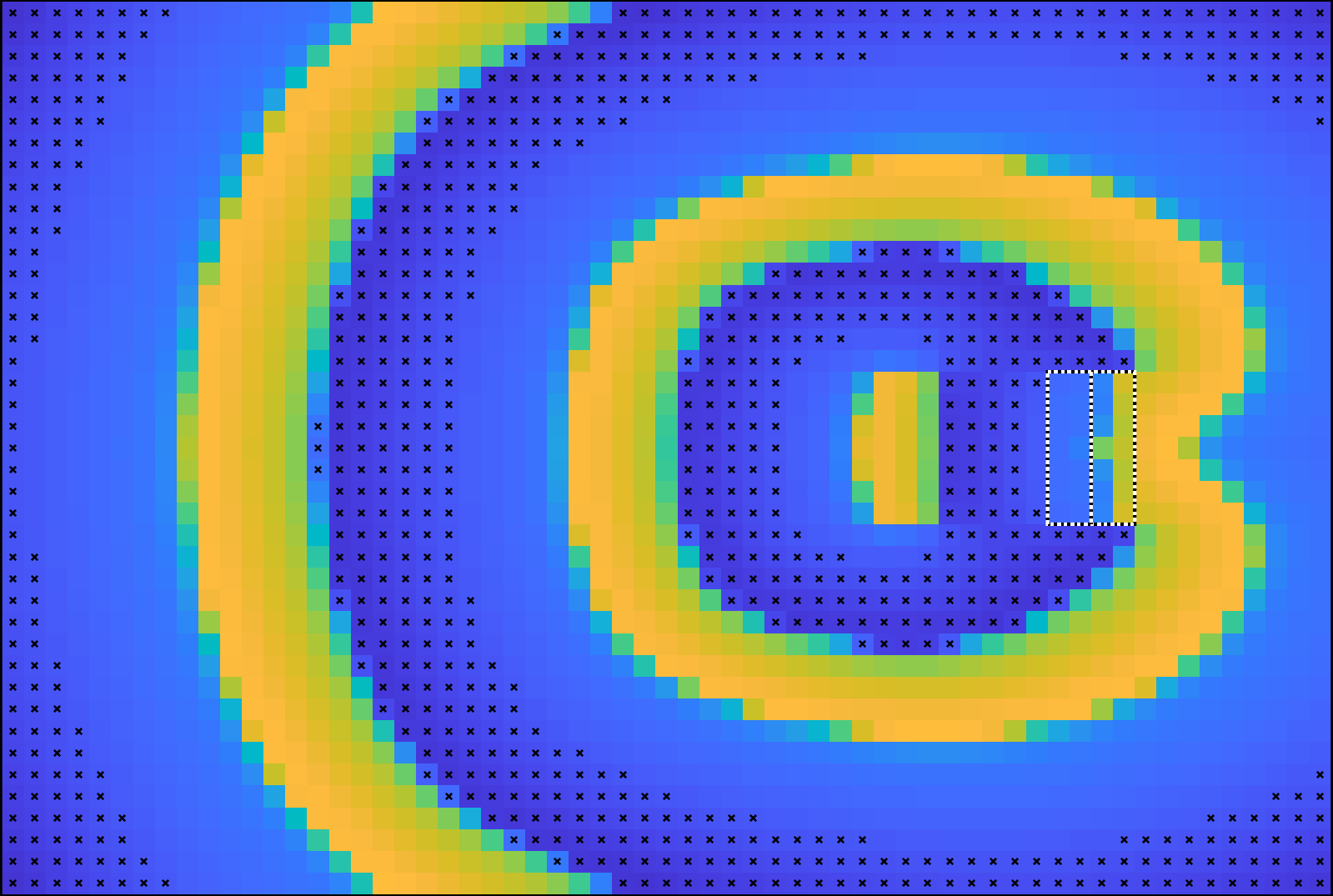}};
        \node[anchor=north] at (0,0) {$t=1800$};
        &
        \node[anchor=south] at (0,0) {\includegraphics[width=1.4\figscale,height=\figscale]{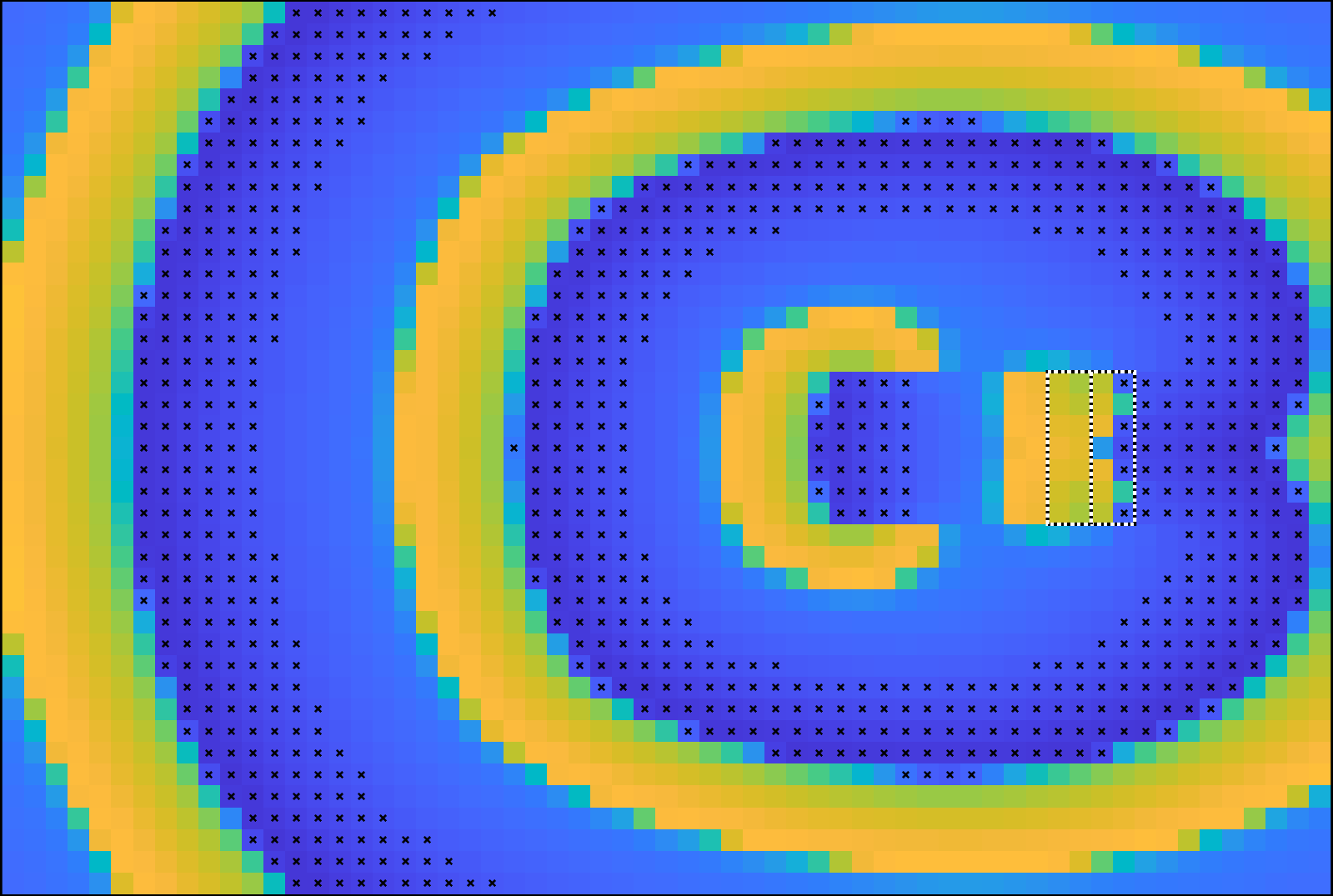}};
        \node[anchor=north] at (0,0) {$t=1900$};
    \\
        \node[anchor=south] at (0,0) {\includegraphics[width=1.4\figscale,height=\figscale]{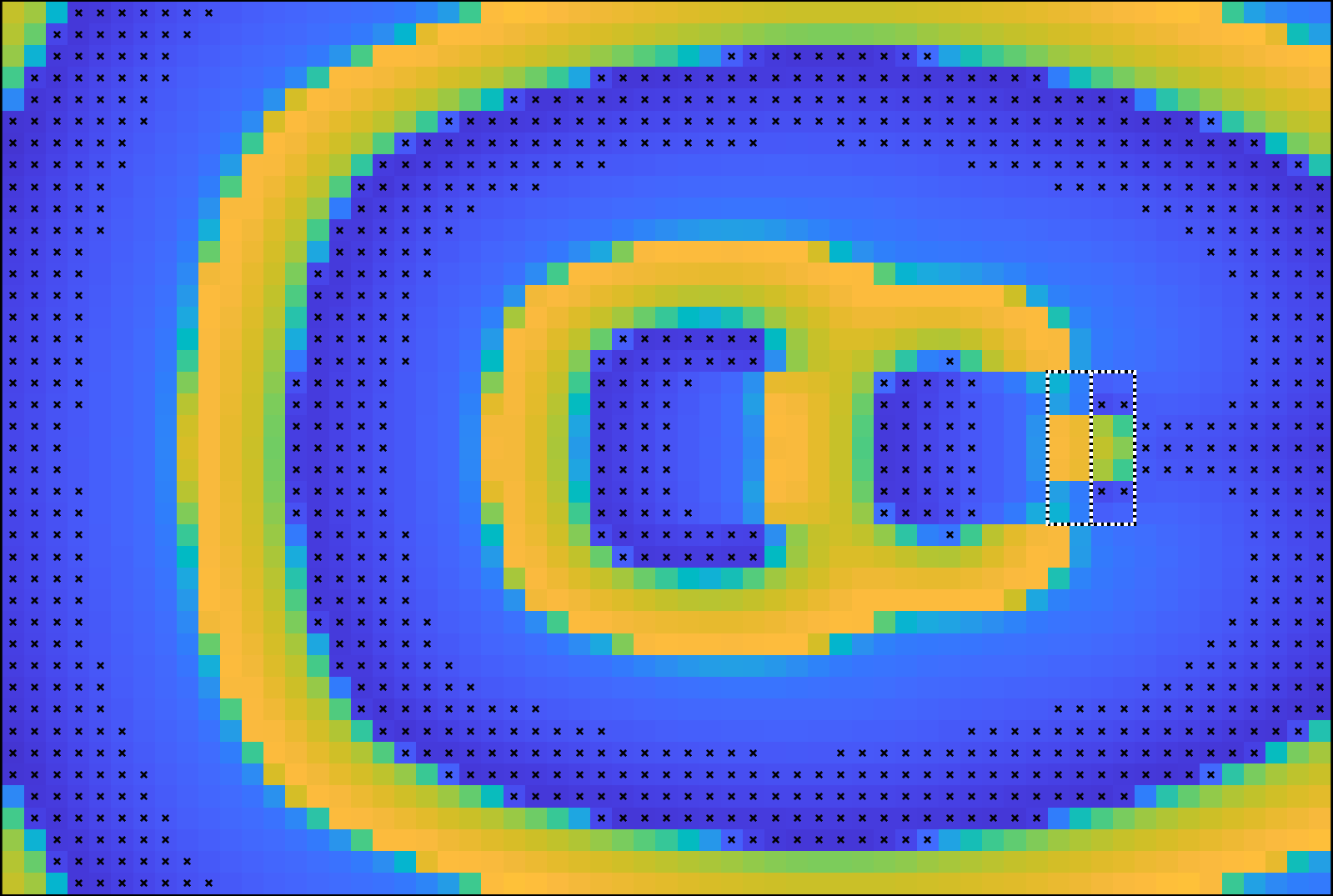}};
        \node[anchor=north] at (0,0) {$t=2850$};
        &
        \node[anchor=south] at (0,0) {\includegraphics[width=1.4\figscale,height=\figscale]{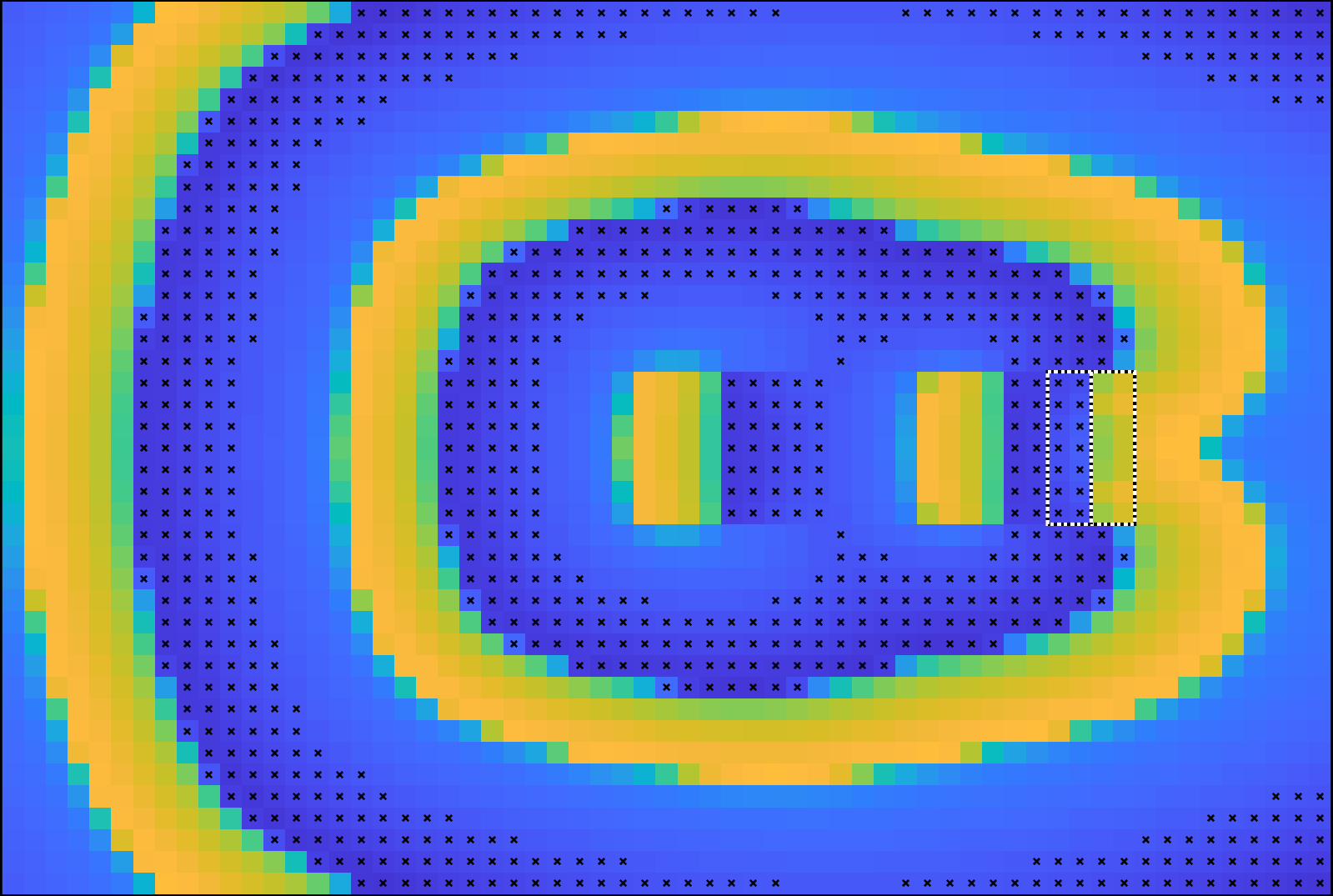}};
        \node[anchor=north] at (0,0) {$t=2950$};
        &
        \node[anchor=south] at (0,0) {\includegraphics[width=1.4\figscale,height=\figscale]{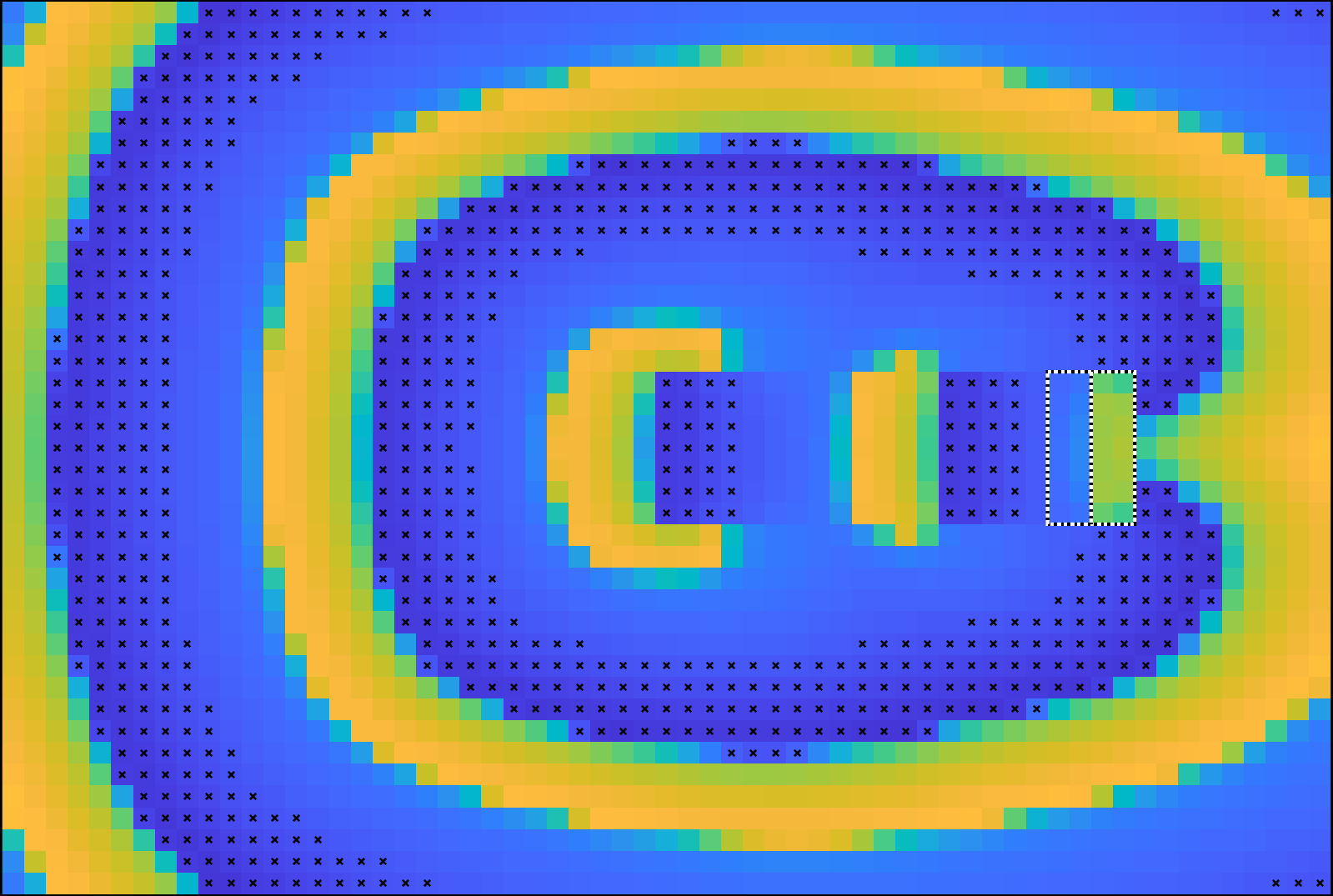}};
        \node[anchor=north] at (0,0) {$t=3000$};
    \\
        };
    \end{tikzpicture}
    \caption{Snapshots of a simulation demonstrating the creation of phantom FIBs. This simulation uses a typical $4 \times 7$ FIB with slow-recovery $b=5.55$ and fast-recovery $\epsilon=0.005$, within a $61 \times 41$ grid of healthy cells. We can observe the creation of dual broken ends on the small left-moving waves which spiral back and form double swirling-back curves.
    }
    \label{fig:phantomfib}
\end{figure}

In \cref{fig:phantomfib}, we show snapshots of a simulation with our standard $4\times 7$ FIB, but in a larger ($61\times 41$) grid and over a longer time.
A video is included in the Supplemental Material \cite{AF-supplement}.
In the first row, a small, directed wave propagates away from the FIB before exhibiting the swirling-back mechanism of a FIB.
Thus, a phantom FIB has been created entirely within healthy cells.
This effect can be compounded across repeated refirings, creating additional phantom FIBs far from the original FIB, as is shown in the second and third rows of \cref{fig:phantomfib}.

In all our simulations, conduction is anisotropic (stronger horizontally), with ratio 1.8.
Thus, generated waves move faster in the horizontal direction than in the vertical direction, so the wavefront is compressed vertically. 
Cells with $v<0$ and $w>0.0225$ are marked with $\times$ to indicate that they are inhibited from firing, but this is an arbitrary cutoff; any residual $v<0$ and $w>0$ acts against firing and tends to delay it.
Each wavefront encounters the aftereffects of the previous wave, becoming even more compressed vertically,  in a positive feedback loop. 
Eventually a wavefront is unable to cause the cells above and below it to fire at all, and becomes a directed wave that can initiate a phantom FIB. 


Within this understanding of the mechanism for creating the phantom FIB, the role of the actual FIB is only to create repeated stimulations.
To verify this understanding, we ran a simulation with repeated external stimulations (pacing) of healthy cells, as shown in \cref{fig:phantomfib_paced}.
A video is included in the Supplemental Material \cite{AF-supplement}.
\begin{figure}
    \centering
    \begin{tikzpicture} 
    \setlength{\figscale}{0.22\textwidth}
    \matrix[row sep=1, column sep=0, cells={scale=1}] at (0,0)
        {
        \node[anchor=south] at (0,0) {\includegraphics[width=\figscale,height=\figscale]{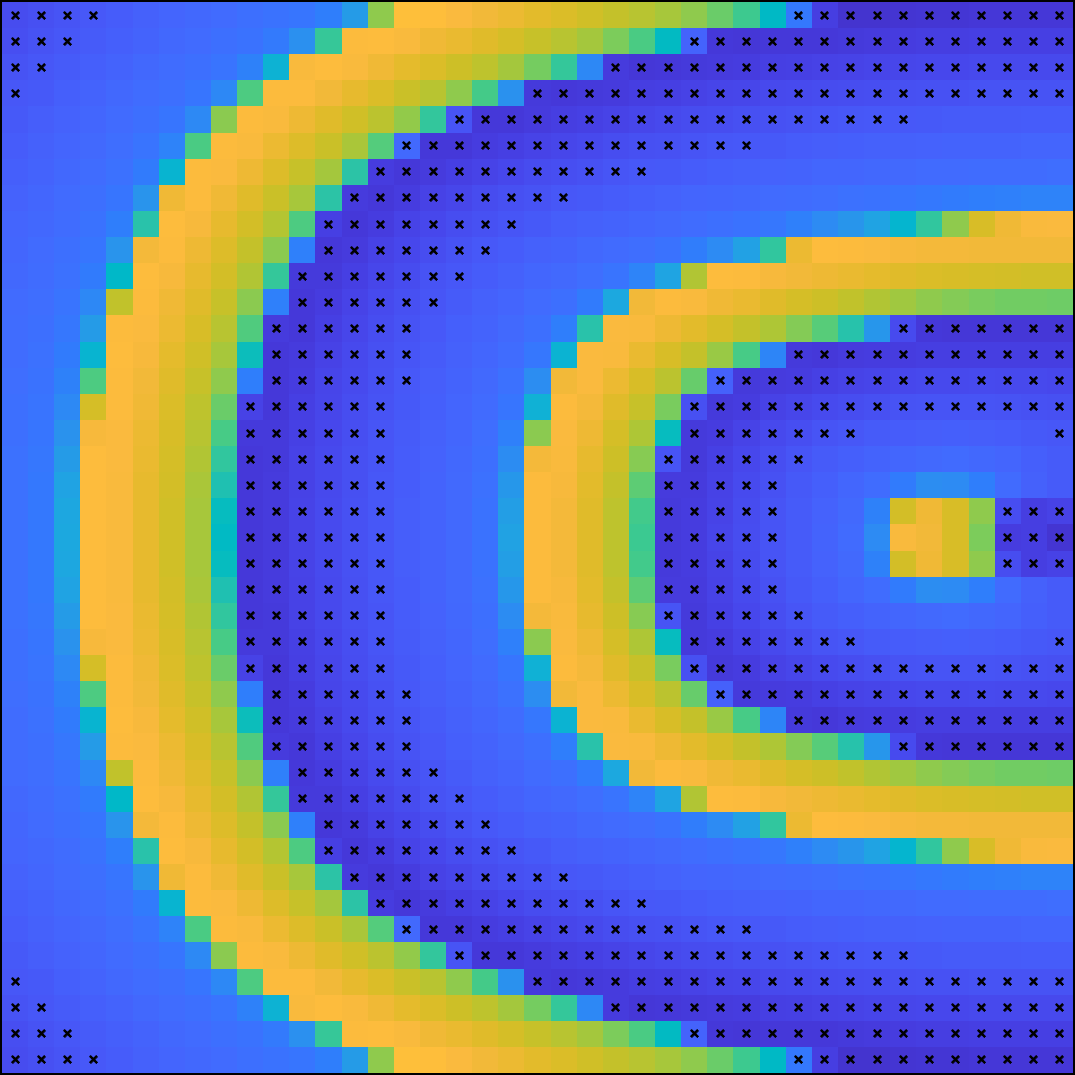}};
        \node[anchor=north] at (0,0) {$t=1300$};
        &
        \node[anchor=south] at (0,0) {\includegraphics[width=\figscale,height=\figscale]{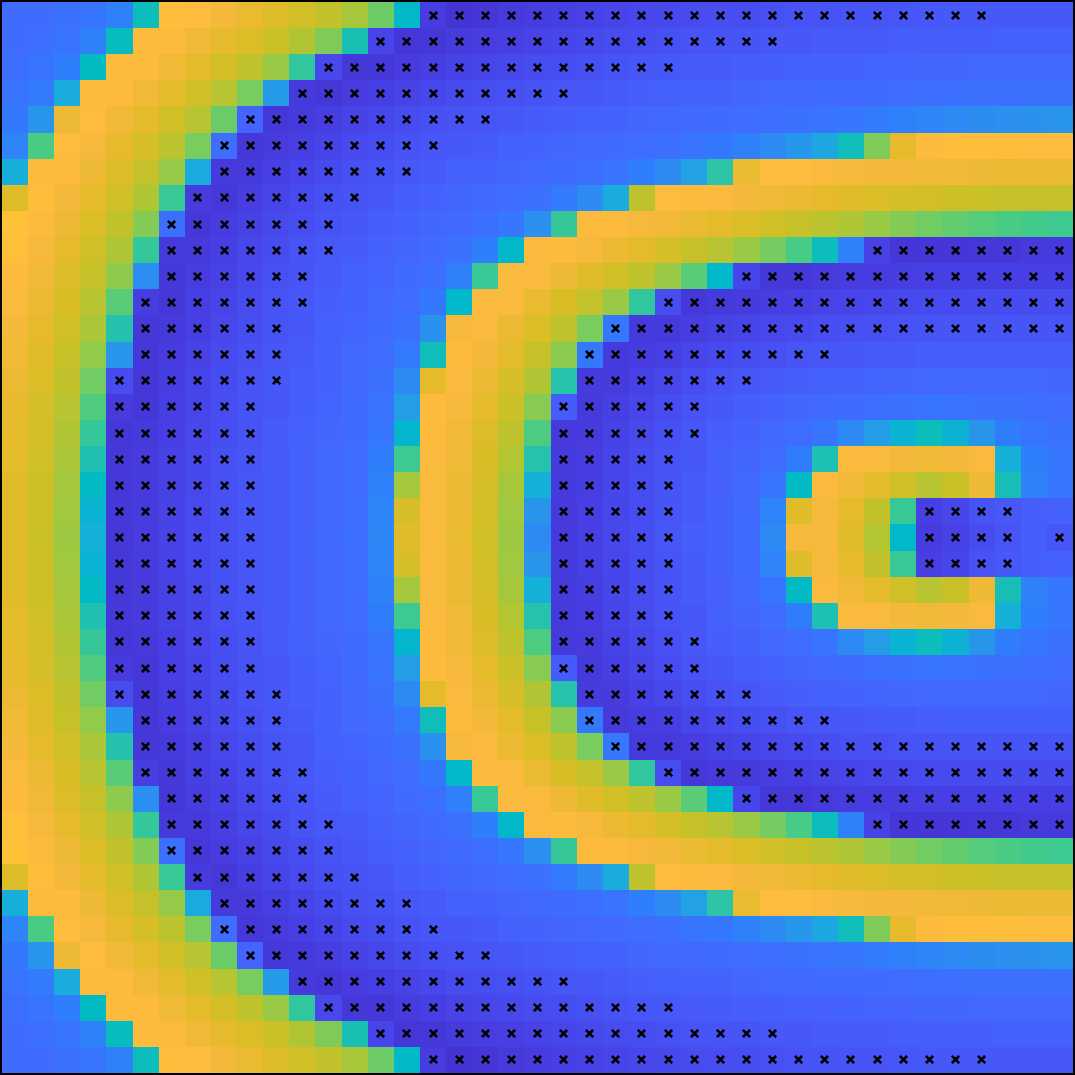}};
        \node[anchor=north] at (0,0) {$t=1360$};
        &
        \node[anchor=south] at (0,0) {\includegraphics[width=\figscale,height=\figscale]{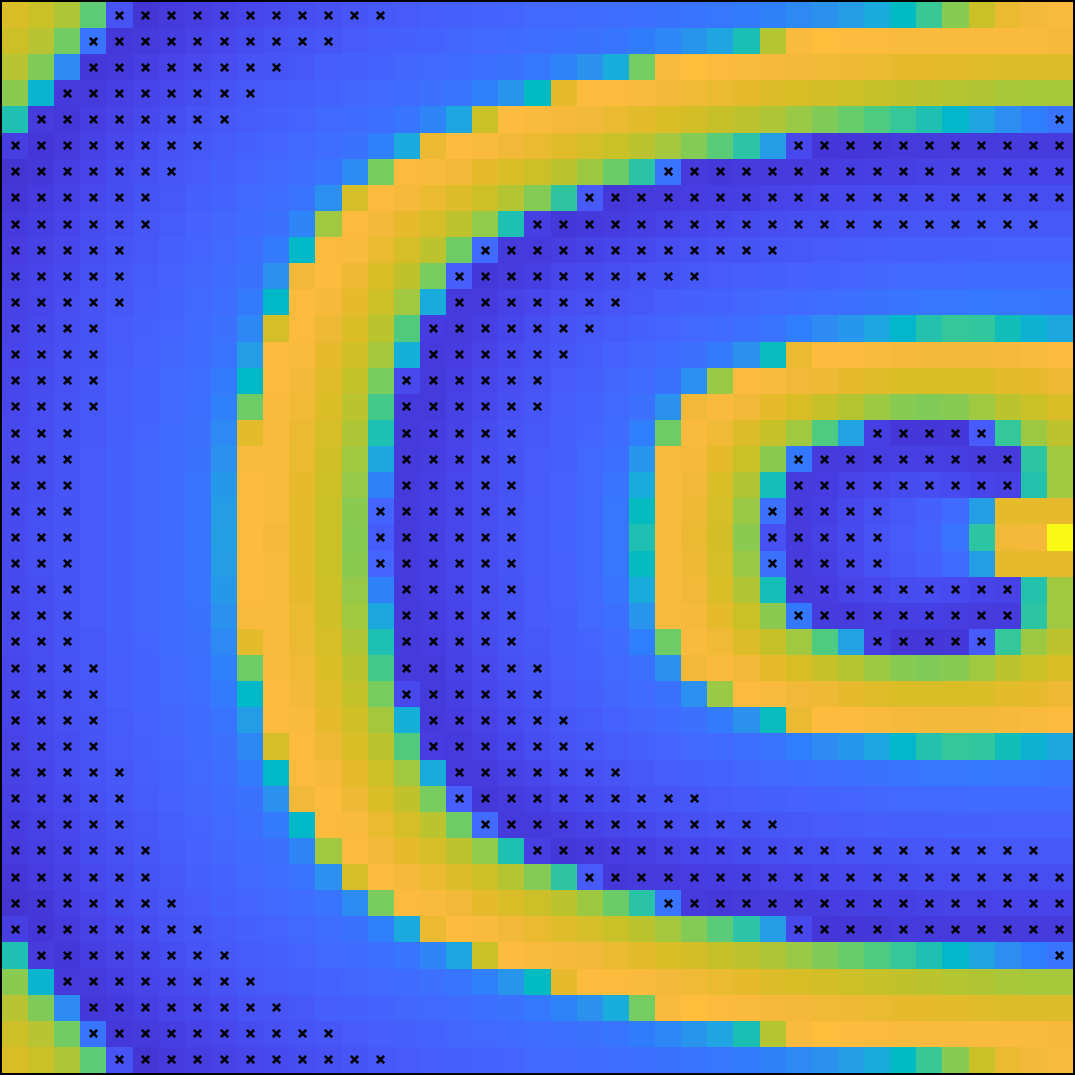}};
        \node[anchor=north] at (0,0) {$t=1440$};
        &
        \node[anchor=south] at (0,0) {\includegraphics[width=\figscale,height=\figscale]{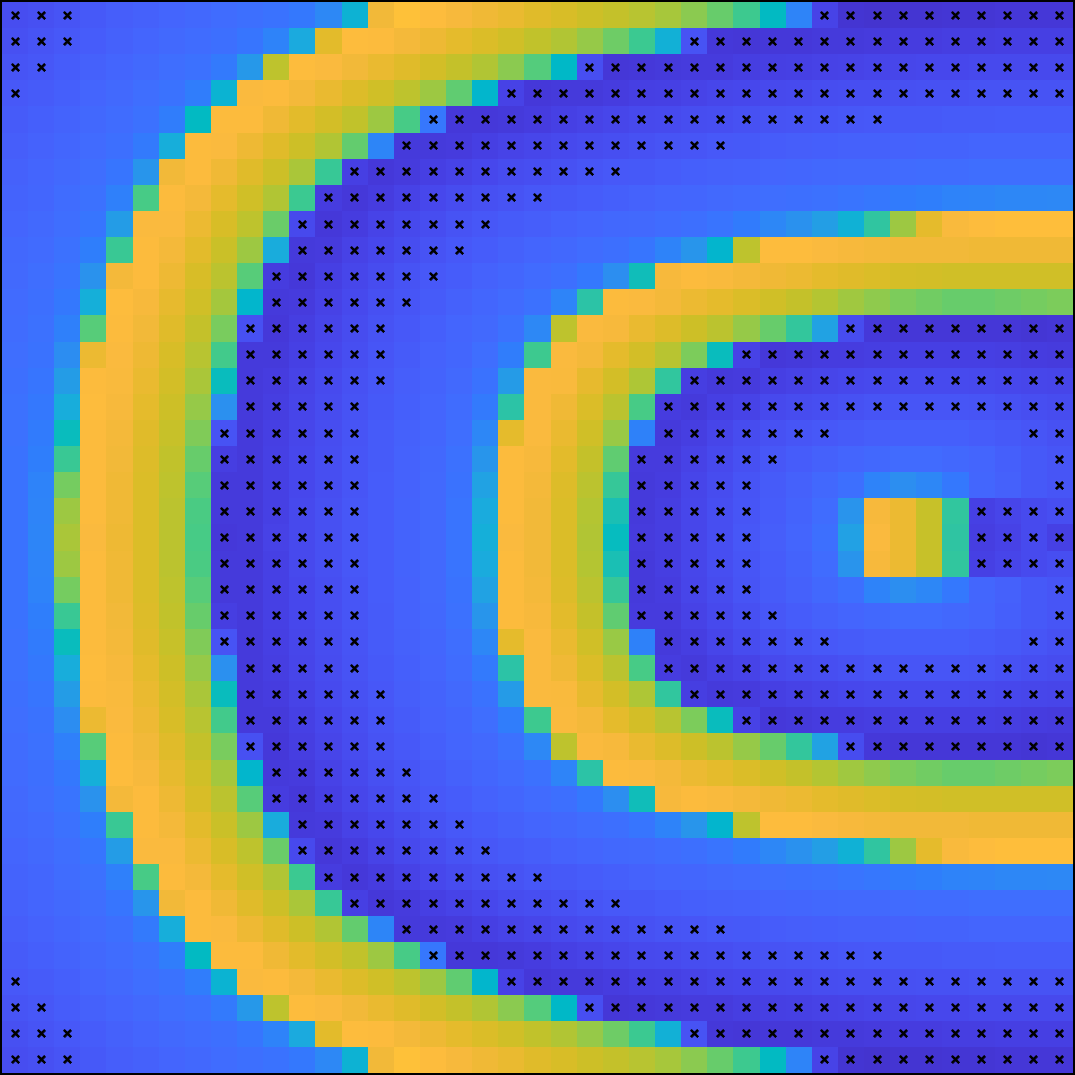}};
        \node[anchor=north] at (0,0) {$t=1520$};
    \\
        \node[anchor=south] at (0,0) {\includegraphics[width=\figscale,height=\figscale]{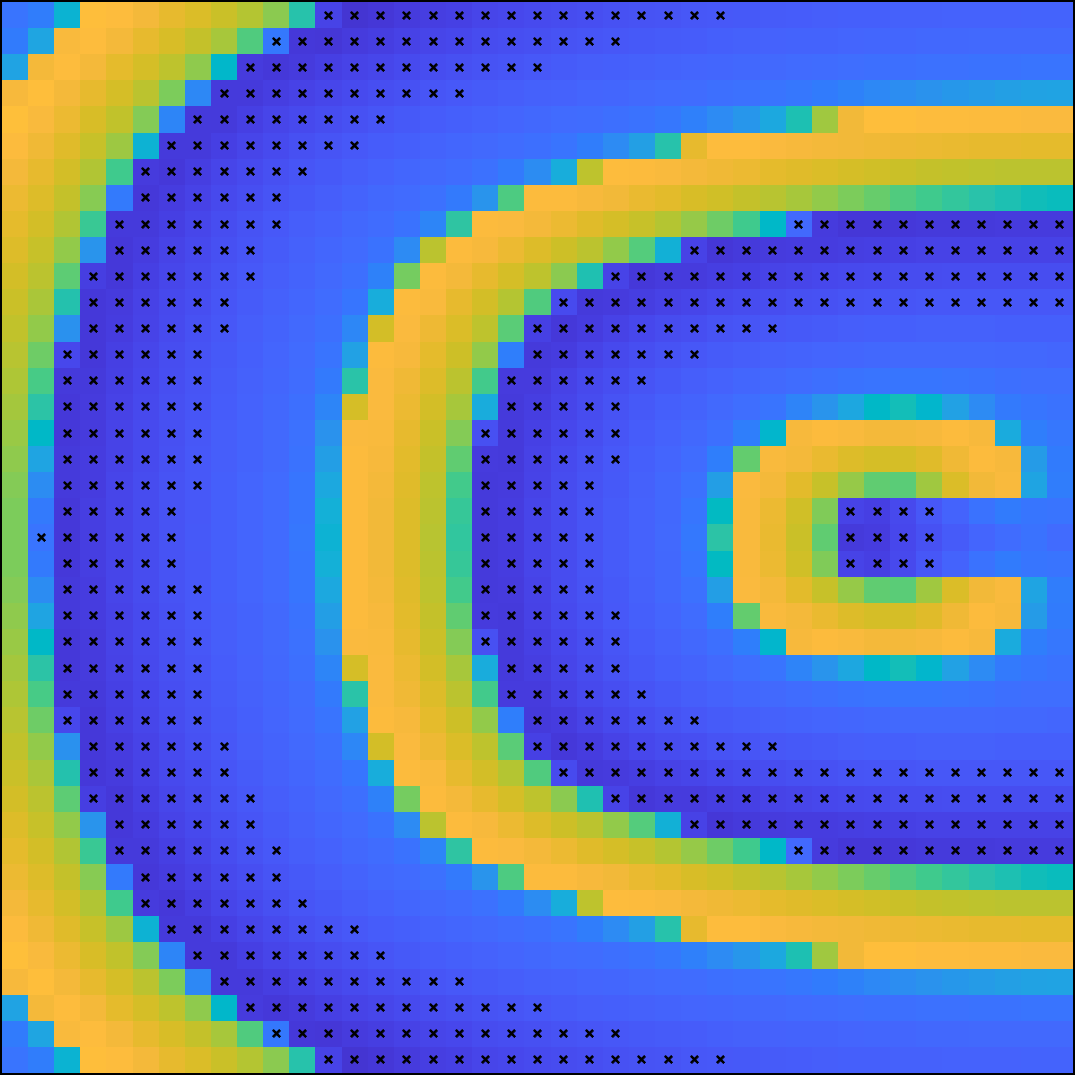}};
        \node[anchor=north] at (0,0) {$t=1600$};
        &
        \node[anchor=south] at (0,0) {\includegraphics[width=\figscale,height=\figscale]{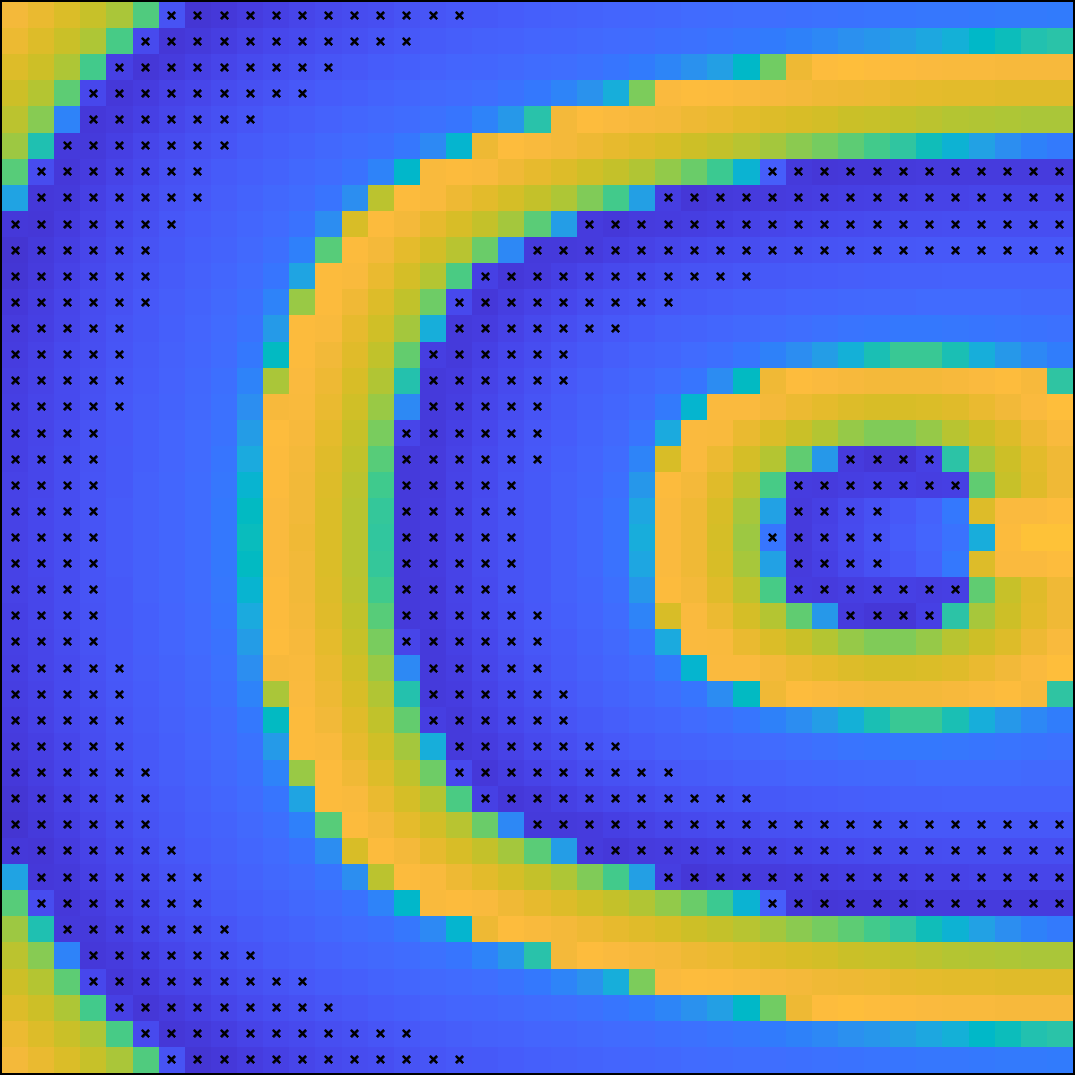}};
        \node[anchor=north] at (0,0) {$t=1640$};
        &
        \node[anchor=south] at (0,0) {\includegraphics[width=\figscale,height=\figscale]{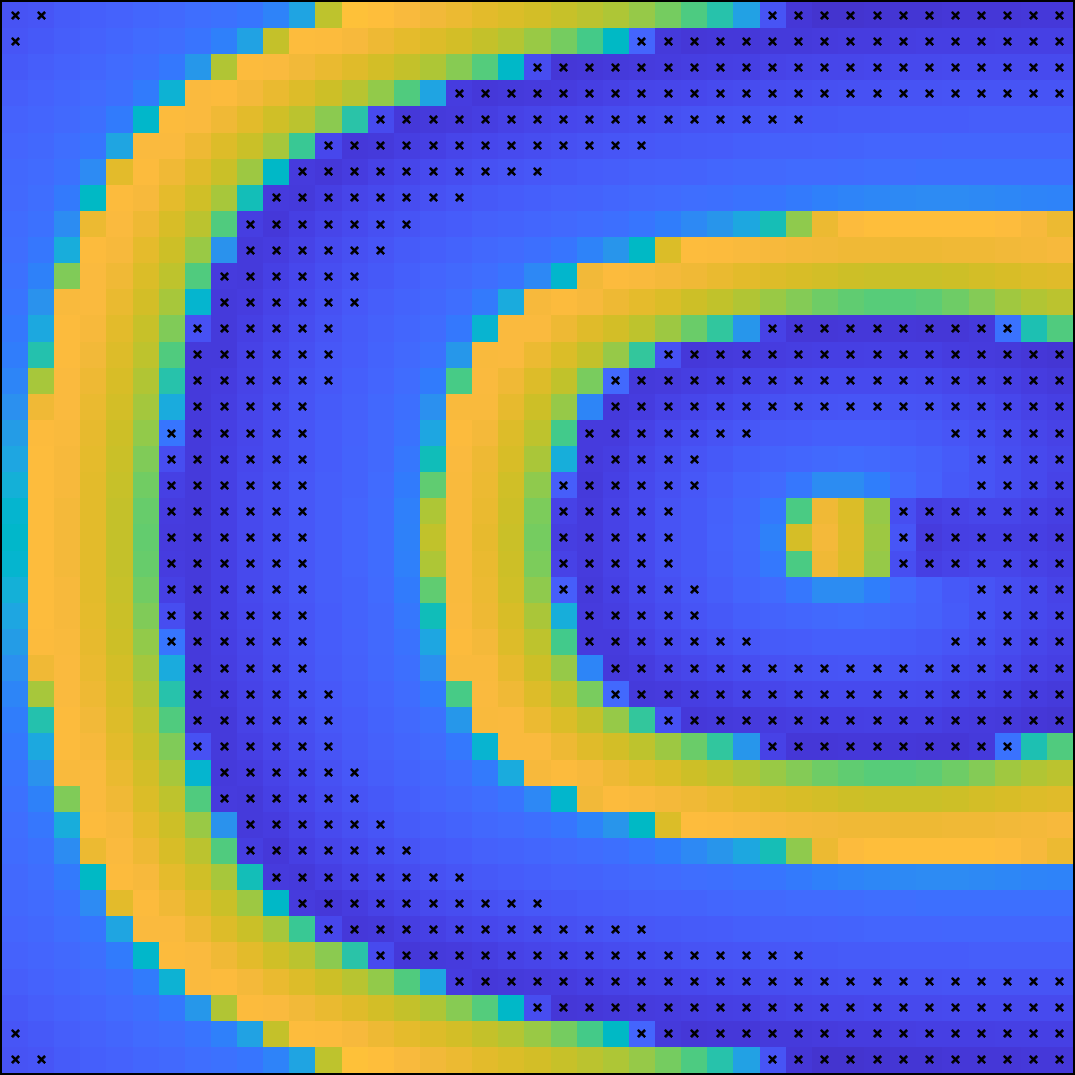}};
        \node[anchor=north] at (0,0) {$t=1960$};
        &
        \node[anchor=south] at (0,0) {\includegraphics[width=\figscale,height=\figscale]{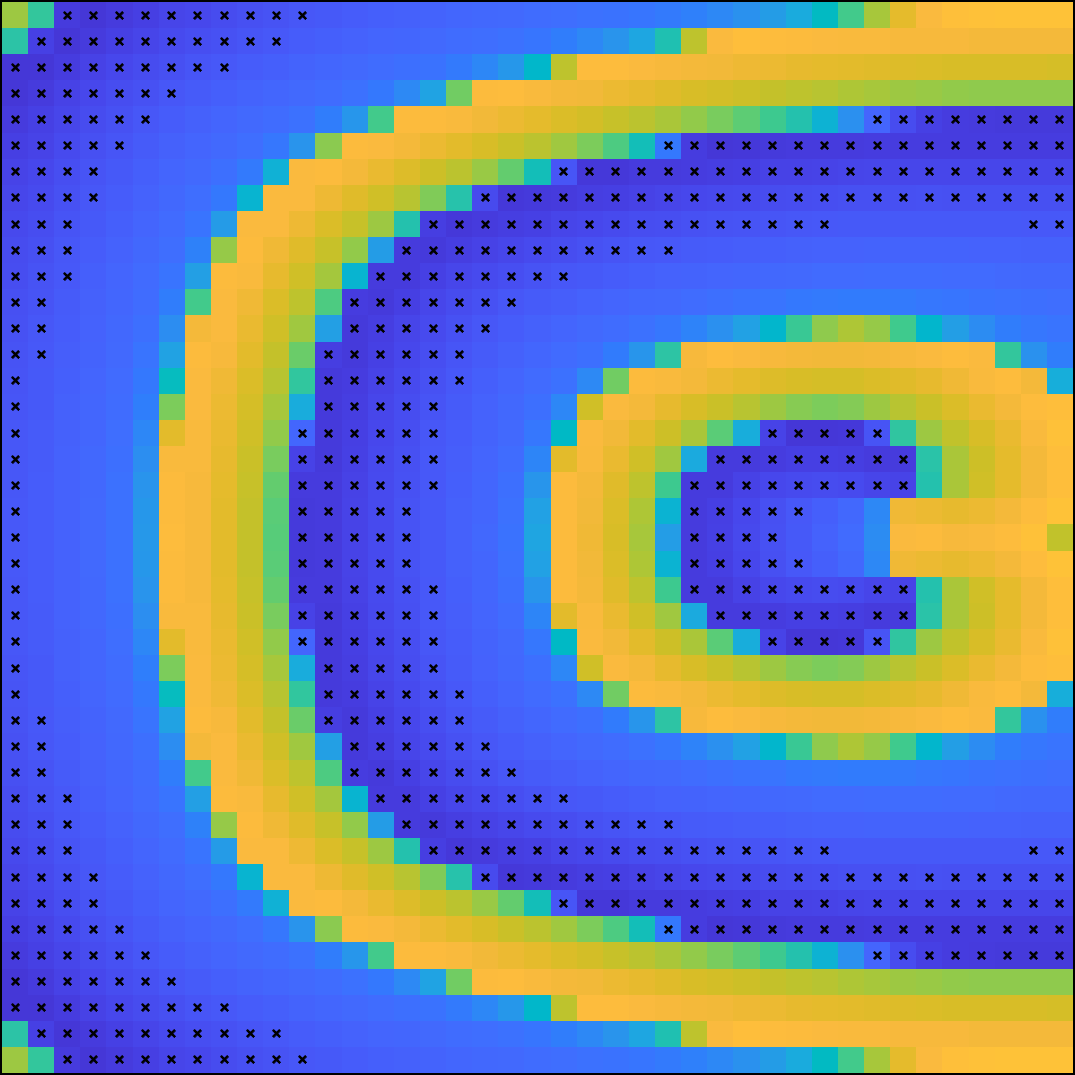}};
        \node[anchor=north] at (0,0) {$t=2100$};
    \\
        };
    \end{tikzpicture}
    \caption{Snapshots of a simulation demonstrating the creation of phantom FIBs via  external pacing. This simulated tissue is paced by applying a voltage to a single cell on the right edge for the first 50 of every 200 time units. The pacing is ended with the stimulus between 1400 and 1450 time units.
    }
    \label{fig:phantomfib_paced}
\end{figure}
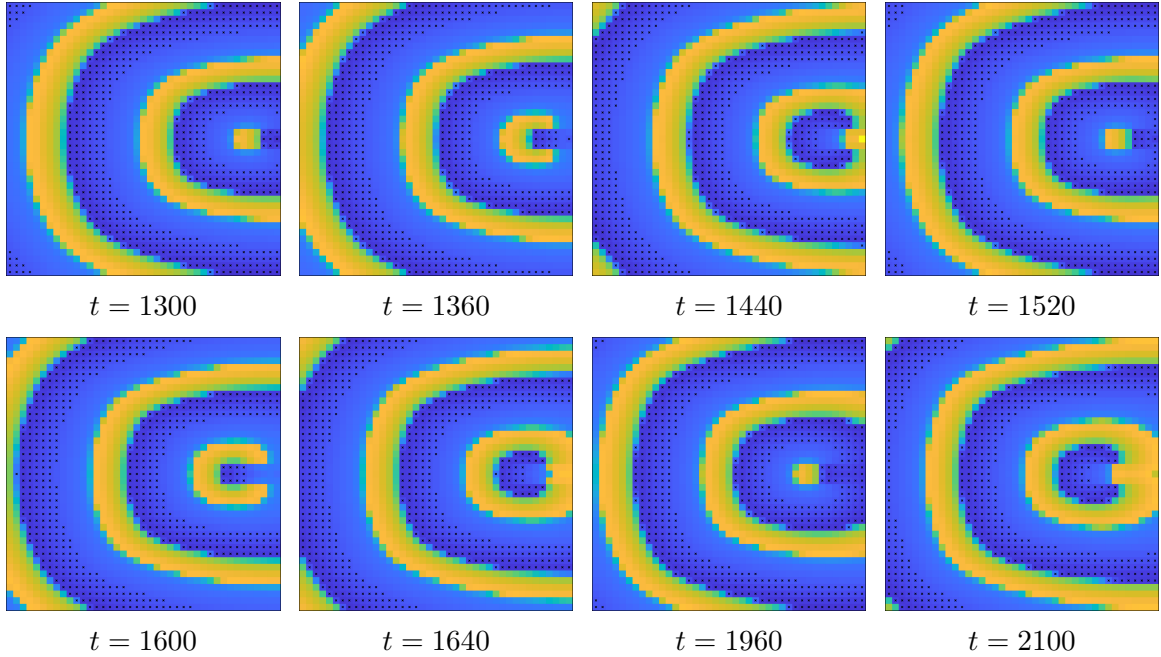
We observe the generation of a phantom fibrillatory center within these healthy cells and that this center is self-sustaining when the external stimulation ceases. 
Note that our FitzHugh-Nagumo model does not include changes in recovery time due to the frequency of stimulation (restitution effects \cite{MO-RH-AB:1964,NOL-DAH:1968}) and the pacing frequency is similar to the natural frequency of the FIB (as opposed to rapid pacing \cite{D-F-V-J-W-C-K-R:2005,O-H-K-L-L-W-C-K:2007}). 

\subsection{Creating Multiple Interacting Wave(let)s using Multiple Small FIBs}

When two FIBs are placed near each other, we observe several different types of interactions between them. 
Sometimes both FIBs are firing, so we observe colliding wave fronts.
Sometimes, one FIB will dominate the other and have greater influence over the general shape of the outgoing wave.
Sustained FIBs with larger fast-recovery $\epsilon$ values tend to dominate by generating waves at higher frequencies.  
Finally, sometimes the dominant FIB ceases firing (spontaneously) and then the other FIB takes over and begins emitting waves, now becoming the dominant FIB. 
Each of these types of interactions can be observed in simulations, as seen in \cref{fig:two_fibs_interact}.
A video is included in the Supplemental Material \cite{AF-supplement}.
\begin{figure}
    \centering
    \begin{tikzpicture} 
    \setlength{\figscale}{0.23\textwidth}
    \matrix[row sep=1, column sep=0, cells={scale=1}] at (0,0)
        {
        \node[anchor=south] at (0,0) {\includegraphics[width=\figscale,height=\figscale]{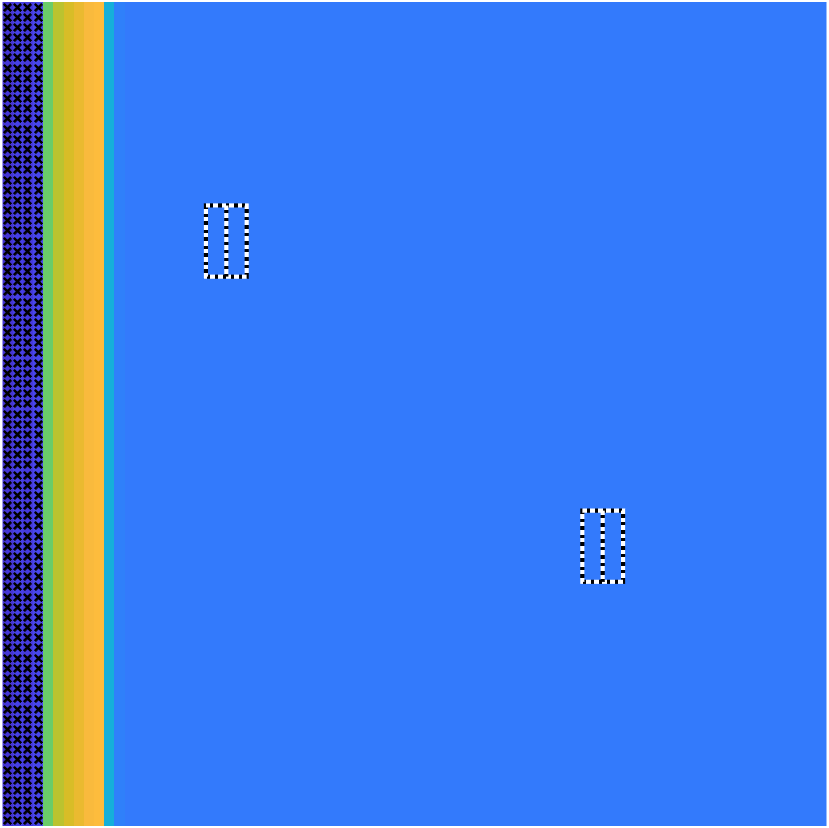}};
        \node[anchor=north] at (0,0) {$t=100$};
        &
        \node[anchor=south] at (0,0) {\includegraphics[width=\figscale,height=\figscale]{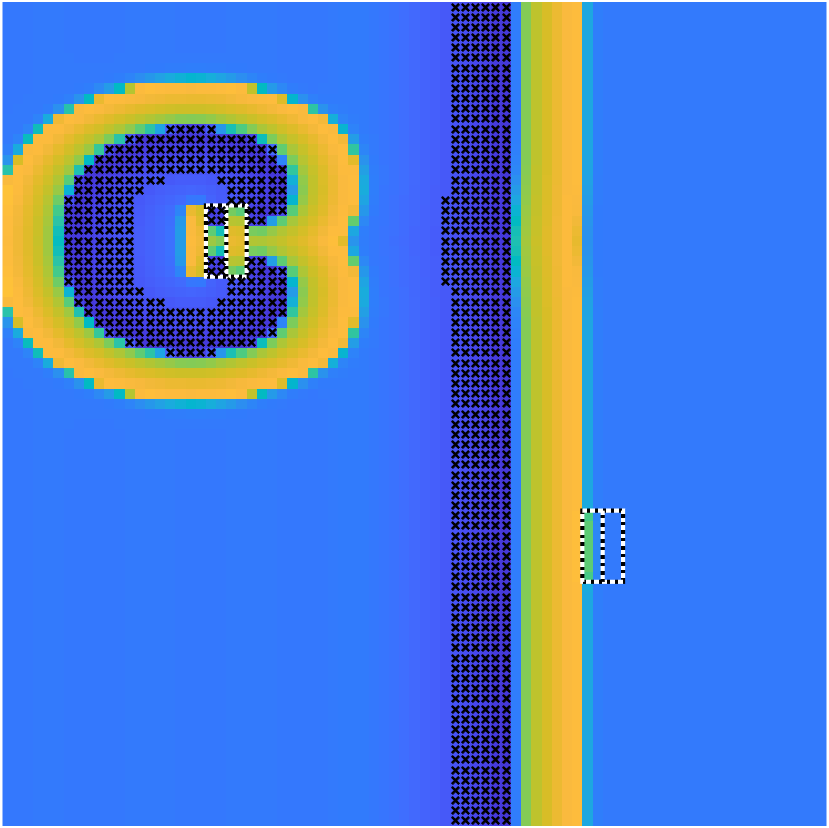}};
        \node[anchor=north] at (0,0) {$t=600$};
        &
        \node[anchor=south] at (0,0) {\includegraphics[width=\figscale,height=\figscale]{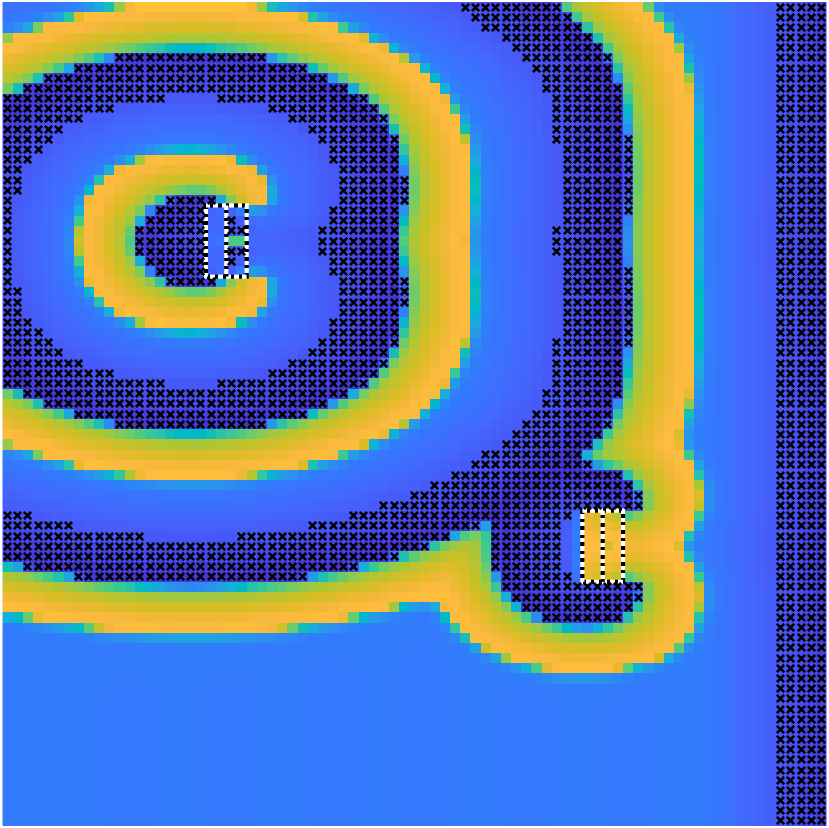}};
        \node[anchor=north] at (0,0) {$t=950$};
        &
        \node[anchor=south] at (0,0) {\includegraphics[width=\figscale,height=\figscale]{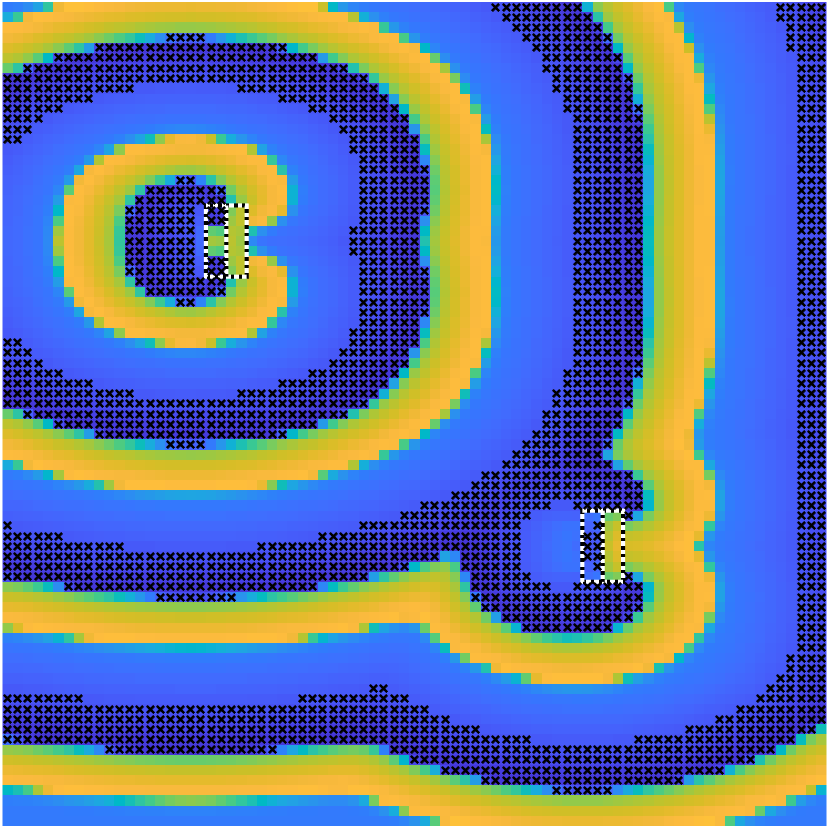}};
        \node[anchor=north] at (0,0) {$t=1200$};
        \\
        \node[anchor=south] at (0,0) {\includegraphics[width=\figscale,height=\figscale]{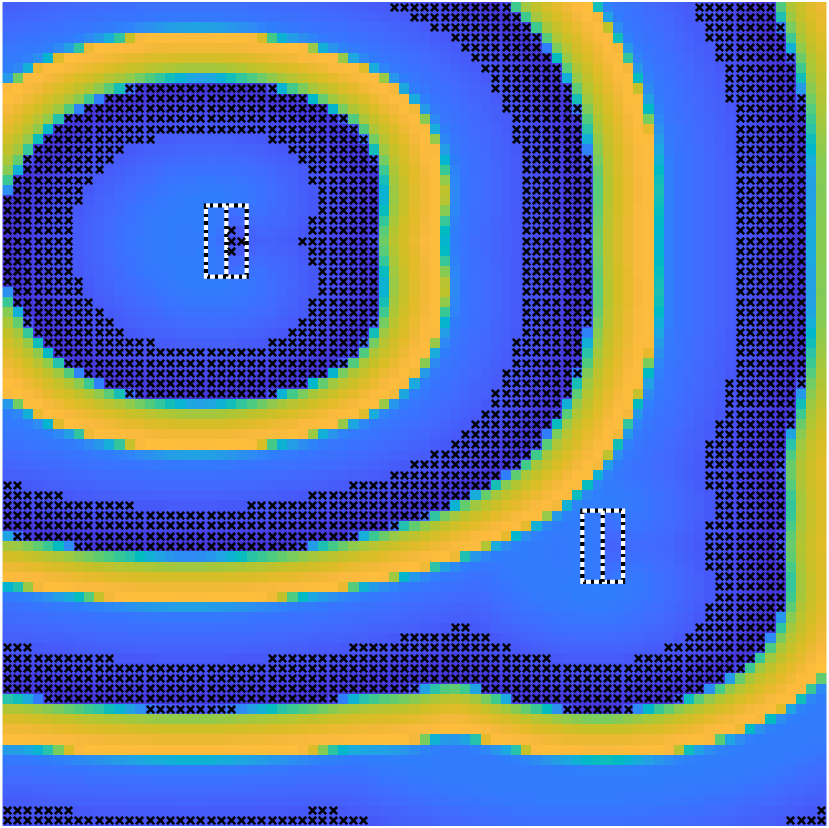}};
        \node[anchor=north] at (0,0) {$t=1600$};
        &
        \node[anchor=south] at (0,0) {\includegraphics[width=\figscale,height=\figscale]{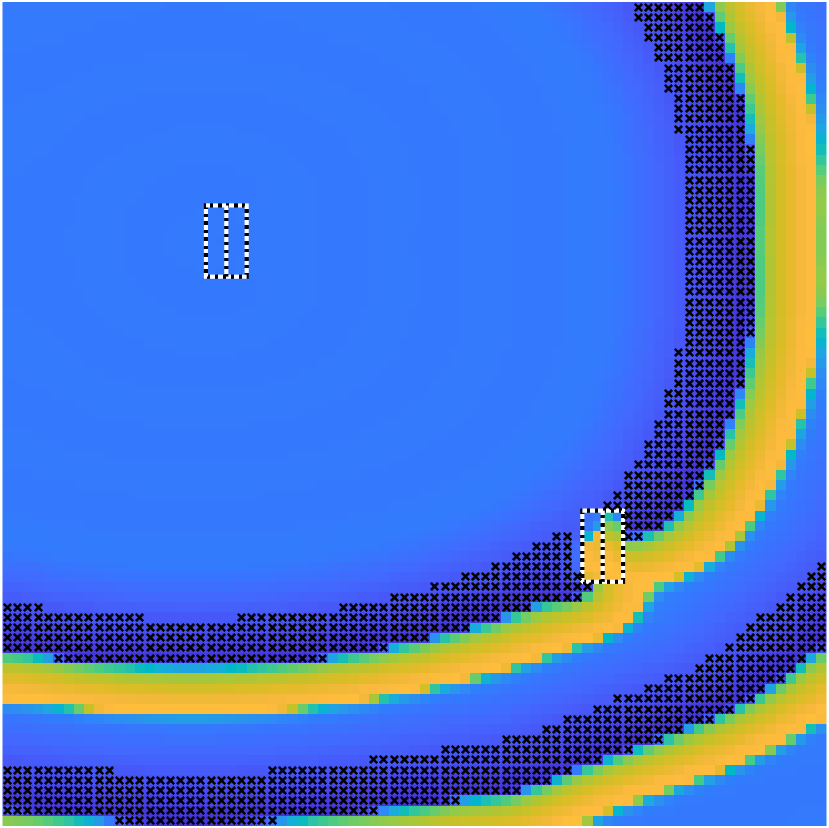}};
        \node[anchor=north] at (0,0) {$t=2000$};
        &
        \node[anchor=south] at (0,0) {\includegraphics[width=\figscale,height=\figscale]{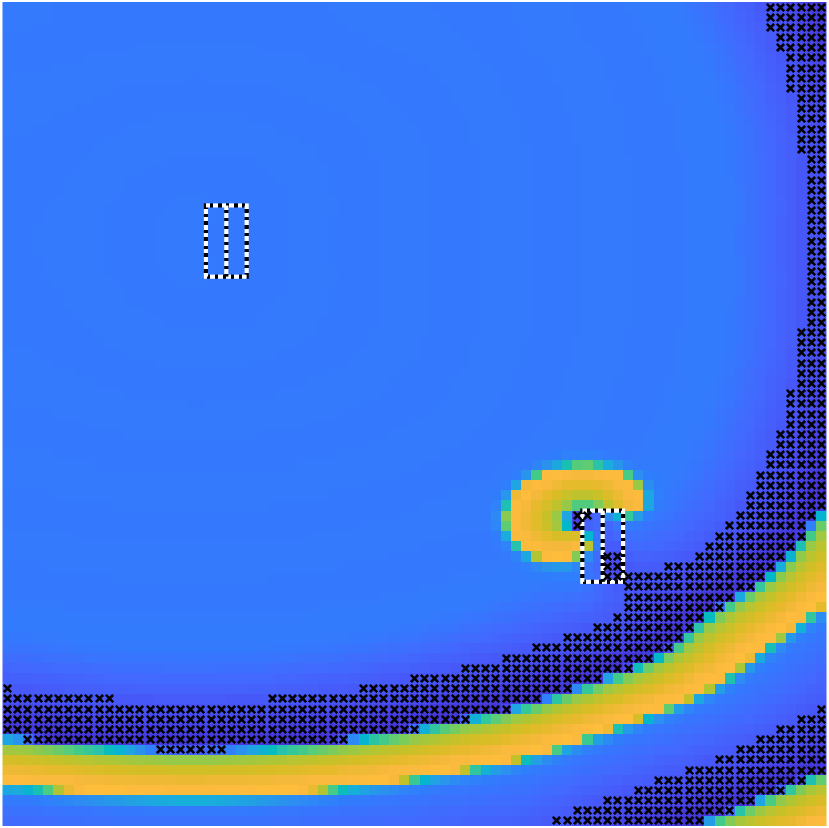}};
        \node[anchor=north] at (0,0) {$t=2125$};
        &
        \node[anchor=south] at (0,0) {\includegraphics[width=\figscale,height=\figscale]{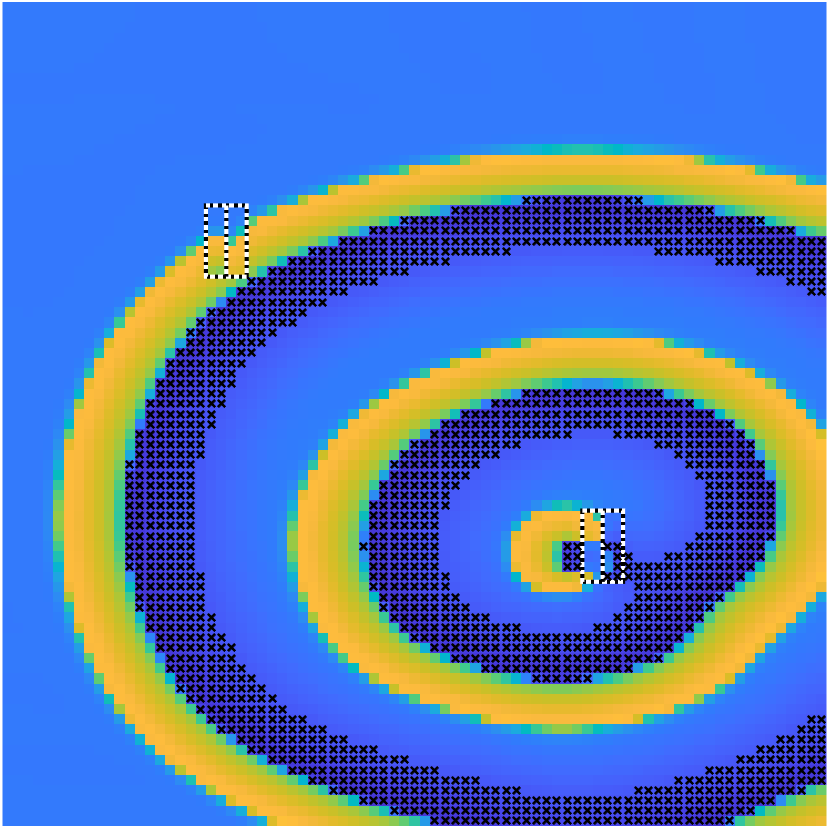}};
        \node[anchor=north] at (0,0) {$t=2600$};
        \\
        \node[anchor=south] at (0,0) {\includegraphics[width=\figscale,height=\figscale]{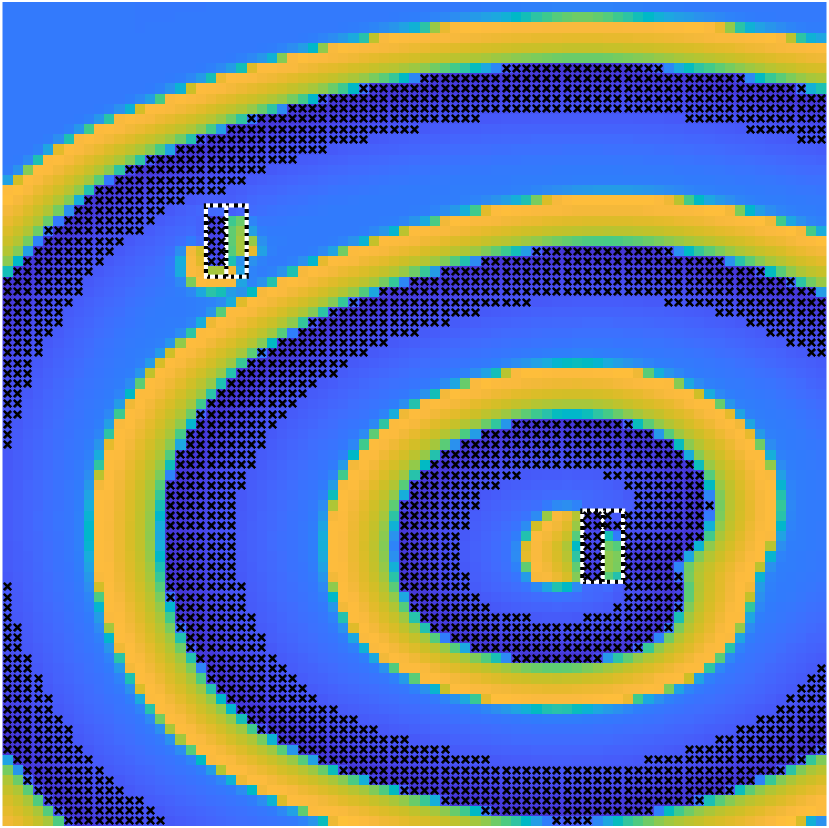}};
        \node[anchor=north] at (0,0) {$t=2800$};
        &
        \node[anchor=south] at (0,0) {\includegraphics[width=\figscale,height=\figscale]{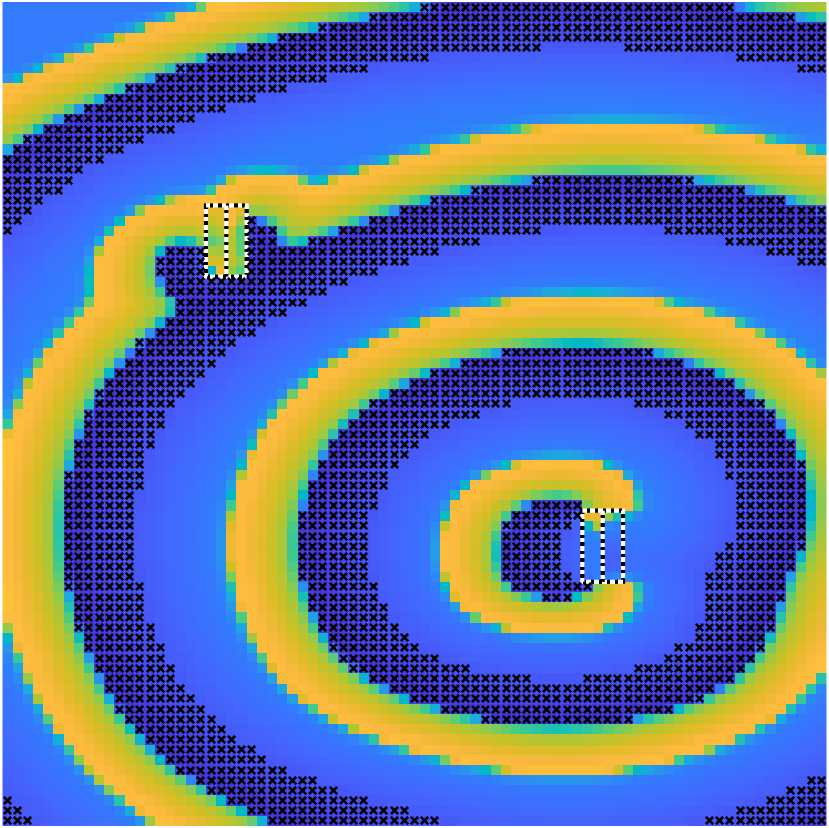}};
        \node[anchor=north] at (0,0) {$t=2900$};
        &
        \node[anchor=south] at (0,0) {\includegraphics[width=\figscale,height=\figscale]{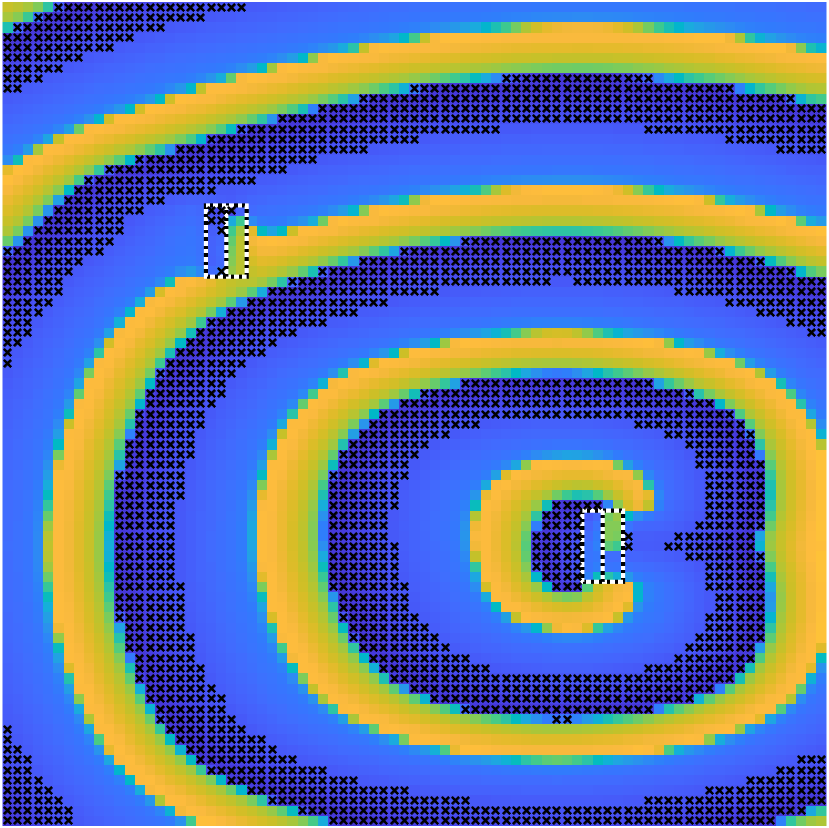}};
        \node[anchor=north] at (0,0) {$t=4000$};
        &
        \node[anchor=south] at (0,0) {\includegraphics[width=\figscale,height=\figscale]{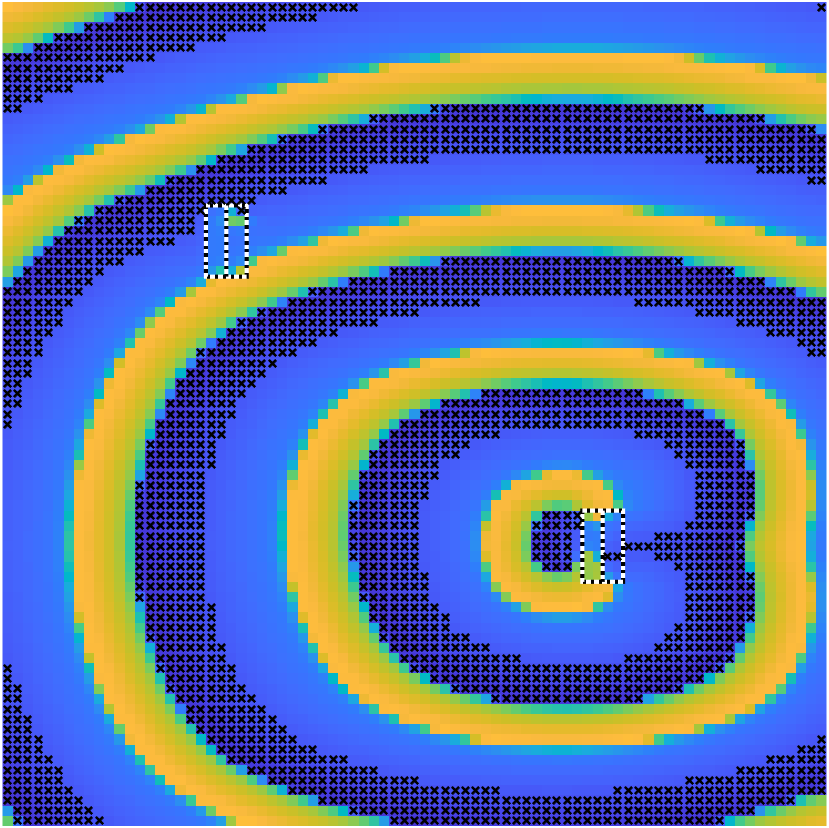}};
        \node[anchor=north] at (0,0) {$t=4200$};
        \\
        };
    \end{tikzpicture}
    \caption{Snapshots of a simulation where two $4 \times 7$ FIBs within a $81 \times 81$ grid interact. The upper-left FIB has slow-recovery $b=5.5$ and fast-recovery $\epsilon=0.007$, while the lower-right FIB has $b=5.3$ and $\epsilon = 0.006$. In the first row, we can observe colliding wave fronts between the FIBs. Next, we can see the dominant upper-left FIB ceasing to fire, and subsequently the lower-right FIB takes over and becomes the dominant FIB. 
    }
    \label{fig:two_fibs_interact}
\end{figure}
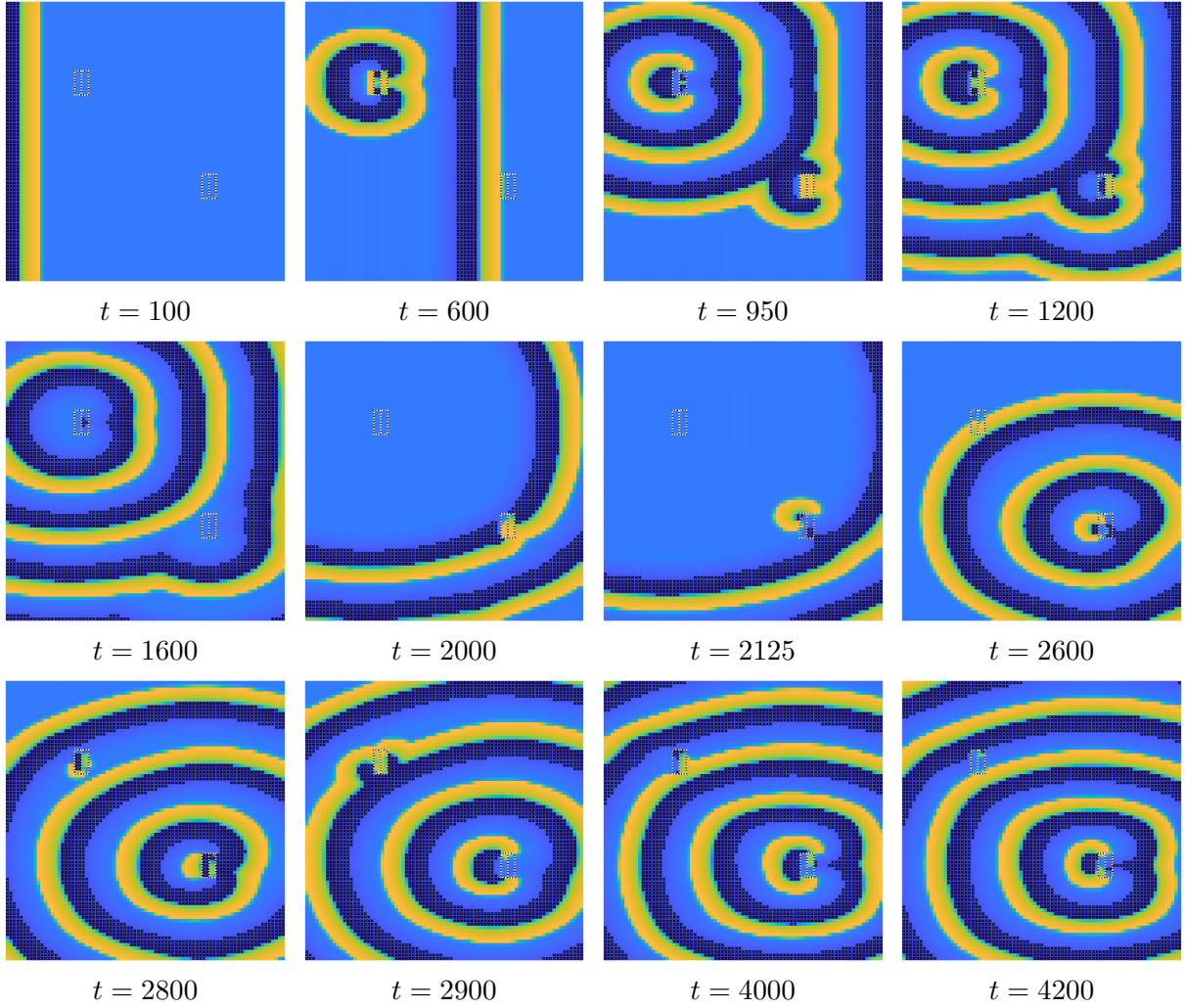

\section{Conclusion}
\label{sec:conclusion}

In this study, we proposed a simple but biologically relevant mechanism, FIBs, which are boundaries between fast-recovery and slow-recovery cells that can initiate fibrillatory-like behavior in two-dimensional networks of excitatory (FitzHugh-Nagumo) cells. 
In a network containing a FIB, fibrillation can arise spontaneously from a single excitation wave traveling normally through the tissue. 
Once initiated, the arrhythmia is self-sustained, at least for a moderate period. 
We provided a mechanistic illustration of the dynamics of the FIB and gave biologically-based explanations of how a FIB can form in heart tissue.
In detailed simulations, we studied the occurrence of fibrillation-like behavior as a function of parameters of underlying model and the geometry of the FIB.
We analyzed the relationship between charge retention time of the slow-recovery cells and time to recovery of the fast-recovery cells in initiating a backward wave.  
The simulations showed that for a small FIB, fibrillatory-like behavior occurs on a thin, fractal-like set of parameter values, but as the FIB grows the set becomes dense and contains large open regions.
The growth of the FIB is consistent with a typical progression in dynamics from periodic to intermittent to chaotic and is consistent with  the gradual onset of symptoms, from healthy to occasional arrhythmia to persistent atrial fibrillation.

Having established and analyzed the main properties of a FIB within an otherwise healthy network of cells, we studied several variations of the geometry of the FIB.  
We concluded first that the dynamics of a FIB are largely robust to changes in geometry as long as the FIB is large enough.  
We showed that spiral-like waves can be produced by setting specific initial conditions within healthy tissue, by the interaction of a FIB with scar tissue, and by differences in wave-propagation speed along a long FIB.
We showed that the center of the fibrillation can become detached from the FIB, giving the perception that the focal point of the arrhythmia is in an area of healthy cells. 
Finally, in simulations with multiple FIBs embedded in healthy tissue, we observed that one FIB can become dominant, thus masking the existence of the other FIB. 
All of these scenarios demonstrate that extreme care should be taken in attempting identify areas in a patient's atrium as candidates for ablation.

Our study reinforces the notion that small local heterogeneity in excitable tissue can initiate and sustain abnormal oscillations.  
We propose that heterogeneous structures similar to the FIBs introduced here are a potential cause of atrial fibrillation.
Our study also potentially illustrates why atrial fibrillation remains a difficult medical puzzle. 

\subsection*{Data Availability Statement}

The code used to generate \cref{fig:basicsnaps,fig:healthyspiral,fig:LFIB,fig:phantomfib,fig:phantomfib_paced,fig:vertical_strip_fib,fig:two_fibs_interact,fig:scar-line,fig:FIBlengths,fig:OSC 7 FIB,fig:pp_trace_3} is available at \cite{ohio-math-cardiac}.

\bibliography{AFbib}

\end{document}